\RequirePackage{fix-cm}
\documentclass[smallextended]{svjour3}       
\smartqed  
\usepackage{graphicx}
\usepackage{graphicx}
\usepackage{amssymb}
\usepackage{graphicx,fancyhdr}
\usepackage{multirow,amsmath}
\usepackage{dsfont}
\usepackage{subfigure}
\usepackage{bbm}
\usepackage{booktabs}
\usepackage{bbm}
\usepackage{amsfonts}
\usepackage{mathrsfs,float}
\usepackage{color,bm}

\begin{document}

\title{Fast Algorithms for Fourier extension  based on boundary interval data \thanks{The research is partially supported by National Natural Science Foundation of China (Nos. 12261131494,12171455, RSF-NSFC 23-41-00002)}
}


\author{Zhenyu Zhao  \and Yanfei Wang \and Anatoly G. Yagola 
}


\institute{Zhenyu Zhao \at
              School of Mathematics and Statistics, Shandong University of Technology, Zibo, 255049, China \\
                          \email{wozitianshanglai@163.com}           
           \and
          Yanfei Wang \at
             Key Laboratory of Deep Petroleum Intelligent Exploration and Development, Chinese Academy of Sciences, Beijing, 100029,China\\
              \email{yfwang@mail.iggcas.ac.cn}
              \and
              Anatoly G. Yagola \at
              Department of Mathematics, Faculty of Physics, Lomonosov Moscow State
University, Vorobyevy Gory, 119991 Moscow, Russia\\
  \email{yagola@physics.msu.ru}
}

\date{Received: date / Accepted: date}

\maketitle

\begin{abstract}{
This paper presents a novel boundary-optimized fast Fourier extension  algorithm for efficient approximation of non-periodic functions. The proposed methodology constructs periodic extensions through strategic utilization of boundary interval data, which is subsequently combined with original function samples to form an extended periodic representation. We develop a parameter optimization framework that preserves superalgebraic convergence while requiring only a few boundary node deployment, resulting in computational complexity marginally exceeding that of standard FFT implementations. Furthermore, we present an improved version of the algorithm tailored for functions exhibiting boundary oscillations. This variant employs grid refinement near the boundaries, which reduces the resolution constant to approximately one-fourth of that in conventional approaches. Comprehensive numerical experiments confirm the efficiency and accuracy of the proposed methods and establish practical guidelines for parameter selection.}
\end{abstract}

\section{Introduction}

Fourier series approximation has been widely used in various  fields for smooth periodic functions. In many physical, engineering and exploration geophysical problems, such as analyzing periodic signals in electrical engineering, studying periodic motion in mechanics such as the oscillation of a pendulum under certain ideal conditions, or decomposing seismic/electromagnetic signals into different frequency components using Fourier series approximation \cite{SteinFA2003,Tikhonov1995,WangIP2011}. The Fourier series approximation offers spectral convergence since more terms are added to the Fourier series, the approximation converges to the actual function in an efficient way in the frequency domain. Moreover, it can be computed numerically via the FFT (Fast Fourier Transform), which is a computationally efficient algorithm. However, the situation changes completely once the periodicity is lost. A jump discontinuity is imposed in the approximation at the domain boundaries, owing to the well-known Gibbs phenomenon. For instance, when we attempt to approximate a  function that is not periodic over a given interval using Fourier series, we will observe this Gibbs phenomenon near the boundaries of the interval. Several methods have been developed to ameliorate the ill effects of the Gibbs phenomenon. For example, methods based on Gegenbauer polynomials, which have certain orthogonality properties that can be exploited to obtain better approximate functions \cite{2006Robust,1992On,2006Optimal}. Methods utilizing Pad\'{e} approximations, which are rational function approximations, can also be used \cite{2000A,Geer1995}. And Fourier extension (or continuation) methods, which extend the function in a certain way to make it more suitable for Fourier-based approximation, are also among the approaches developed to address this problem \cite{Boyd2002A,Bruno2022,Huybrechs2010,Lyon2011,2012Sobolev,Matthysen2016FAST,Matthysen2018,Webb2020}.

The idea of the existing discrete Fourier extension  can be described as follows: Let $f$ be a non-periodic, smooth and sufficiently differentiable  function over $[-1,1]$,  which is sampled in a predefined set of points $\{t_{\ell}\}\subset[-1,1]$, the discrete Fourier extension $C(f)$ of $f$ is defined by
 \begin{equation}
   C(f)=\sum_{k=-N}^N{\bf c}_ke^{\text{i}\frac{k\pi}{T}t},\quad T>1,
 \end{equation}
where the coefficients $\{{\bf c}_k\}_{k=-N}^N$ are chosen to minimize the error
\begin{equation}\label{minsol}
  \min\sum_{\ell}\left|f(t_{\ell})-\sum_{k=-N}^N{\bf c}_ke^{\text{i}\frac{k\pi}{T}t_{\ell}}\right|^2.
\end{equation}

    To obtain a solution to the minimization problem \eqref{minsol},  a severely ill-conditioned least-squares system has to be solved. In \cite{Boyd2002A}, Boyd   introduced   truncated singular value decomposition (TSVD) to  stabilize the calculation process, with an accuracy almost up to machine precision. A convergence analysis by Adcock et al. states that for analytic functions the method is at least superalgebraically convergent \cite{adcock2014numerical}. Utilizing randomized algorithms, Lyon   provided a fast algorithm  with $T=2$ \cite{Lyon2011}. In \cite{Matthysen2016FAST}, Matthysen et al. presented  fast algorithms for arbitrary $T$, exploiting it as a specific variant of bandlimited extrapolation. This idea was further developed in \cite{Matthysen2018} to form a fast algorithm for function approximation in 2D domains.

\begin{figure}
	\begin{center}
		{\resizebox*{7cm}{5cm}{\includegraphics{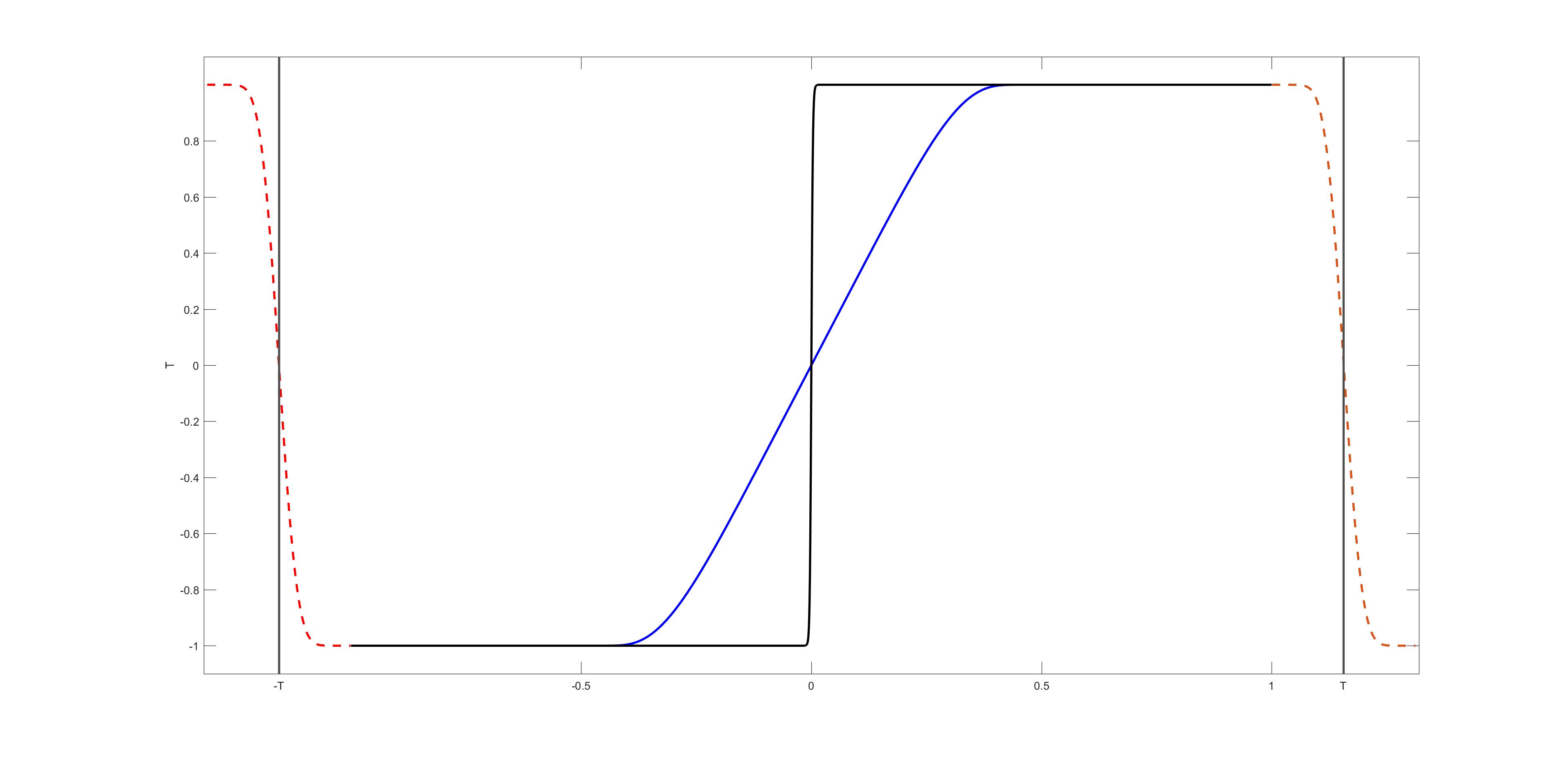}}}
		{\caption{Periodic extensions of $f_1(t)$ and $f_2(t)$. \label{fig1}}}
	\end{center}
\end{figure}

    It should be noted that the goal of the extension process is to effectively connect the data at both ends of the function. The following two functions were constructed:
    \begin{equation}
    \begin{aligned}
      f_1(t)=\left\{
      \begin{aligned}
        &-1,&x<-\frac{1}{2},\\
        &\tanh\left(\tan \pi t\right),&-\frac{1}{2}\leq t\leq \frac{1}{2},\\
        &1,&t>\frac{1}{2},
      \end{aligned}\right.\\
       f_2(t)=\left\{
      \begin{aligned}
        &-1,&t<-\frac{1}{2},\\
        &\tanh\left(100\tan \pi t\right),&-\frac{1}{2}\leq t\leq \frac{1}{2},\\
        &1,&t>\frac{1}{2}.
      \end{aligned}\right.
    \end{aligned}
          \end{equation}
Because the data of the two functions near the endpoints are the same, the same connection data can be used in the intervals $[-T, 1)$ and $(1, T]$. An illustration of this process is shown in Fig. \ref{fig1}. This raises the issue of whether the data inside the interval really help in the calculation of the extension. In fact, a very efficient method, known as the FC-Gram method has been developed earlier, which uses only boundary data (projecting it into a basis of orthogonal polynomials) for extension calculations \cite{bruno2010high}. Since then, a lot of research has been done on further improvement and application of the FC-Gram algorithm \cite{2016An,2020A,Bruno2022}. In \cite{plonka2018numerical}, the authors proposed a method to calculate the Fourier extension using the derivative on the boundary combined with a two-point Taylor interpolation polynomial. The extension part calculated by these two methods are more stable, but neither of them is superalgebraically convergent. When the boundary interval is more volatile than the interior, the number of nodes required to achieve near-machine precision by these two methods was  higher than that of the previous method.  { In \cite{bruno2010high}, the authors identified a critical limitation in directly applying the original extension algorithm with relatively few  nodes in the boundary interval: the inherent difficulty in achieving satisfactory approximation accuracy.  { To address this, they introduced a pre-processing step involving Gram polynomial basis generation. Existing theoretical analyses reveal that Fourier extension accuracy is co-determined by two interdependent parameters - extension length and sampling ratio. While previous studies prioritized increasing sampling ratios for accuracy improvement (as extension length adjustments offered limited benefits in full-data algorithms), this work demonstrates that strategically increasing extension length can effectively reduce boundary node requirements. Specifically, we establish that through optimal parameter combinations, accurate extension computation can be achieved with a few boundary node deployment, where the required node count becomes fixed and solely dependent on target accuracy.

This paper presents a novel fast Fourier extension algorithm with the following key innovations:}
\begin{itemize}
  \item{ Parameter optimization enables extension computation using minimal boundary nodes while maintaining FFT-compatible efficiency. Unlike FC-Gram approaches, our method eliminates polynomial basis requirements.}
  \item Because the calculation of the extension part  uses only the data of the small interval near the boundary, the influence of oscillations within the interval on the approximation is significantly weakened. Therefore, for a function that is smooth in the neighborhood of the interval boundary and oscillates within the interval, the proposed algorithm has a smaller resolution constant than that of the original full-data algorithm. This also applies to the FC-Gram algorithm.
      \end{itemize}}
Furthermore, for functions that oscillate in the neighborhood of the interval boundary, we propose an improved algorithm based on boundary interval grid refinement, which can significantly reduce the resolution constant of this type of function. The proposed algorithms exhibit good performance and provide  clear conclusions regarding the settings of various parameters. Therefore, we believe that it is of substantial  interest to those who use FEs  in practice.

\subsection{Overview of the paper} In Section \ref{SEC2}, we outline the algorithm and explain its structure. The specific calculation method of the extension part based on boundary interval data is presented in Section \ref{SEC3}. Section \ref{SEC4} present the computational complexity and resolution constant analysis of the algorithm. In Section \ref{SEC5}, we  systematically tests the key parameters, identifies their influencing factors, and  provide an optimization setting scheme. In Section \ref{SEC6}, we present  numerical tests of the final algorithm and compare it with three other fast algorithms. Furthermore, for the boundary oscillation function, we propose an improved algorithm based on the boundary interval grid refinement.

\section{Outline of the boundary interval algorithm\label{SEC2}}
In this paper, we restrict ourselves to the case where the sampling points are uniform, that is,
\begin{equation}
  t_{\ell}=\frac{\ell}{M},\quad \ell=-M,\ldots,M.
\end{equation}
{ Let $m_{\Delta}\in  \mathbb{N}^+$ denote the number of boundary interval nodes, with corresponding index sets
\begin{equation}\begin{aligned}
  &S_l=\{-M,-M+1,\ldots,-M+m_{\Delta}-1\},\\
  &S_r=\{M-m_{\Delta}+1,M-m_{\Delta}+2,\ldots,M\}.
\end{aligned}
\end{equation}
The calculation of the extension part  is performed in the interval $[0,2\pi]$. Define $T_{\Delta}>1$ as the scaling ratio between the extended interval $[0, 2\pi]$ and the data-containing interval, analogous to parameter $T$ in the original full-data algorithm.  Further determine the number of nodes $L_{\Delta}$ used in the extension calculation and the corresponding discrete node $\{x_j\}$, specifically}
\begin{equation}
  \label{defLd}
  \begin{aligned}
&L_{\Delta}= 2\times \left\lceil T_{\Delta}\times (m_{\Delta}-1)\right\rceil,\\&
 h=\frac{2\pi}{L_{\Delta}}, \quad x_j=(j-1)h,\quad j=1, 2, \ldots, L_{\Delta},
 \end{aligned}
\end{equation}
where $\lceil\cdot\rceil$ is the rounding symbol.
 Then let
\begin{equation}
\begin{aligned}
&  J_{\Delta,1}=\{1,2,\ldots,m_{\Delta}\},~ J_{\Delta,2}=\{\frac{L_{\Delta}}{2}+1,\frac{L_{\Delta}}{2}+2,\ldots,\frac{L_{\Delta}}{2}+m_{\Delta}\}, \\
& J_{\Delta}=J_{\Delta,1}\cup J_{\Delta,2}.
\end{aligned}
\end{equation}
{ They indicate the locations of the boundary data during the extension calculation process.}
For $j\in J_{\Delta}$, take
\begin{equation}\label{defg}
  g(x_j)=\left\{
  \begin{aligned}
    &f(t_{M-m_{\Delta}+j}),& j\in J_{\Delta,1},\\
    &f(t_{-M+j-{L_{\Delta}}/{2}-1}),&j\in J_{\Delta,2}.
  \end{aligned}\right.
\end{equation}
Then we calculate the extension function $g_c$  defined on the interval $[0,2\pi]$ of $g$, and obtain its values at $x_j, j=1,2,\ldots, L_{\Delta}$  according to the method described in Section \ref{SEC3}.
Finally, we   obtain the function $f_c$ with a period of $2+\lambda$, where
 \begin{equation}\label{deflambda}
   \lambda=\frac{\lceil T_{\Delta}-1\rceil\times(m_{\Delta}-1)}{M},
 \end{equation}
 and the data of the extension function $f_c(t_{\ell}),\ell=-M,\ldots,M, M+1,\ldots, M+\frac{L_{\Delta}}{2}-m_{\Delta}$ in one period $[-1,1+\lambda)$ is
\begin{equation}
  f_c(t_{\ell})=\left\{
  \begin{aligned}
    &f(t_{\ell}),&|\ell|\leq M,\\
&g_c(x_{m_{\Delta}+\ell-M}),&\ell>M.
  \end{aligned}\right.
\end{equation}

{
\begin{figure}
			\begin{center}
\subfigure[\label{1a}Given data and selected boundary data] {
\resizebox*{5.5cm}{!}{\includegraphics{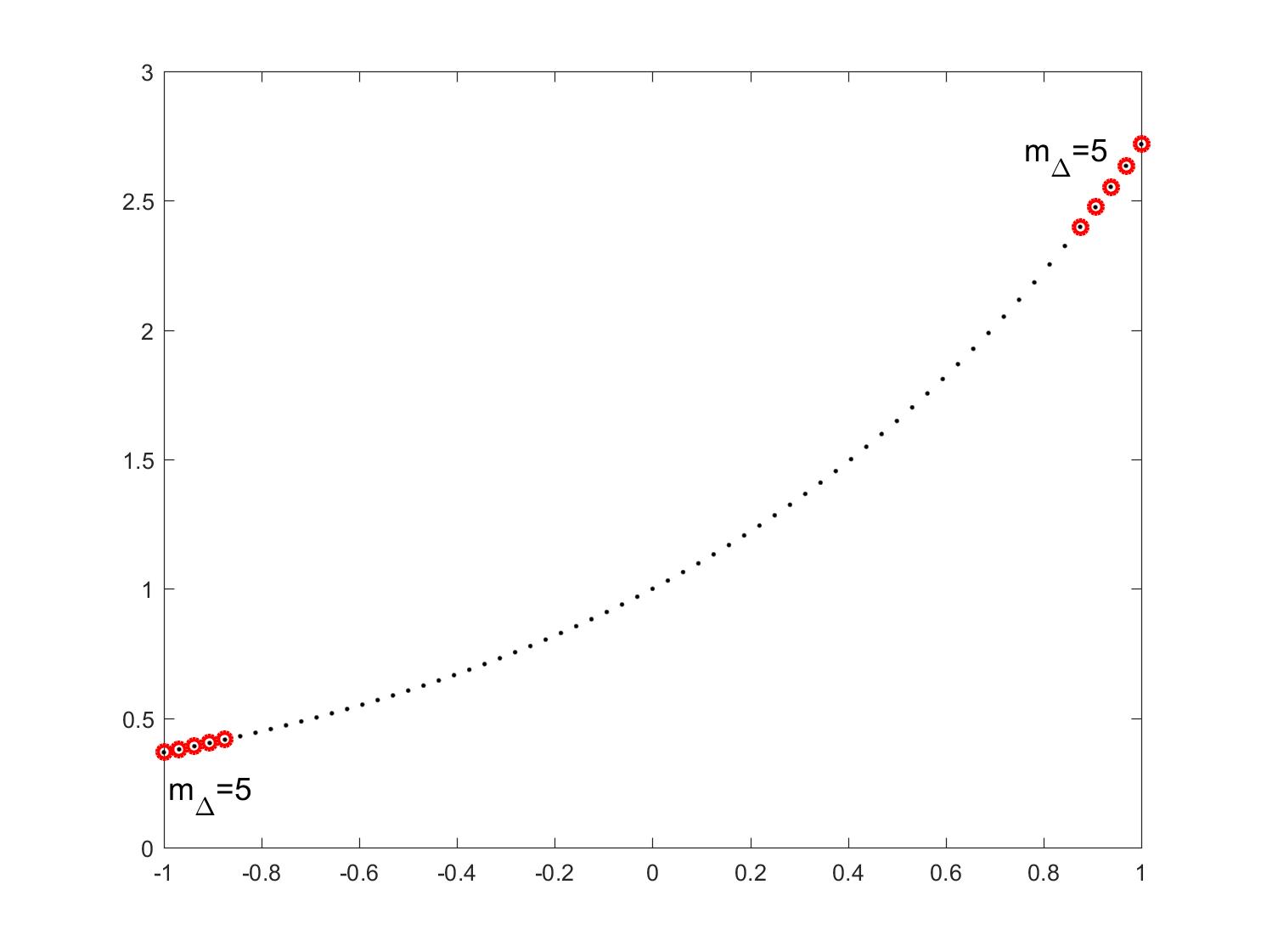}}
}%
\subfigure[\label{1b}$g_c(x_j)$] {
\resizebox*{5.5cm}{!}{\includegraphics{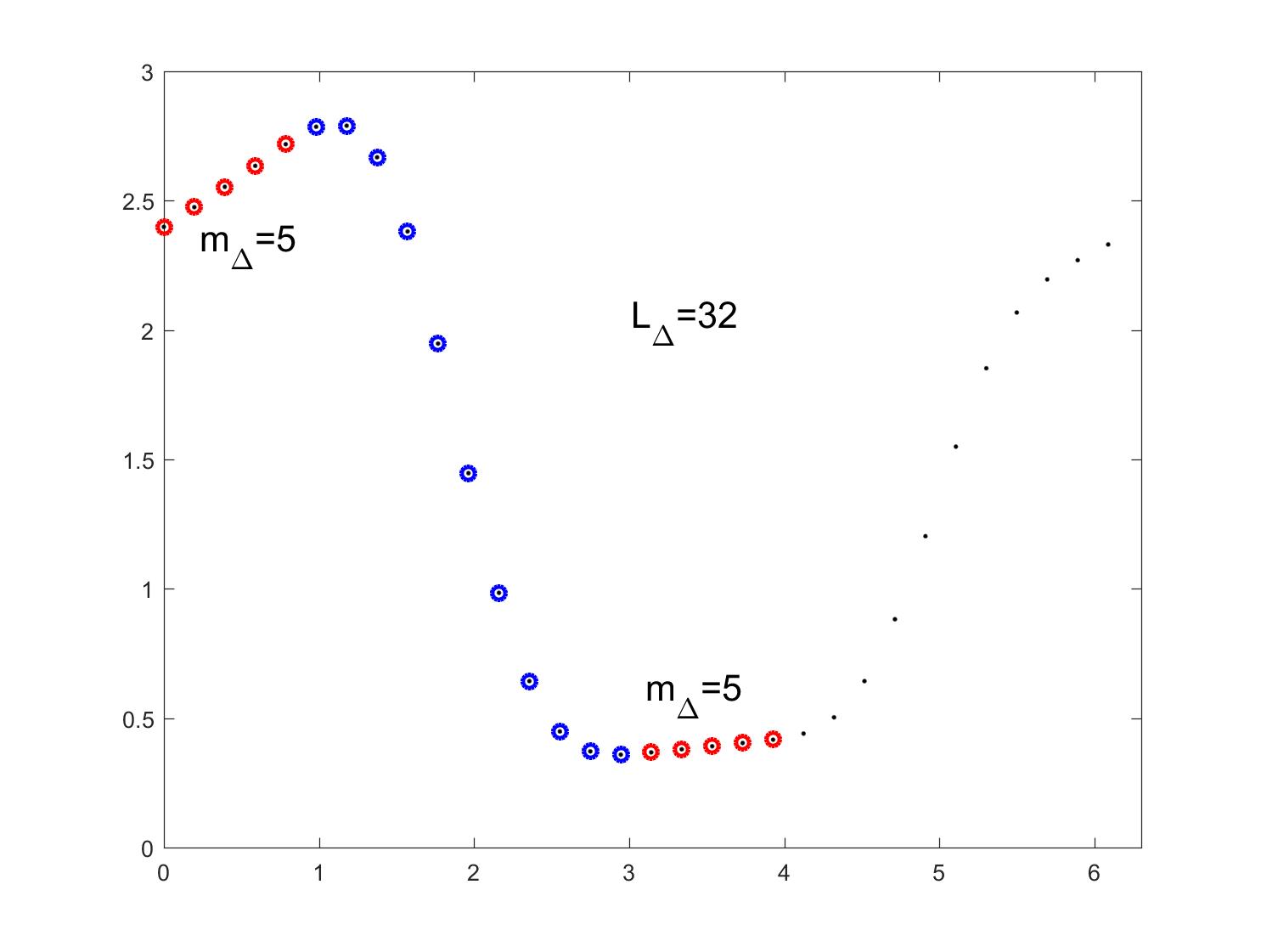}}
}%

 \subfigure[\label{1c}$f_c(t_{\ell})$] {
\resizebox*{5.5cm}{!}{\includegraphics{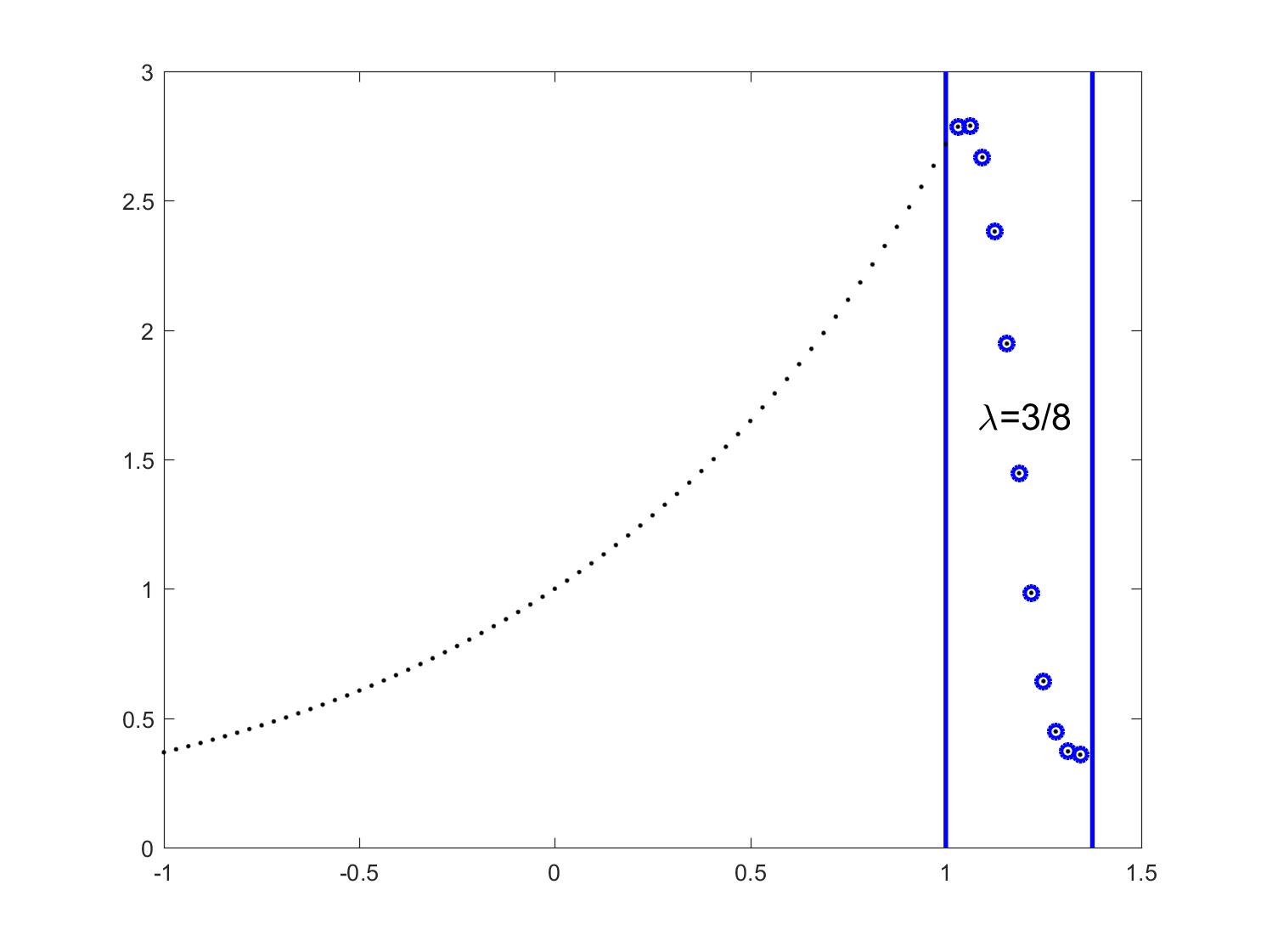}}
}%
		{\caption{Calculation process of a periodic extension of $f(x)=\exp(x)$ ($T_{\Delta}=4, M=32$). \label{fig2}}}
	\end{center}
\end{figure}}

Fig. \ref{fig2} shows the basic process of the boundary interval algorithm  using the function {$f(x)=\exp(x)$ }as an example. Here,  to make the graph  clearer, we use a relatively small number of nodes. Fig. \ref{1a} shows the given data (indicated by the black dotted line) and   boundary data selected for the calculation extension (indicated by the red circle). Fig. \ref{1b} shows the $g_c(x_j)$ calculated from the boundary data. The data represented by   blue circles represent the subsequent extension data required. Fig. \ref{1c} shows the data of $f_c$ in one period formed by combining the original data with the calculated extension data. It can be observed from the calculation process that the calculation of the extension part  depends  only on the boundary interval data.

After obtaining the data in one period, the subsequent algorithm can be implemented using FFT. Therefore, whether the algorithm can achieve satisfactory efficiency depends on the computational complexity of the extension part, that is, the size of the parameters $T_{\Delta}$ and $m_{\Delta}$. We present their selection schemes through testing in Section \ref{SEC5}.

\section{ The calculation of $g_c$\label{SEC3}}
Including the appropriate normalization, we let
\begin{equation}
  \phi_k(x)=\frac{1}{\sqrt{L_{\Delta}}}e^{\text{i}kx}, \quad k=-n_{\Delta},\ldots,n_{\Delta},\quad x\in[0,2\pi].
\end{equation}
Assuming  $g(x_j)$ is given by \eqref{defg}, then the  extension function $g_c$ of $g$  is defined as
\begin{equation}
  g_c(x)=\sum_{k=-n_{\Delta}}^{n_{\Delta}} {\bf c}_k\phi_k(x).
\end{equation}
where the coefficients ${\bf c}_k$ are determined by
\begin{equation}\label{lqfb}
  \min_{({\bf c}_k)}\sum_{j\in J_{\Delta}}\left|\sum_{k=-n_{\Delta}}^{n_{\Delta}}{\bf c}_k\phi(x_j)-g(x_j)\right|^2.
\end{equation}
If we let
\begin{equation}
  \begin{aligned}
   &A_{l,k}&=&\left\{
  \begin{aligned}
    &\phi_k(x_l),& 1\leq l\leq m_{\Delta},\\
    &\phi_k(x_{\frac{L_{\Delta}}{2}+l-m_{\Delta}}),& m_{\Delta}<l\leq 2m_{\Delta},
  \end{aligned}
  \right.\\
 & b_l&=&\left\{
  \begin{aligned}
    &\frac{1}{m_{\Delta}}f(x_l),& 1\leq l\leq m_{\Delta},\\
    &\frac{1}{m_{\Delta}}f(x_{\frac{L_{\Delta}}{2}+l-m_{\Delta}}),& m_{\Delta}<l\leq 2m_{\Delta},
  \end{aligned}
  \right.\\
  &&&l=1, 2, \ldots, 2m_{\Delta},k=-n_{\Delta}, \ldots, n_{\Delta}.
  \end{aligned}
\end{equation}
\begin{figure}
			\begin{center}
\subfigure[\label{3a}$n_{\Delta}=100$, $T_{\Delta}=6$] {
\resizebox*{6cm}{!}{\includegraphics{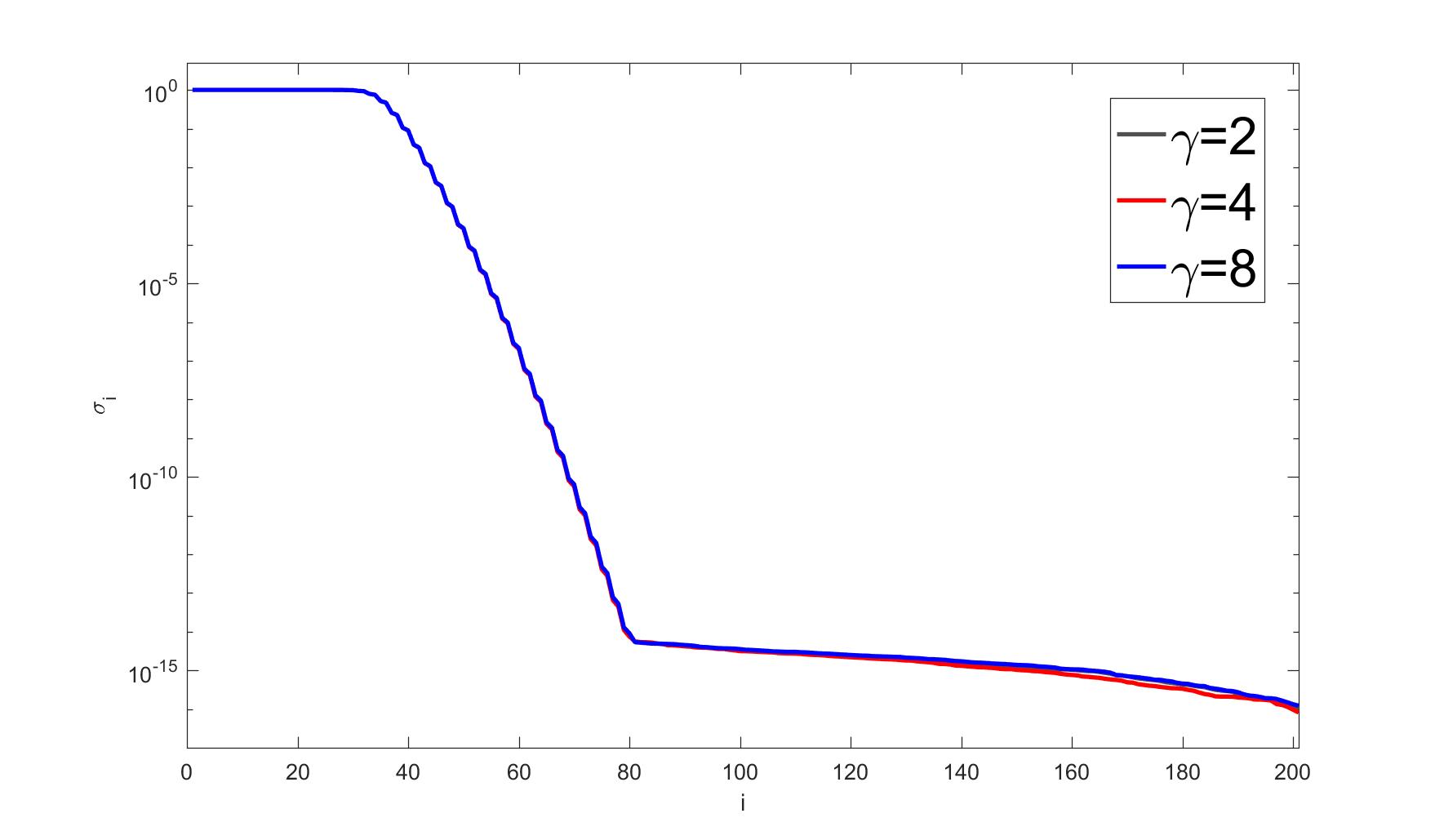}}
}%
\subfigure[\label{3b}$m_{\Delta}=100$, $T_{\Delta}=20$] {
\resizebox*{6cm}{!}{\includegraphics{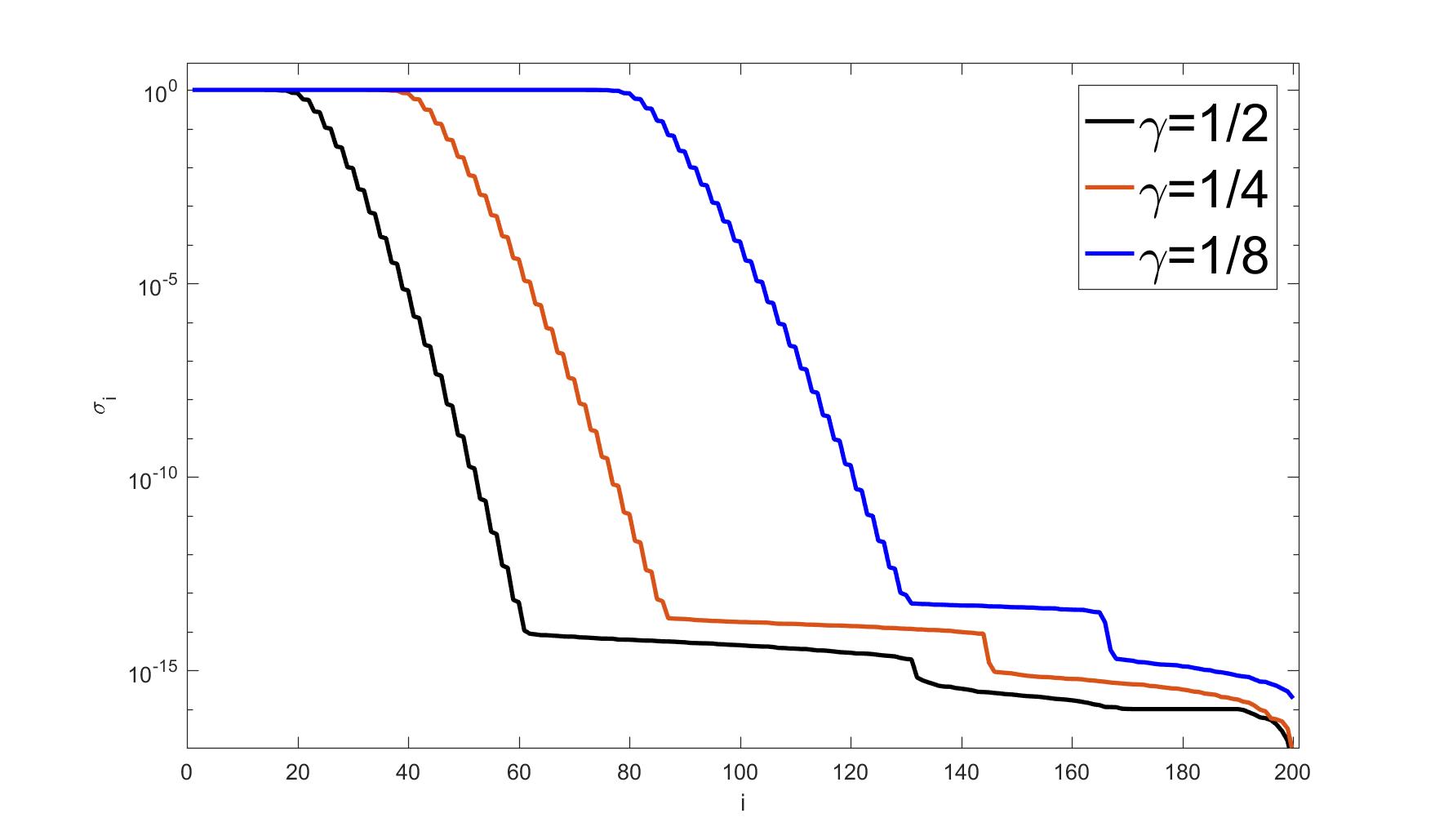}}
}%

\subfigure[\label{3c}$m_{\Delta}=100$, $n_{\Delta}=100$] {
\resizebox*{6cm}{!}{\includegraphics{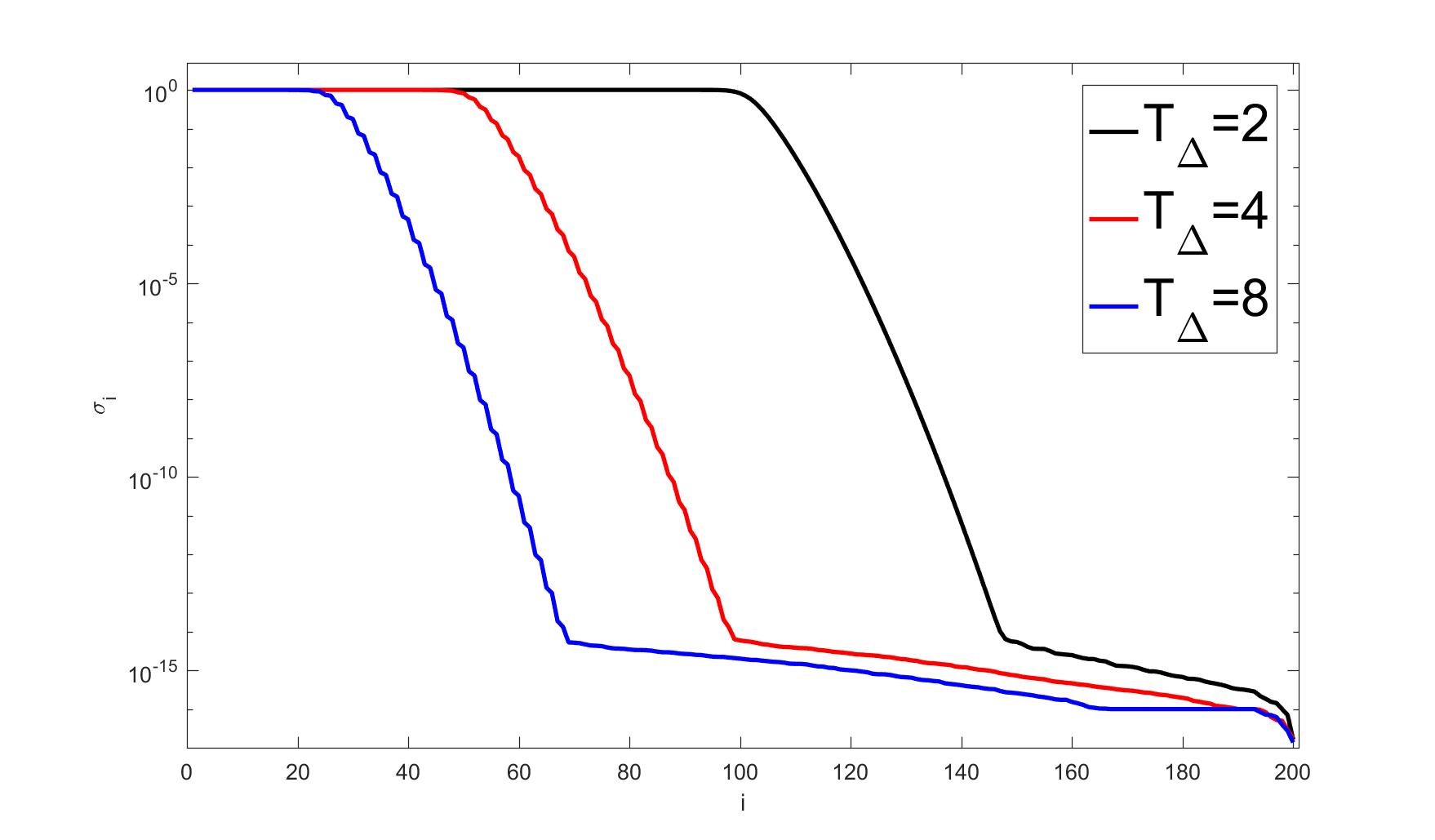}}
}%

		{\caption{The singular values of $A$ for various $m_{\Delta}$, $n_{\Delta}$ and $T_{\Delta}$. \label{fig3}}}
	\end{center}
\end{figure}

Then the solution of \eqref{lqfb} can be obtained  by solving the rectangular system
\begin{equation}\label{lseq}
  A{\bf c}\approx b
\end{equation}
in a least squares sense. If we let
\begin{equation}
  N_{\Delta}=\min\{2m_{\Delta},2n_{\Delta}+1\},
\end{equation} then  the SVD of $A$ can be given as
\begin{equation}
  A=\sum_{i=1}^{N_{\Delta}}{\bf u}_i\sigma_i{\bf v}^T_i,
\end{equation}
where $\sigma_i$ are the singular values of $A$ such that
\begin{equation}
  \sigma_1\geq \ldots\geq \sigma_{N_{\Delta}}\geq 0,
\end{equation}
while ${\bf u}_i$ and ${\bf v}_i$ are the left and right singular vectors of $A$, respectively. We define
\begin{equation}
  \label{defgamma}
  \gamma_{\Delta}=\frac{m_{\Delta}-1}{n_{\Delta}}.
\end{equation}
Fig. \ref{fig3}  shows the singular values of $A$ calculated by  SVD techniques for various $m_{\Delta}$, $n_{\Delta}$ and $T_{\Delta}$ in MATLAB. It can be observed that the first approximately $\mathcal{O}(n_{\Delta}/T_{\Delta})$ are close to $1$ and then the singular values
rapidly decay until the true singular values are masked by machine precision errors.

The situation of singular values reflects the ill-posedness of \eqref{lseq}, regularization must be introduced. Referring to the treatment methods in \cite{adcock2014numerical}, the truncated SVD least-squares solution to  \eqref{lseq} can be given as
\begin{equation}
  {\bf c}_{\tau}=\sum_{\sigma_i>\tau} \frac{{\bf u}_i^Tb}{\sigma_i}{\bf v}_i,
\end{equation}
where   parameter $\tau$ is usually chosen to be close to the machine precision.
{ \section{Computational complexity  and resolution constant \label{SEC4}}
The computational complexity of the proposed algorithm decomposes into three principal components:
\begin{itemize}
  \item {\bf Calculation of $c_\tau$}. The computational effort in this step is mainly concentrated on the singular value decomposition of   matrix $A$, which is related to $m_{\Delta}$ and $\gamma_{\Delta}$, and can be estimated as
\begin{itemize}
  \item If $\gamma _{\Delta}< 1$, the computational complexity is
\begin{equation}
  \mathcal{O}\left(m_{\Delta}^2n_{\Delta}\right)=\mathcal{O}\left(\frac{1}{\gamma_{\Delta}}m_{\Delta}^3\right).
\end{equation}
  \item If $\gamma_{\Delta} \geq 1$, the computational complexity is
\begin{equation}
  \mathcal{O}\left(m_{\Delta}n^2_{\Delta}\right)=\mathcal{O}\left(\frac{1}{\gamma_{\Delta}^2}m_{\Delta}^3\right).
\end{equation}
\end{itemize}
  \item {\bf Calculation of $g_c$}. This step can be achieved by IFFT, the computational complexity is
\begin{equation}
  O\left(L_{\Delta}\log(L_{\Delta})\right).
\end{equation}
  \item {\bf Calculation of Fourier coefficients of $f_c$.} This can be achieved by FFT, and the computational complexity is
\begin{equation}
  O\left(M\log(M)\right).
\end{equation}
\end{itemize}
Next, we analyze the factors that affect these parameters in detail. First, we provide the following concept for the  resolution constant, which comes from the literature  \cite{adcock2014parameter}:
\begin{definition}
  Let $\left\{F_M\right\}_{M\in\mathbb{N}}$ be a sequence of approximations such that $F_M(f)$ depends
only on the values of $f$ on an equispaced grid of $2M + 1$ points. For $\omega>0,~~0<\delta<1,$ let
$$\mathcal{R}(\omega,\delta)=\min\left\{M\in\mathbb{N}:\|\exp(\text{i}\pi\omega t)-F_M(\exp(\text{i}\pi\omega t))\|_{\infty}\leq\delta\right\},$$
then we say that $F_M$ has resolution constant $0 < r < \infty$ if
\begin{equation}
 R(\omega, \delta) \sim r\omega, \quad \omega\rightarrow\infty,
\end{equation}
for any fixed $\delta$.
\end{definition}

The resolution constant of the present method  is influenced  by two aspects:
\begin{itemize}
  \item First, we need to consider the approximation ability of $g_c$ to the boundary interval data of the function $f$. For the function $f(x)=\exp(\text{i}\pi\omega t)$ on $[-1,1]$, since the calculation of $g_c$ is performed on $[0,2\pi]$ after intercepting the boundary interval data, the function corresponding to $g$ in \eqref{defg} is $$g(x)=\exp(\text{i}\omega_{\Delta}\pi t),$$
      where
\begin{equation}\label{estrs1}
  \omega_{\Delta}=\frac{m_{\Delta}-1}{M}T_{\Delta}\omega.
\end{equation}
Obviously, if $g_c$ is to be able to effectively approximate $g$, it must satisfy
\begin{equation}\label{estrs2}
  n_{\Delta}\geq \omega_{\Delta}.
\end{equation}
Therefore, according to \eqref{estrs1}, \eqref{estrs2} and \eqref{defgamma}, it can be obtained that the number of nodes $M$ have to satisfy
\begin{equation}\label{estrrs3}
  M\geq  \gamma_{\Delta}T_{\Delta}\omega.
\end{equation}
  \item Second, it is necessary to consider the number of nodes required for using $f_c$ to calculate the overall approximation. In this step we use FFT to approximate a periodic function, according to the basic conclusion of Fourier approximation, this step requires the number of nodes $M$ to satisfy
      \begin{equation}\label{estrrs4}
        M\geq \omega.
      \end{equation}
\end{itemize}
Note that in the above analysis process, we consider the case where the internal and boundary frequencies are the same. In the actual calculation process, when facing general functions, it can be seen from the above analysis  that if the frequency of the function at the boundary is significantly lower than that of the interior, the number of nodes required for calculation will be greatly reduced. This also applies to the FC-Gram algorithm. We will verify this in the subsequent numerical experiments.

Through (4.1)-(4.3) and the above analysis, we can see that whether the proposed method can have a smaller computational complexity depends on whether there is a smaller $m_{\Delta}$ value that can enable $g_c$ to accurately approximate the function in the boundary interval. In the original full-data Fourier extension algorithm, most literatures recommend the use of oversampling mode. The results show that the accuracy of the calculation can be improved with the increase of the sampling ratio $$\gamma=\frac{M}{N}.$$
In fact, the increase of $\gamma$ plays a role in reducing the discrete step length on $[-T,T]$, which can also be achieved by increasing $T$, as mentioned in the previous literature on the estimation of the resolution constant of the full-data algorithm \cite{adcock2014parameter}:
\begin{equation}\label{rcfulldata}
  r\sim T\gamma.
\end{equation}
In previous studies on full-data algorithms, researcher tend to increase $\gamma$ because no benefits were observed with increasing $T$.  However, for the algorithm proposed in this paper, increasing $T$ seems to be a better choice because it  significantly reduces the number of nodes required for the boundary interval. In the next section, we will test these parameters more comprehensively and provide clear results and parameter setting suggestions, so that the algorithm can complete the calculation with lower computational complexity.}

\section{Numerical  investigation and determination of key parameters\label{SEC5}}
The proposed algorithm involves multiple parameters, and their selection has a substantial influence on the performance of the algorithm. This section   tests and analyzes them and then provides an optimized  scheme.

\subsection{Testing the parameter $T_{\Delta}$}

 From the analysis in the previous section, we know that the extension length $T_{\Delta}$ has a substantial influence on the convergence rate and the resolution power. Generally, a larger extension length can achieve a higher convergence rate but also means a larger resolution constant.

\begin{figure}
			\begin{center}
\subfigure[\label{4a} $\gamma_{\Delta}=1$, $M=500$, $m_{\Delta}=100$.] {
\resizebox*{6cm}{!}{\includegraphics{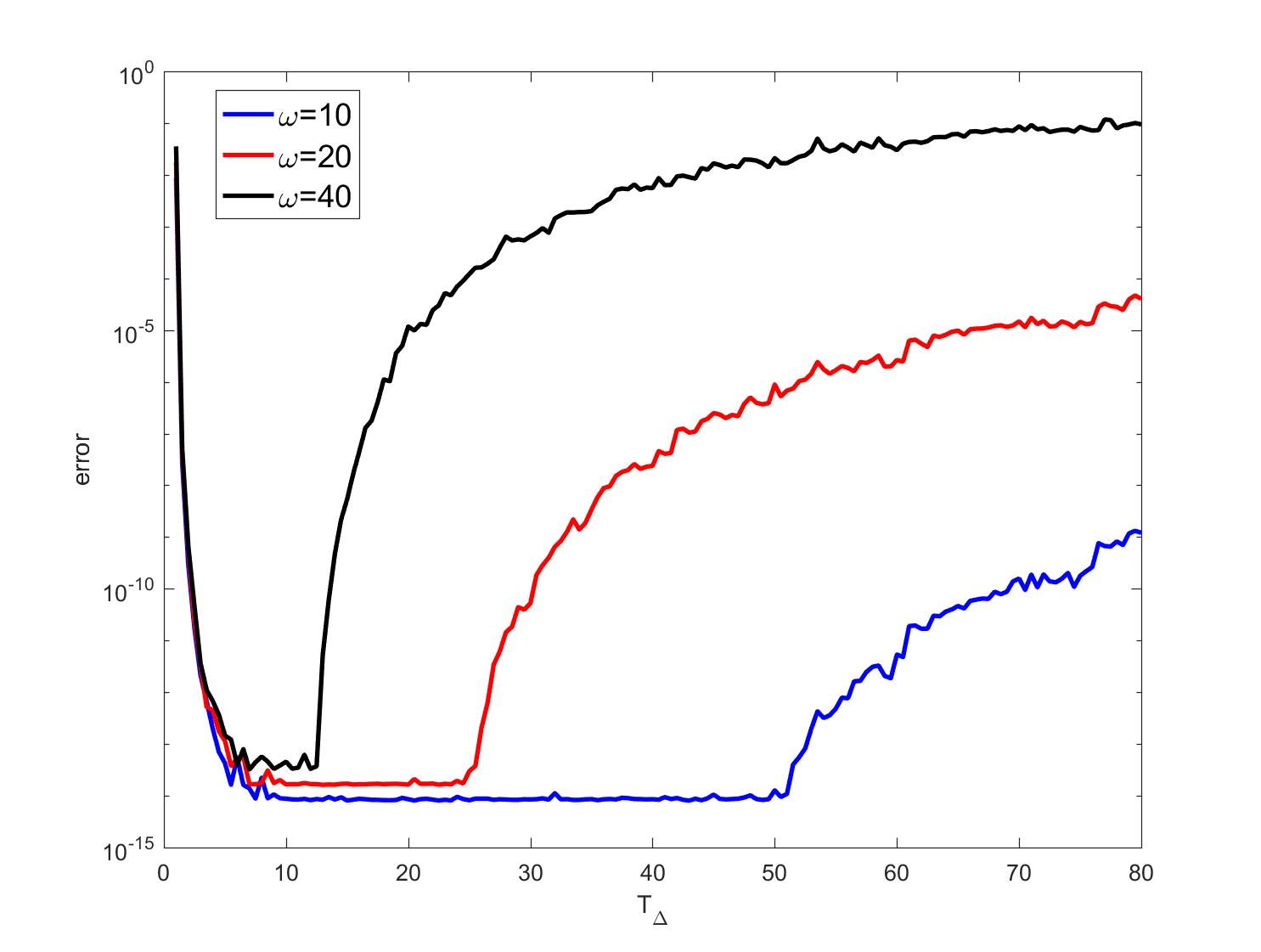}}
}%
\subfigure[\label{4b} $\gamma_{\Delta}=1$, $M=500$, $\omega=20$.] {
\resizebox*{6cm}{!}{\includegraphics{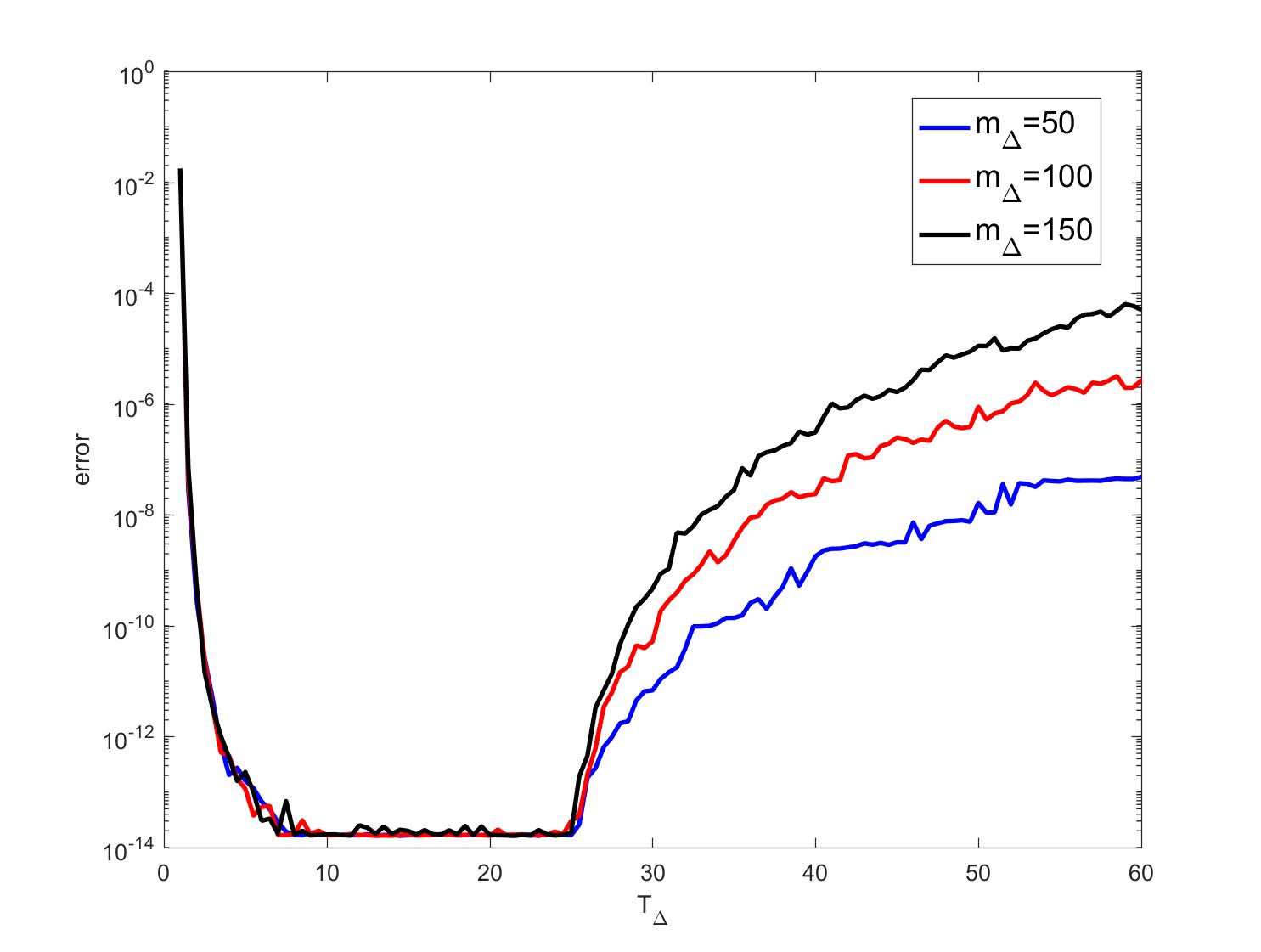}}
}%

\subfigure[\label{4c} $M=500$, $m_{\Delta}=100$, $\omega=20$.] {
\resizebox*{6cm}{!}{\includegraphics{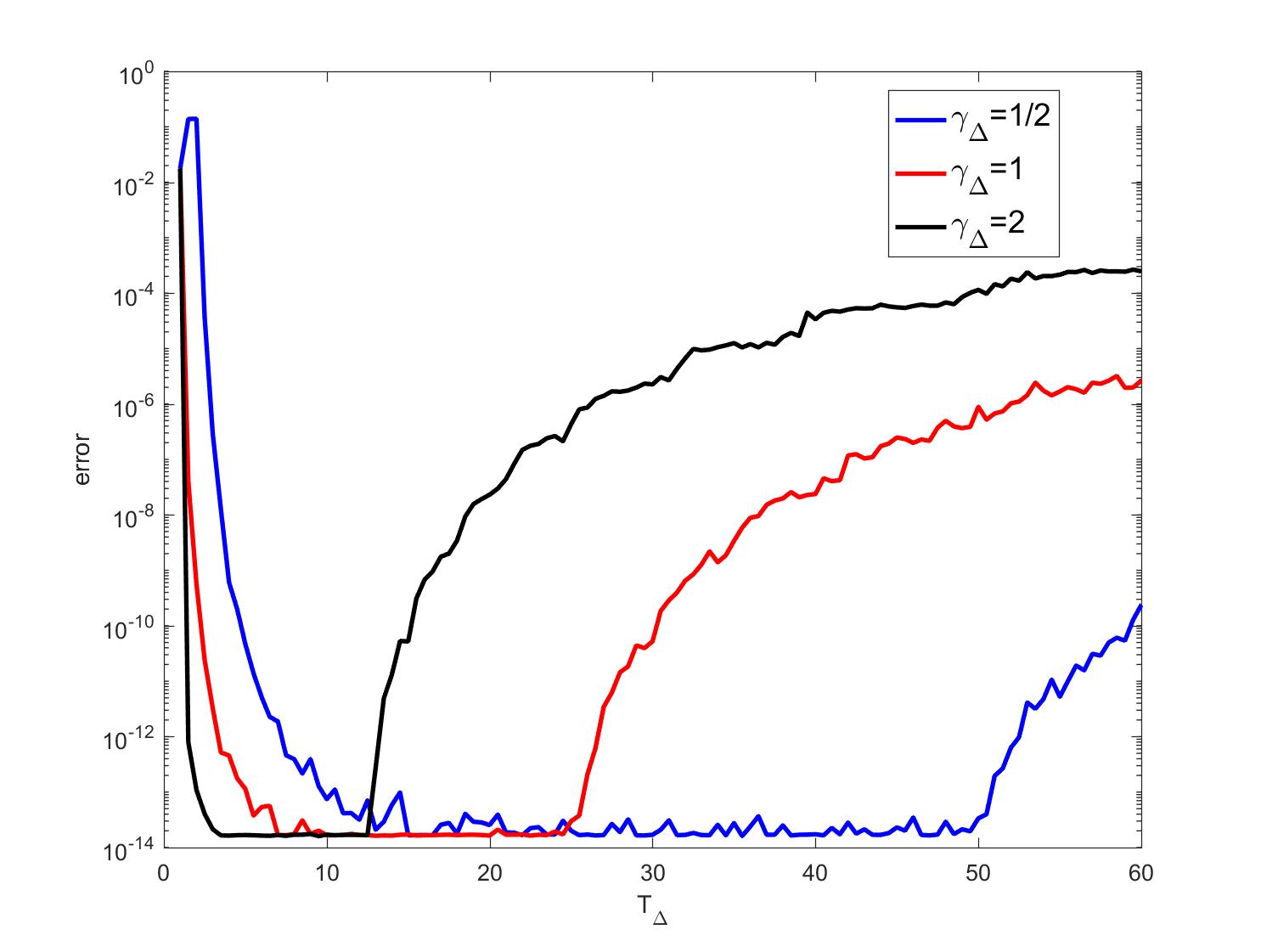}}
}%
\subfigure[\label{4d} $\gamma_{\Delta}=2$, $m_{\Delta}=100$, $\omega=20$.] {
\resizebox*{6cm}{!}{\includegraphics{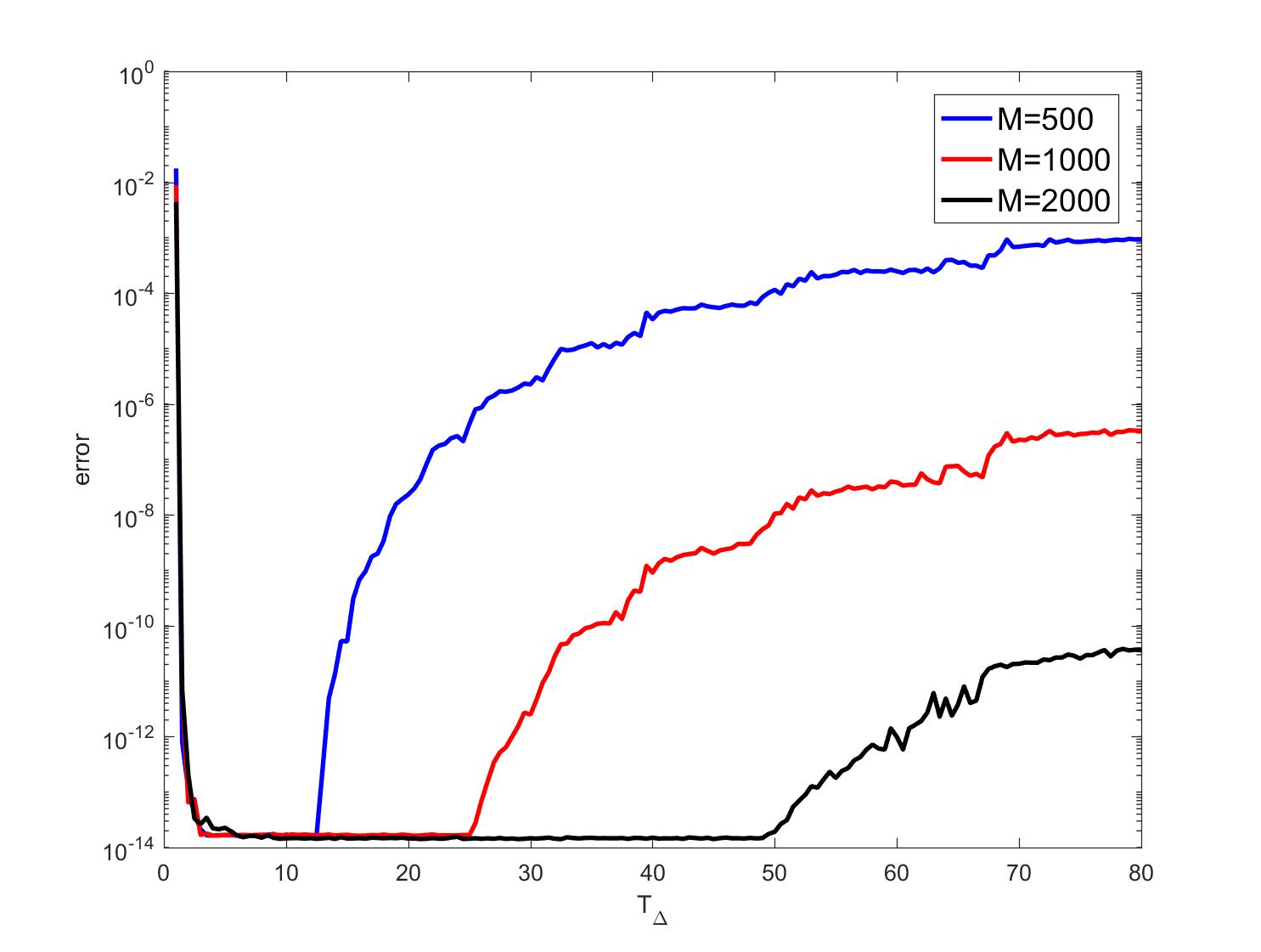}}
}%

		{\caption{The approximation  error against $T_{\Delta}$  for various parameters.}\label{Fig4}}
	\end{center}
\end{figure}

 In this section, we numerically observe the impact of $T_{\Delta}$ on the performance of the algorithm. Now we take
 \begin{equation*}
   f(t)=\exp({\text{i}}\pi\omega t)
 \end{equation*}
 to test the interaction between $T_{\Delta}$ and the other parameters. In Fig \ref{Fig4}, we present the results  of the approximation error against $T_{\Delta}$ for different parameters. Here and what follows, the approximation error is evaluated as the maximum pointwise error over an equi-sampled grid that is ten times denser than that used for construction and the parameter $\tau=1e-14$ (we also tested with $\tau=1e-12,1e-13,1e-15$ and obtained results similar to those shown here.) Three distinct regions
are visible in all graphs in Fig. \ref{Fig4}:
\begin{itemize}
  \item Region $I_1:=\{T_{\Delta}:1<T_{\Delta}<\hat{T}_{1}\}$, where the error decreases rapidly to   the machine accuracy as $T_{\Delta}$ increases.
    \item Region $I_2:=\{T_{\Delta}:\hat{T}_{1}<T_{\Delta}<\hat{T}_{2}\}$, where the error  is maintained close to the machine accuracy.
  \item Region $I_3:=\{T_{\Delta}>\hat{T}_{2}\}$, where the error rebounds, and the increase in $T_{\Delta}$  causes the error to increase instead.
\end{itemize}

\begin{table}
\begin{center}
  \caption{The approximation value of $\hat{T}_{1}$  for various $\gamma$. \label{tab1} }
  \small
{\begin{tabular*}{\textwidth}{@{\extracolsep\fill}cccccccccccccc} \toprule
$\gamma$ &$8$&$4$&$2$&$1$&$\frac{1}{2}$&$\frac{1}{4}$&$\frac{1}{8}$\\\midrule
$\hat{T}_{1}$&1.09&1.20&2.3&5.9&12.5&24.6&48.9\\

\bottomrule
  \end{tabular*}}
\end{center}
\end{table}

It can be seen that the value of $\hat{T}_{1}$  is influenced by $\gamma_{\Delta}$, but not by $\omega$, $m_{\Delta}$, and $M$, while $\hat{T}_{2}$ is influenced by $\gamma_{\Delta}$, $M$, and $\omega$.
For several different values of $\gamma_{\Delta}$, we give the approximate values of $\hat{T}_{1}$ (we tested the functions of $\omega$ from $1$ to $50$, and took the mean of the corresponding $T_{\Delta}$ when their error was less than $1e-13$ for the first time as the estimate of $\hat{T}_1$), which are listed in Table \ref{tab1}. For $\hat{T}_2$, we can  estimate:
\begin{equation}\label{Tlimt}
  \hat{T}_{2}\sim\frac{ M}{\gamma_{\Delta}\omega}.
\end{equation}
This is consistent with the result obtained in equation \eqref{estrrs3}.

\begin{figure}

			\begin{center}
\subfigure[\label{5a} $\gamma 1=$, $M=100$, $m_{\Delta}=50$.] {
\resizebox*{6cm}{!}{\includegraphics{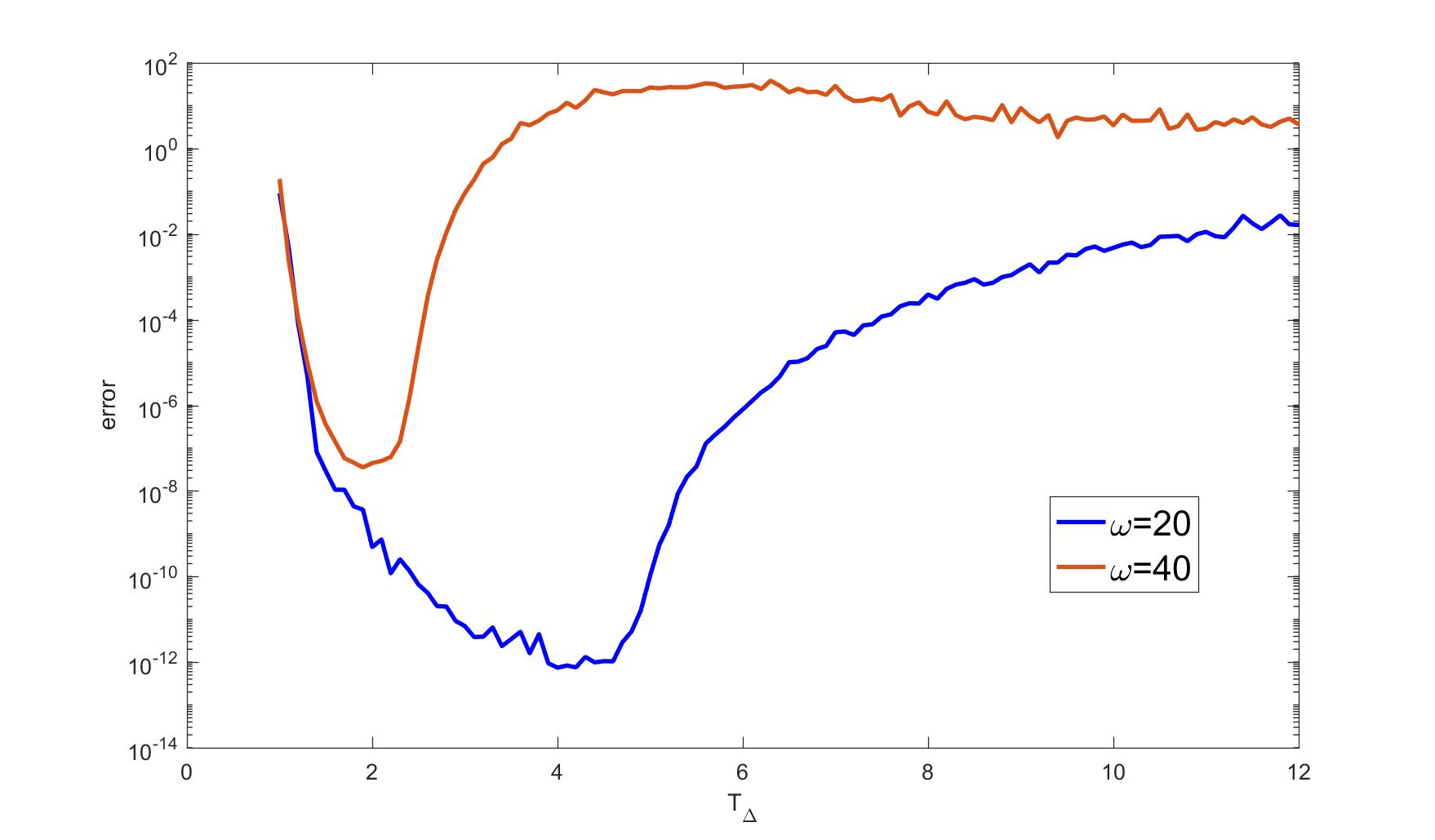}}
}%
\subfigure[\label{5b} $\gamma=2$, $M=400$, $m_{\Delta}=100$.] {
\resizebox*{6cm}{!}{\includegraphics{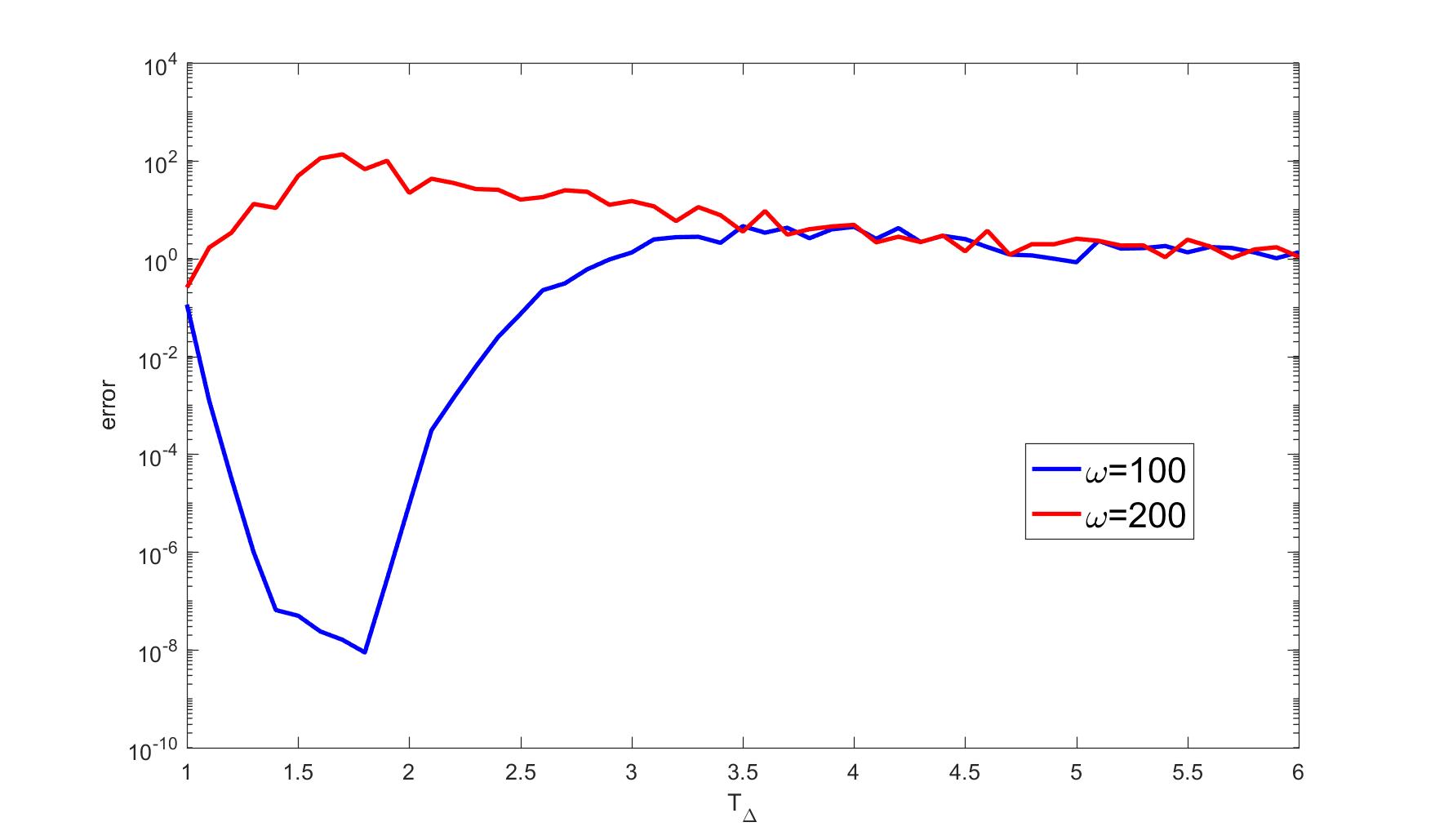}}
}%

		{\caption{The approximation  error against $T_{\Delta}$  for $\hat{T}_{2}<\hat{T}_{1}$.}\label{Fig5}}
	\end{center}
\end{figure}

When the number of nodes $M$ is not sufficiently large,  regions $I_1$ and $I_2$ are merged into a single region. At this time, the error reached the optimal value near  $\hat{T}_{2}$. This phenomenon is shown in   Fig. \ref{Fig5}.

\begin{figure}[H]
			\begin{center}
\subfigure[\label{6a} $\omega=20$, $\gamma=1$,  $m_{\Delta}=50$.
] {
\resizebox*{6cm}{!}{\includegraphics{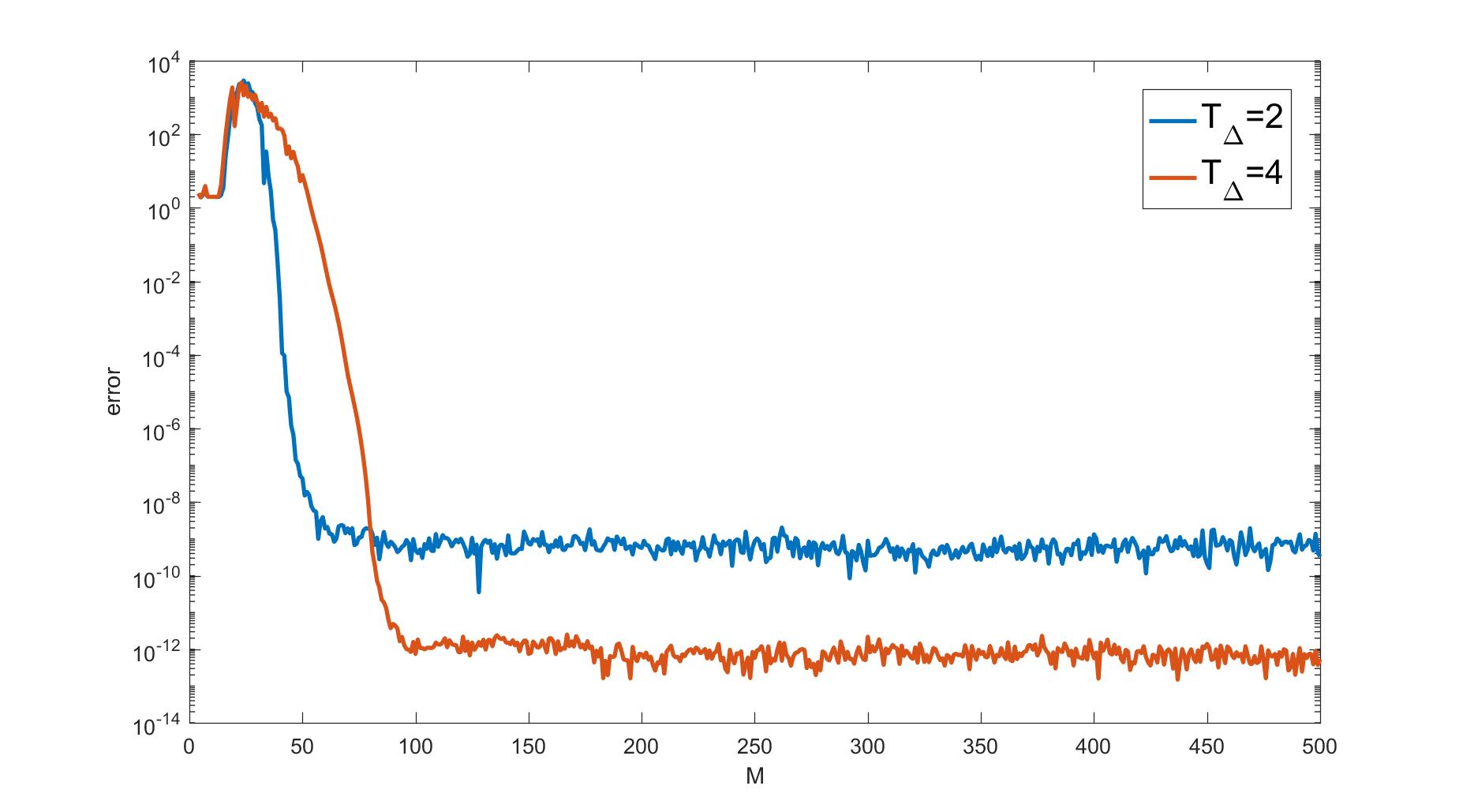}}
}%
\subfigure[\label{6b} $\omega=20$, $\gamma=1$, $M=500$.] {
\resizebox*{6cm}{!}{\includegraphics{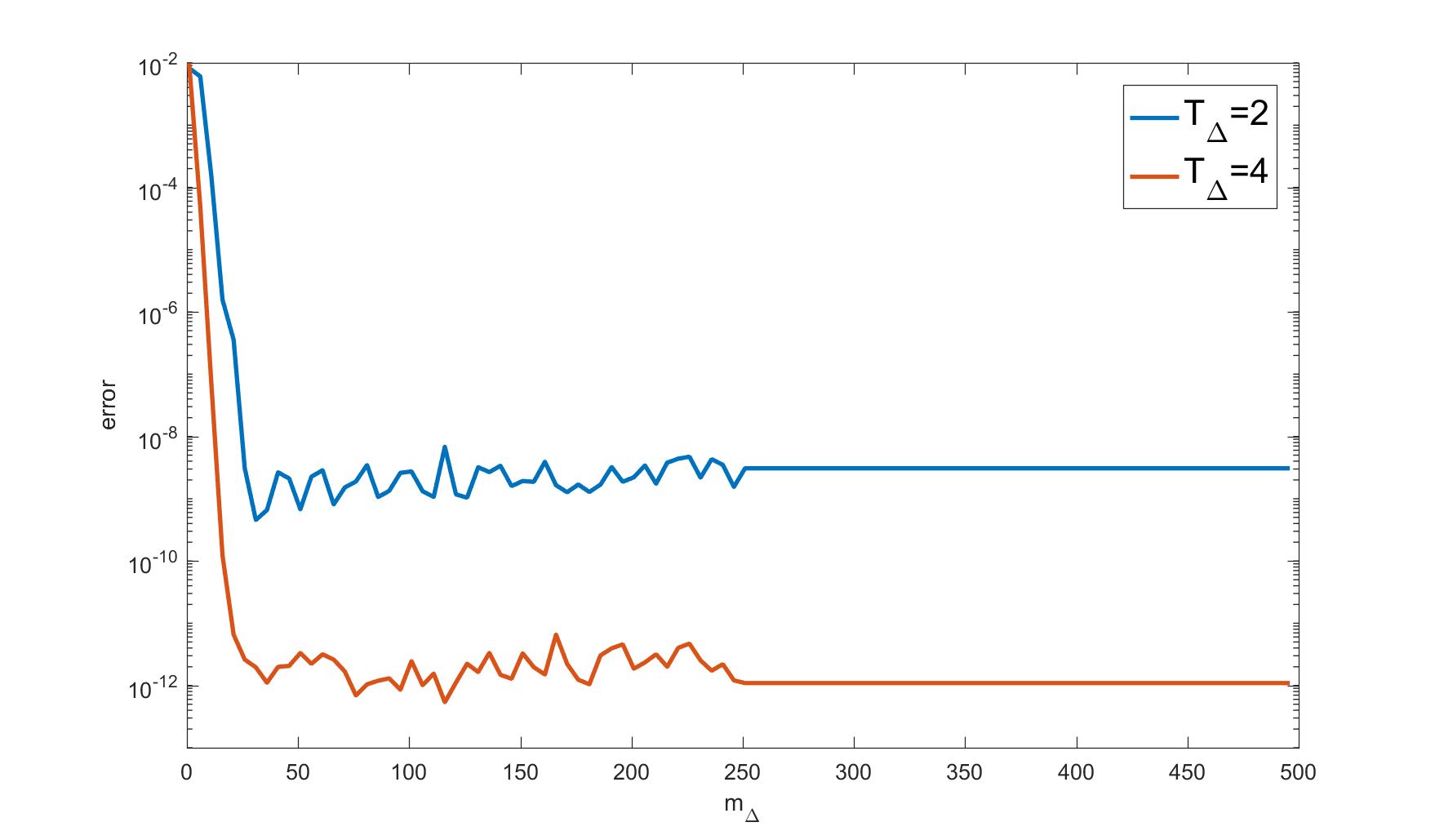}}
}%

		{\caption{The approximation  error against $T_{\Delta}$  for ${T}_{\Delta}<\hat{T}_{1}$ (If $m_{\Delta}>M$, we take $m_{\Delta}=M$).}\label{Fig6}}
	\end{center}
\end{figure}

\begin{figure}[H]
			\begin{center}
\subfigure[\label{7a} $\omega=20$, $\gamma=1$.] {
\resizebox*{6cm}{!}{\includegraphics{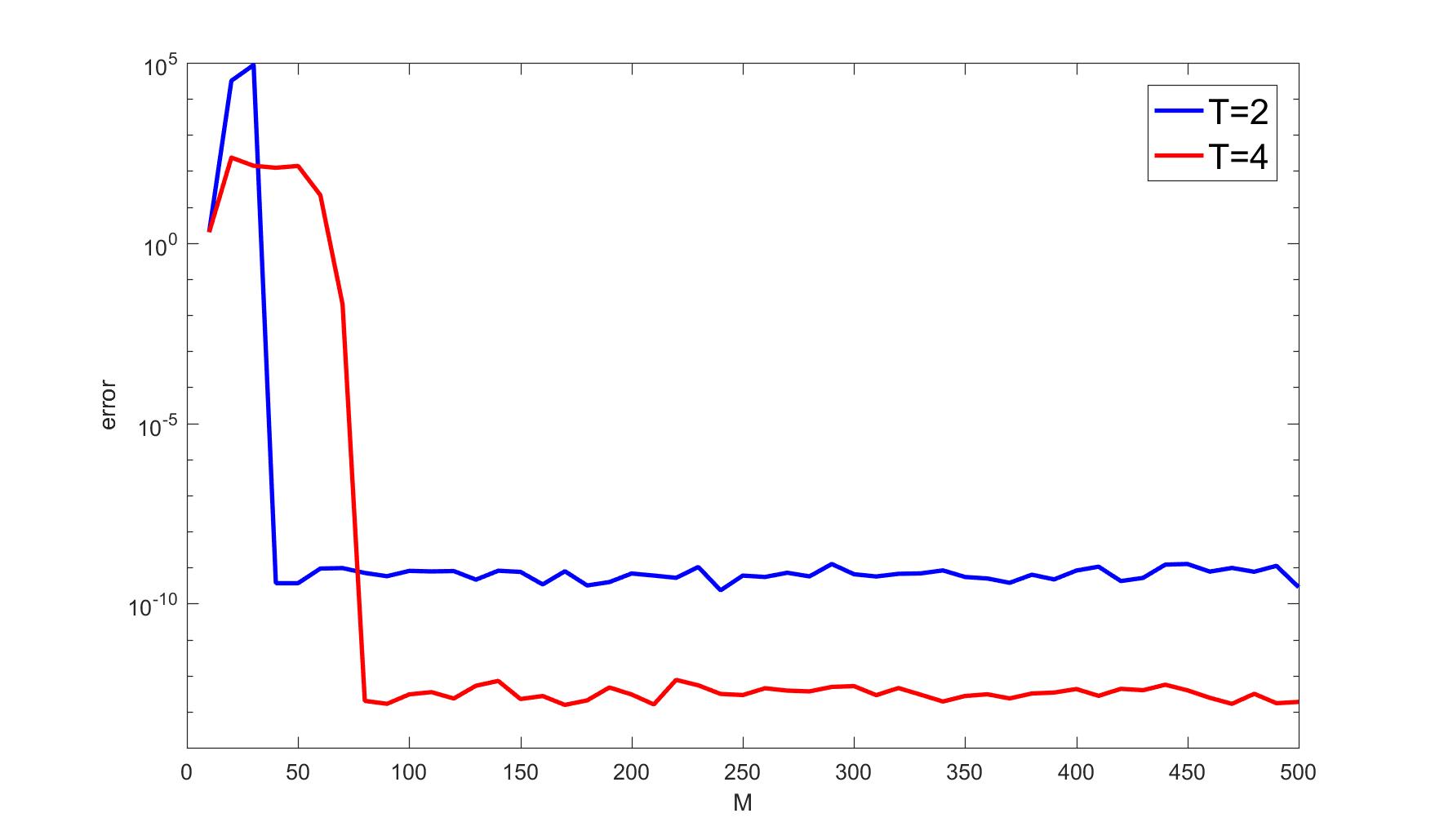}}
}%
\subfigure[\label{7b} $\omega=20$, $\gamma=1/2$.] {
\resizebox*{6cm}{!}{\includegraphics{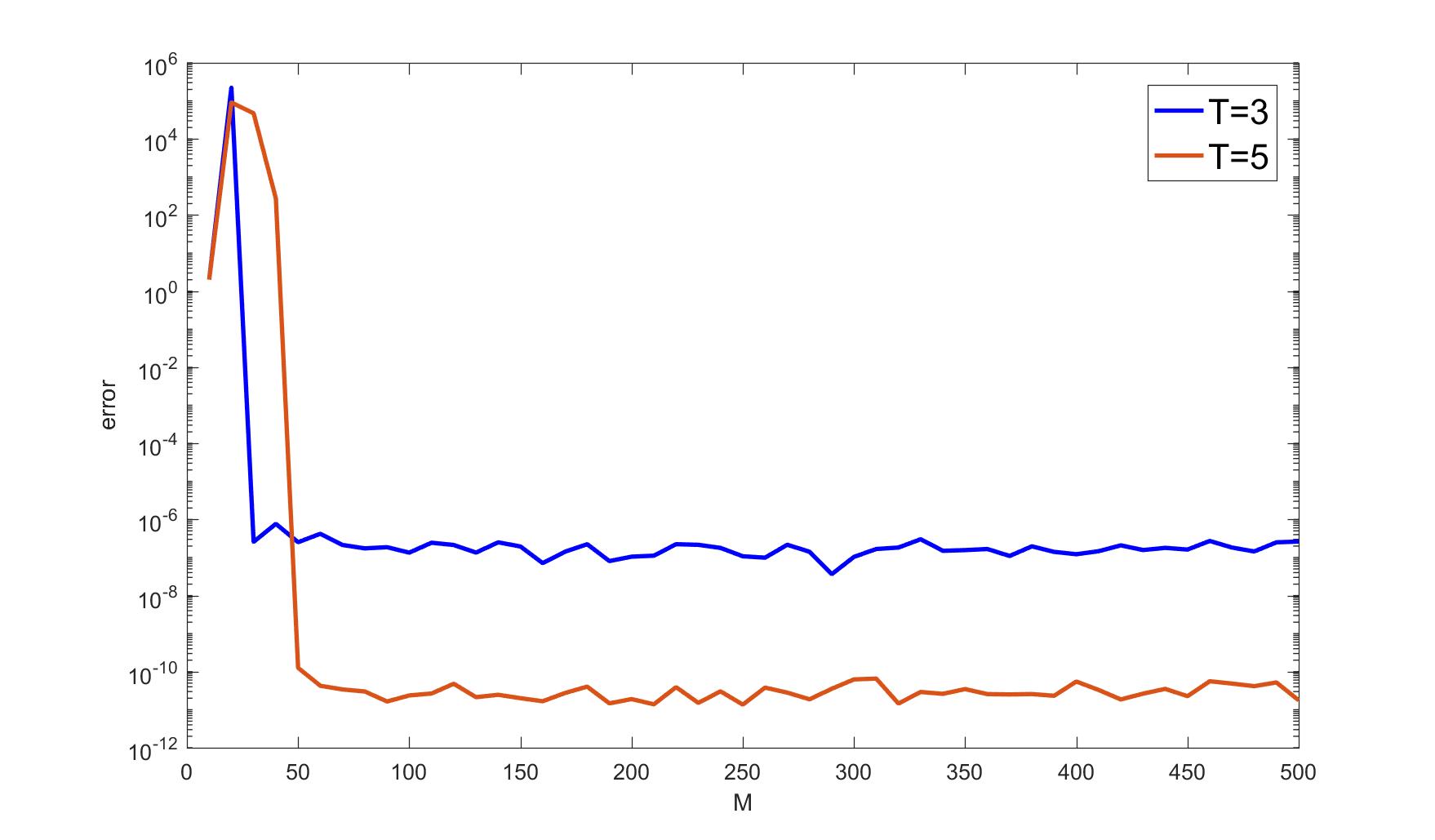}}
}%

		{\caption{The approximation  error against $M$  for full data FE with $T<\hat{T}_1$.}\label{Fig7}}
	\end{center}
\end{figure}

What will happen if the $T_{\Delta}$ value selected  is less than $\hat{T}_{1}$? Figure \ref{Fig6} shows this. We can see that the error will not decrease after reaching a certain threshold, which cannot be improved by increasing the numbers of $M$ and $m_{\Delta}$. It should be noted that this situation also applies to the Fourier extension algorithm using the full data in previous literatures. In Fig. \ref{Fig7}, we show the results of the algorithm using full data. It can be seen that when $\gamma$ is fixed and the value of $T$ is smaller than the corresponding $\hat{T}_1$, the method cannot achieve machine accuracy. Similarly, this cannot be improved by increasing the number of $M$.

From the  test results above, we can see that when $\gamma_{\Delta}$ is fixed, the selection of $T_{\Delta}$ is crucial to   the approximation accuracy.  The approximation can obtain satisfactory
results only when $T_{\Delta}$ is within the interval $[\hat{T}_1, \hat{T}_2]$. The left endpoint $\hat{T}_1$ of the interval is determined by   parameter $\gamma_{\Delta}$, whereas the right endpoint $\hat{T}_2$ is affected by   parameters $\gamma_{\Delta}$, $\omega$ and $M$.

\begin{figure}
			\begin{center}
\subfigure[\label{8a} $M=1000$, $\gamma_{\Delta}=2$, $T_{\Delta}=6$.] {
\resizebox*{6cm}{!}{\includegraphics{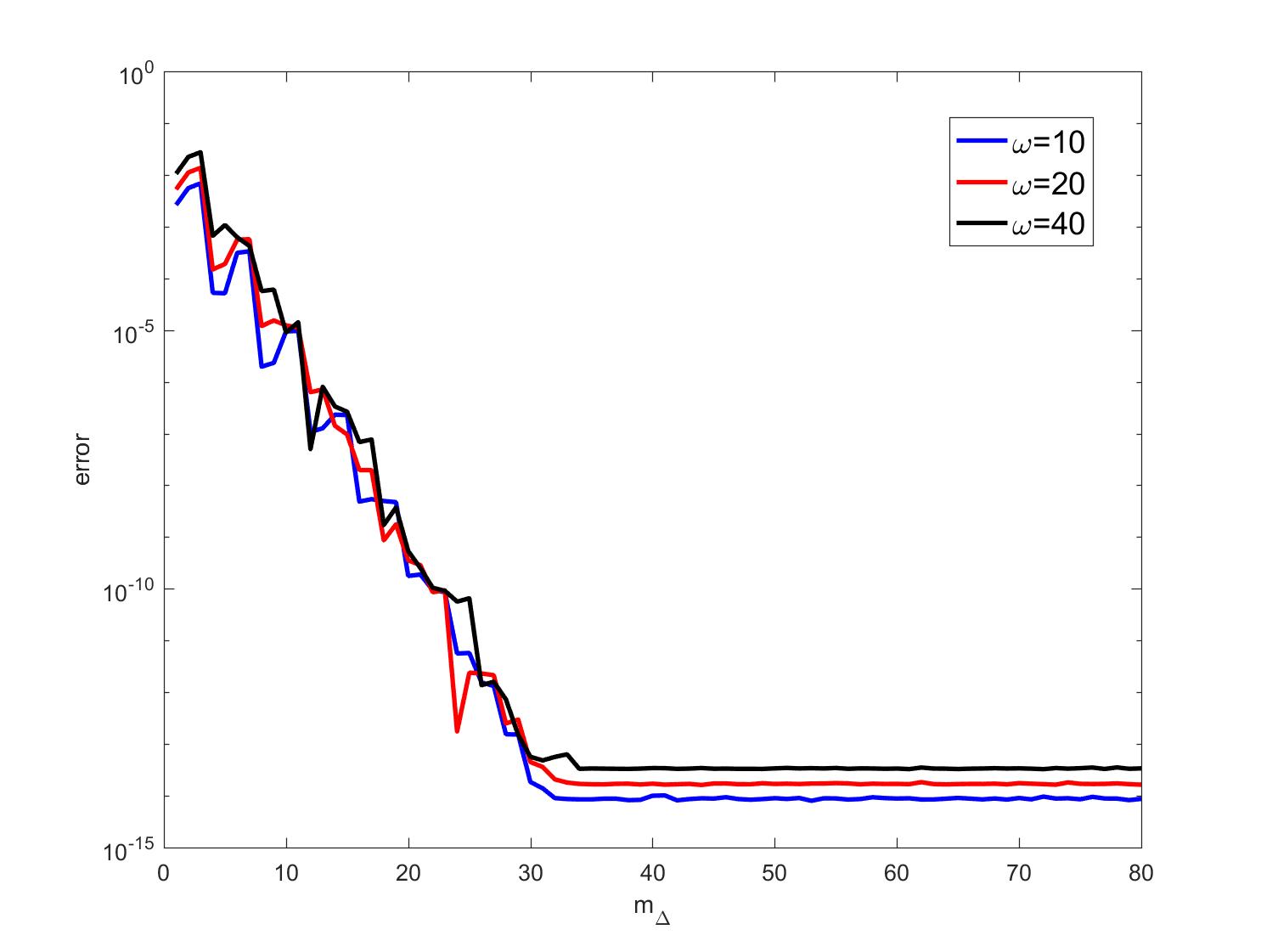}}
}%
\subfigure[\label{8b}  $\omega=20$,  $\gamma_{\Delta}=2$, $T_{\Delta}=6$.] {
\resizebox*{6cm}{!}{\includegraphics{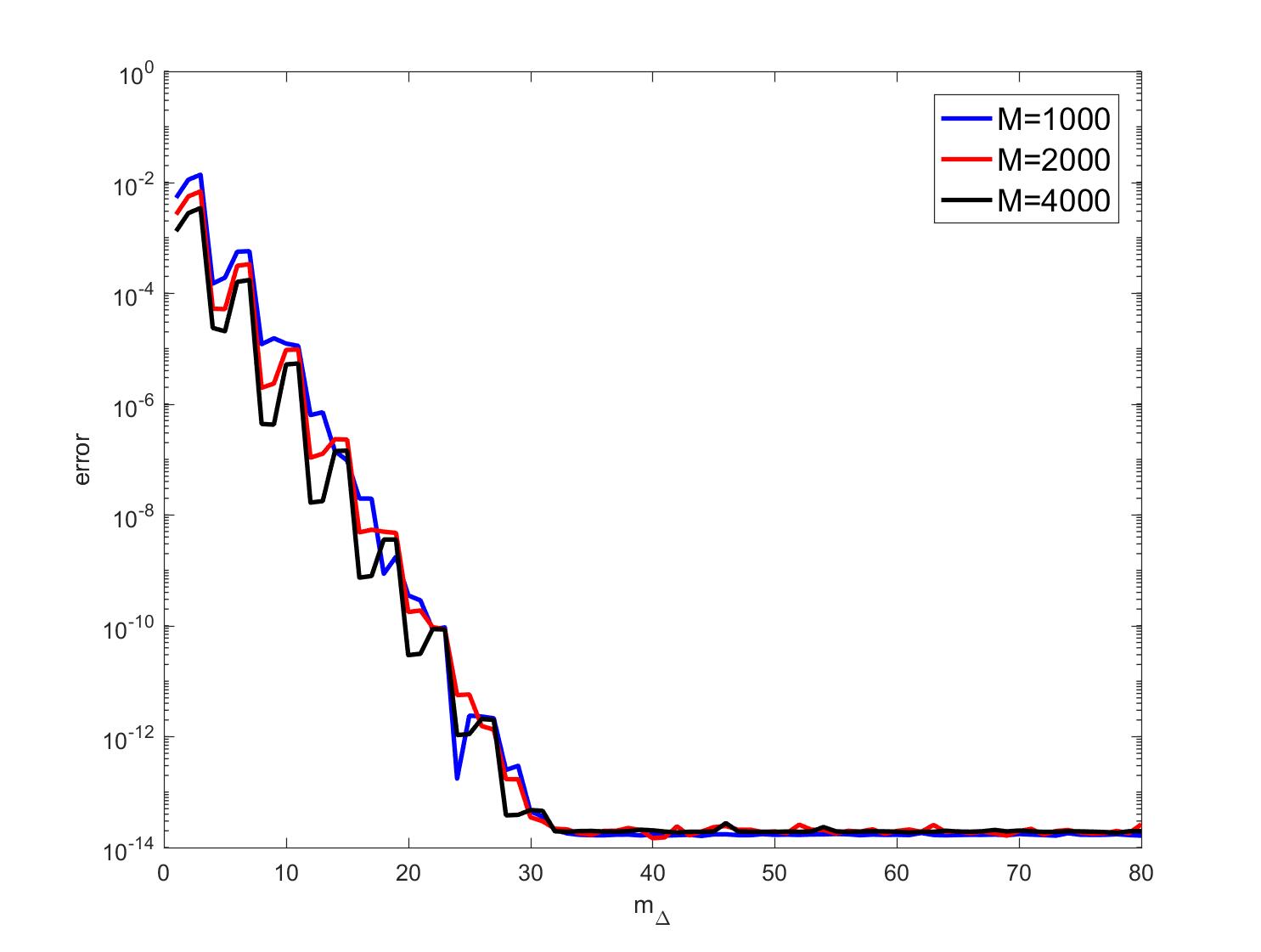}}
}%
		{\caption{The approximation  error against $m_{\Delta}$  for various $\omega$ and $M$.}\label{Fig8}}
	\end{center}
\end{figure}
 \begin{figure}
			\begin{center}
\subfigure[\label{9a} $M=2000$, $\gamma_{\Delta}=8$, $\omega=20$.] {
\resizebox*{6cm}{!}{\includegraphics{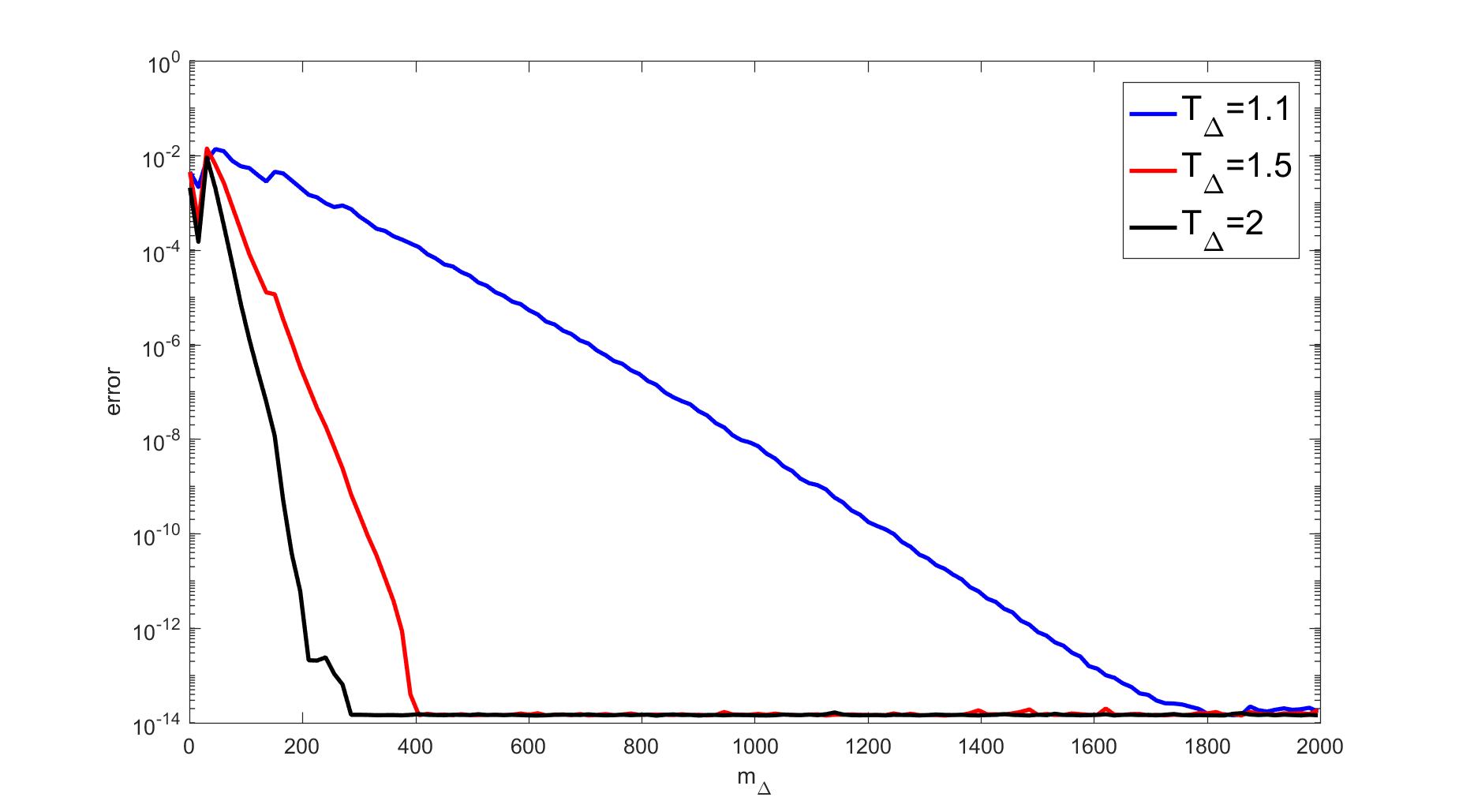}}
}%
\subfigure[\label{9b} $M=1000$, $\gamma_{\Delta}=4$, $\omega=20$.] {
\resizebox*{6cm}{!}{\includegraphics{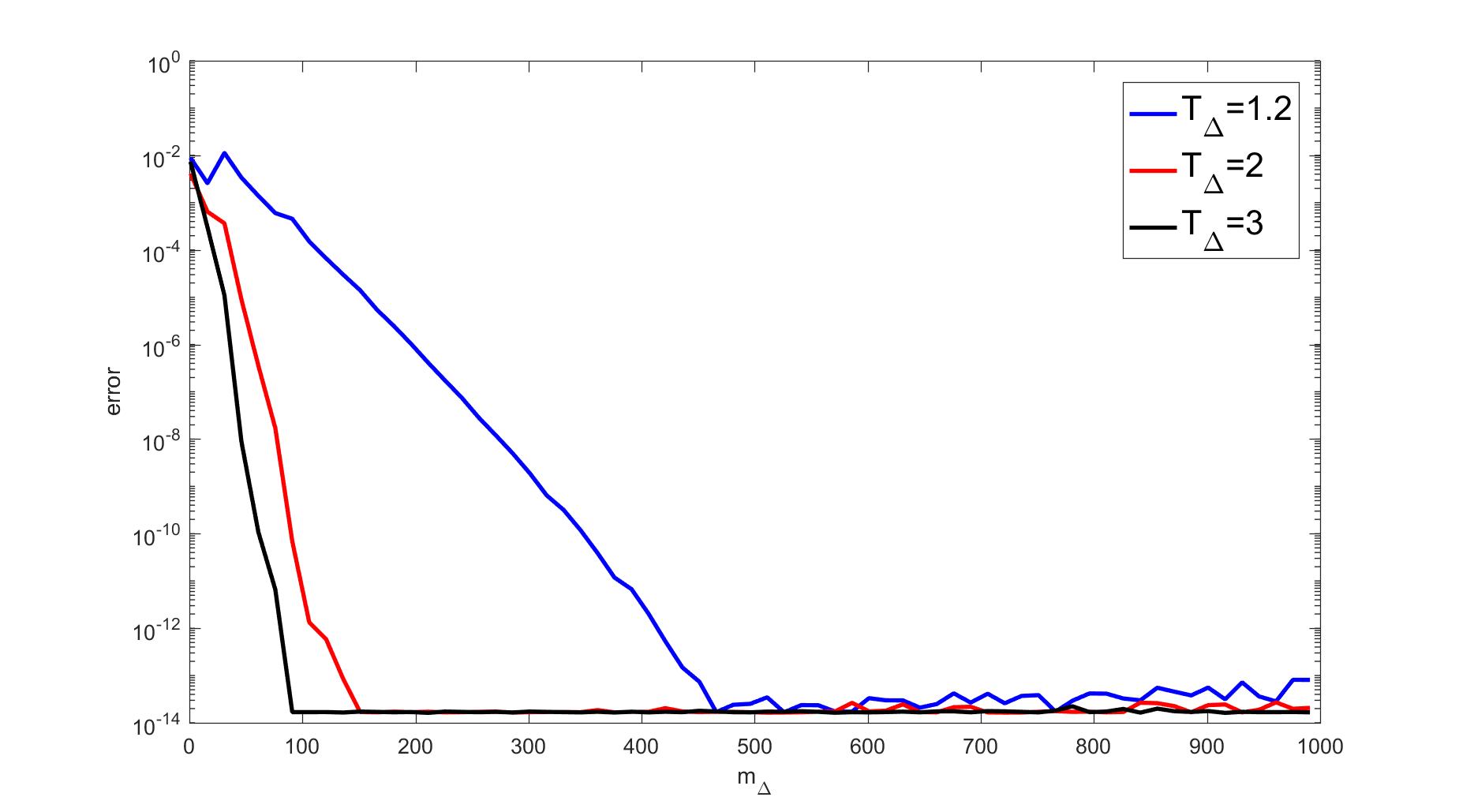}}
}%

\subfigure[\label{9c}  $M=500$, $\gamma_{\Delta}=2$, $\omega=20$.] {
\resizebox*{6cm}{!}{\includegraphics{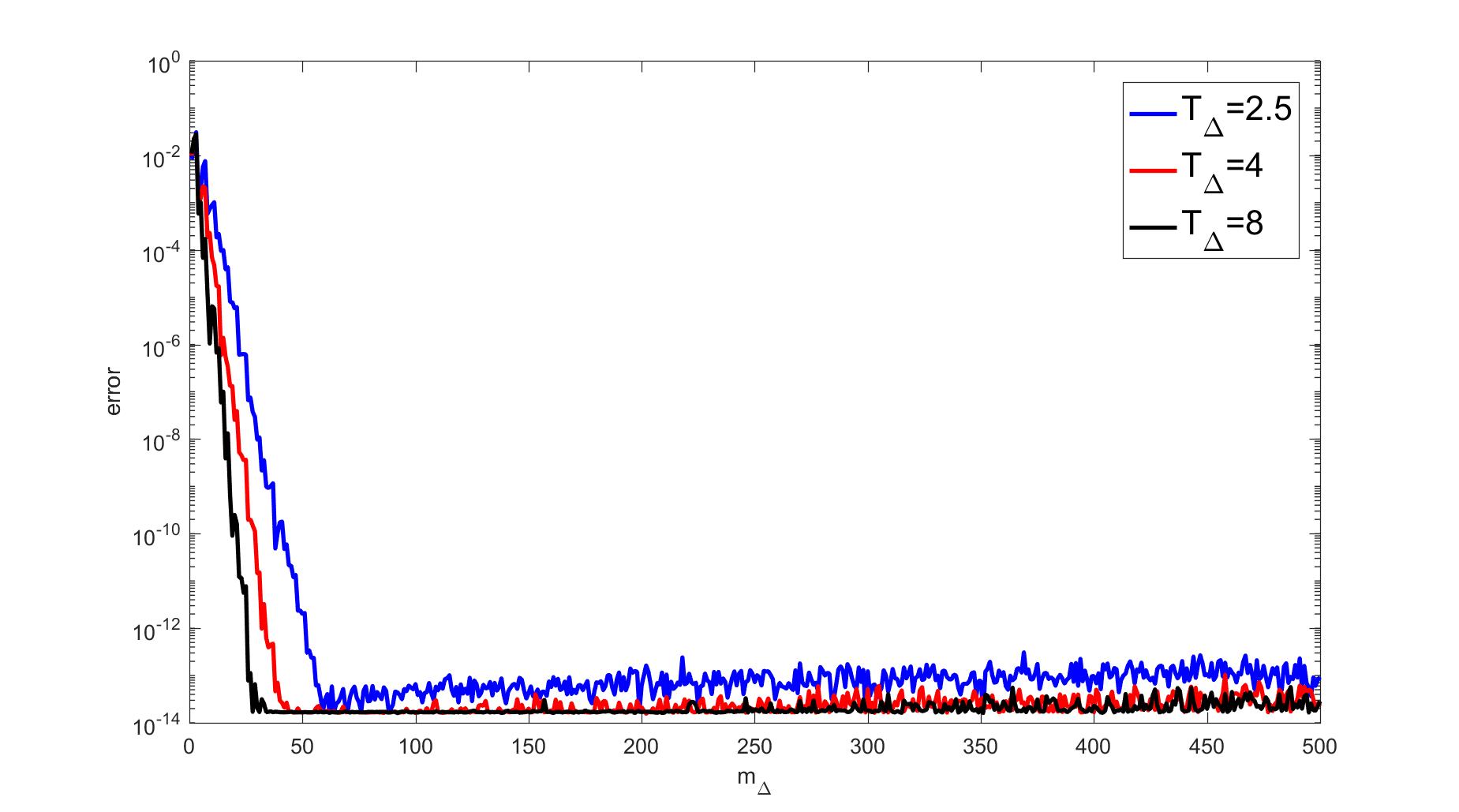}}
}%
\subfigure[\label{9d}  $M=500$, $\gamma_{\Delta}=1$, $\omega=20$.] {
\resizebox*{6cm}{!}{\includegraphics{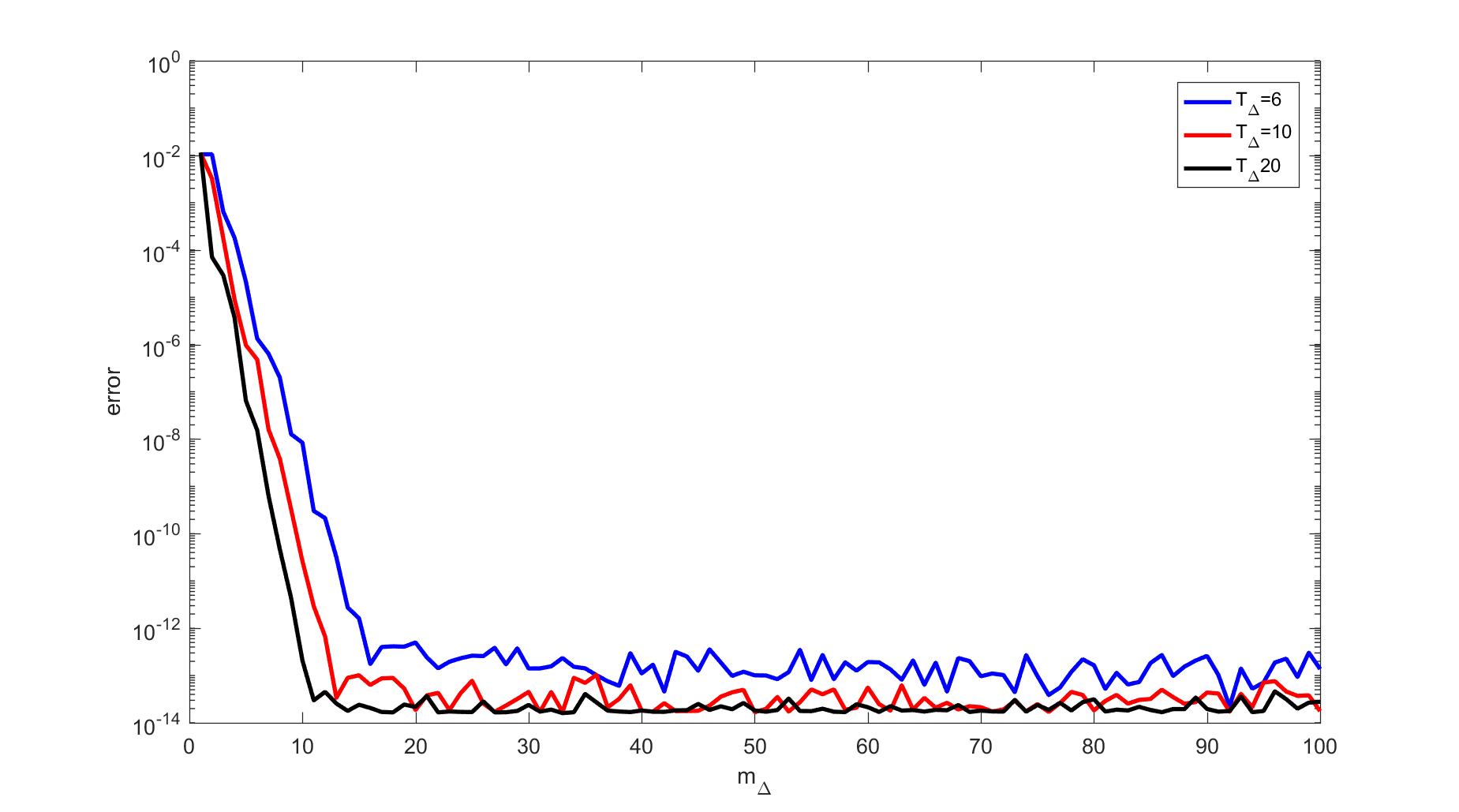}}
}

\subfigure[\label{9e}  $M=500$, $\gamma_{\Delta}=1/2$, $\omega=20$.] {
\resizebox*{6cm}{!}{\includegraphics{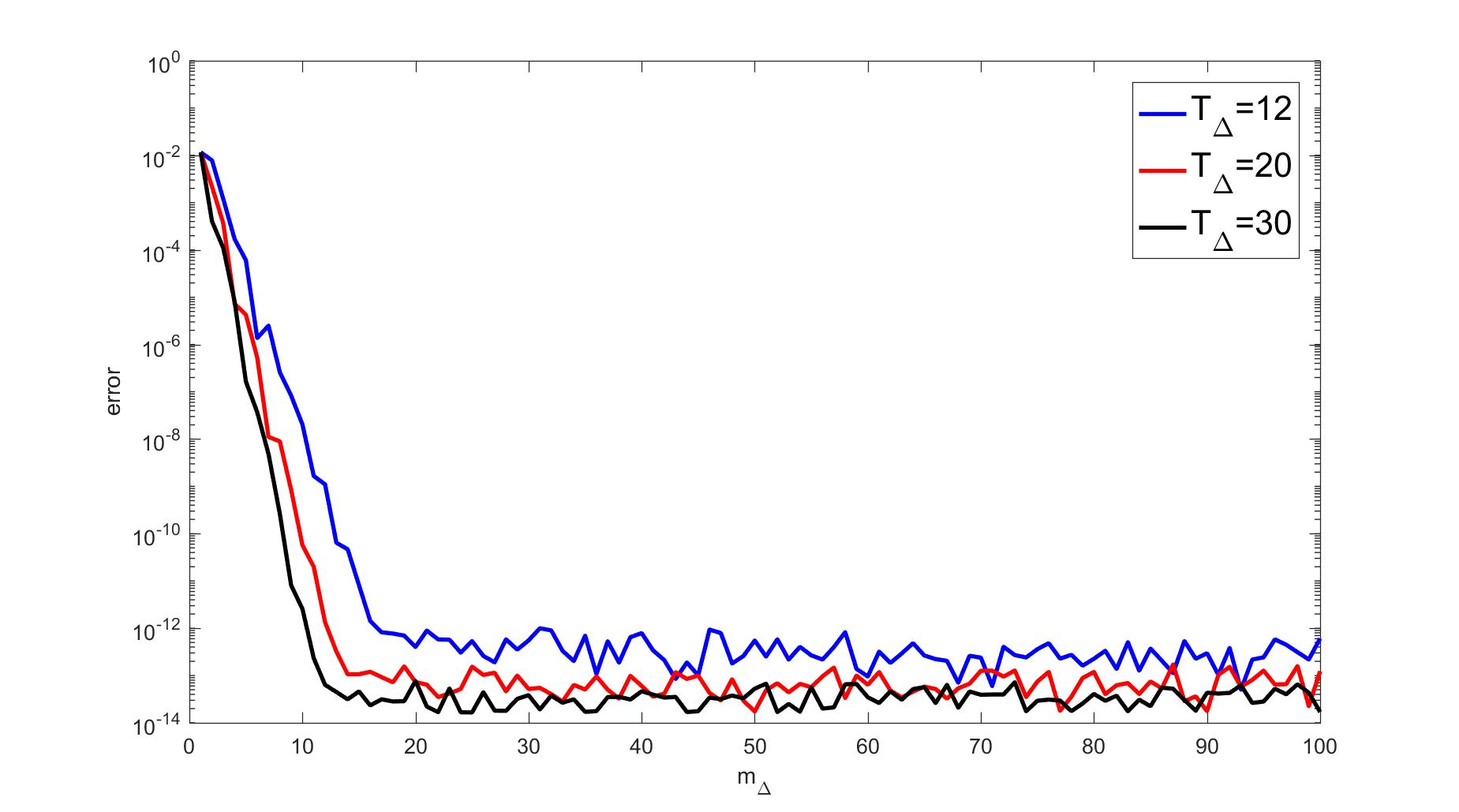}}
}\subfigure[\label{9f}  $M=500$, $\gamma_{\Delta}=1/4$, $\omega=20$.] {
\resizebox*{6cm}{!}{\includegraphics{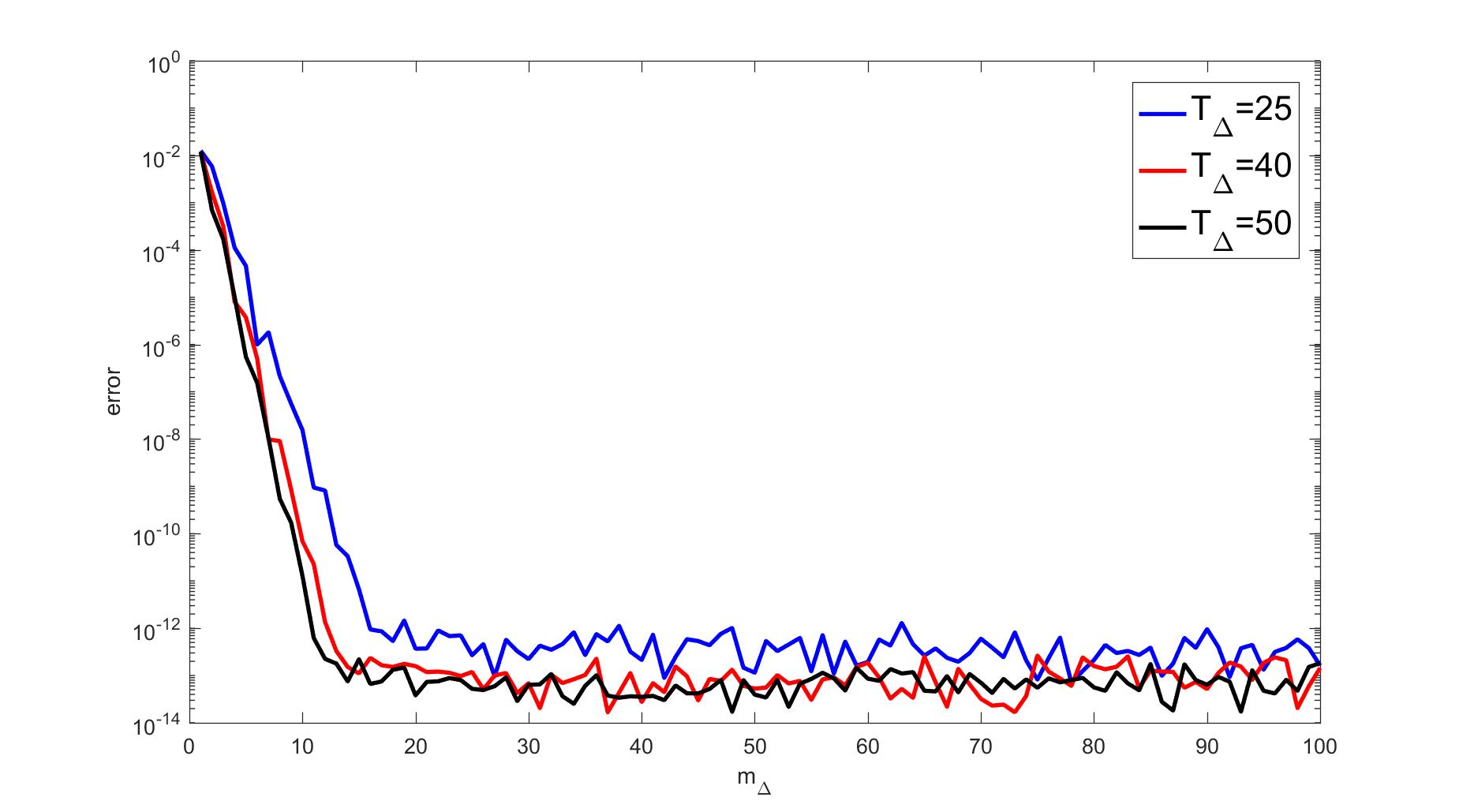}}
}
	{\caption{The approximation  error against $m_{\Delta}$  for various $\gamma_{\Delta}$ and $T_{\Delta}$.}\label{Fig9}}
	\end{center}
\end{figure}

\subsection{Testing the parameter $m_{\Delta}$}
The size of the parameter $m_{\Delta}$ determines whether this study is meaningful. This subsection tests the influence of  the other parameters on it. First, we take $M>{T}_{\Delta}\gamma_{\Delta}\omega$ for observation. As shown in Fig. \ref{Fig8}, as the value of $m_{\Delta}$ increases, the error decreases rapidly and reaches the machine precision after a certain value $\hat{m}_{\Delta}$, while the value of $\hat{m}_{\Delta}$  is independent of $M $ and $\omega$.

 In Fig. \ref{Fig9}, we further show how the error changes with $m_{\Delta}$ for different values of $\gamma_{\Delta}$ and $T_{\Delta}$. It can be seen that the value of $\hat{m}_{\Delta}$ decreases with the increase of $T_{\Delta}$. This is consistent with the previous analysis, because the step size $h=\frac{2\pi}{L_{\Delta}}$ in the $g_c$ calculation process is determined by the product $(m_{\Delta}-1)T_{\Delta}$. When $\gamma_{\Delta} \leq  1$, the value of $\hat{m}_{\Delta}$ at the corresponding $\hat{T}_1$ is very close. When $\gamma _{\Delta}> 1$, the value of $\hat{m}_{\Delta}$ at the corresponding $\hat{T}_1$ increases rapidly with the increase of $\gamma_{\Delta}$. From the perspective of use, we are more concerned about the $\hat{m}_{\Delta}$ value near $\hat{T}_1$. We tested this and the results are listed in Table \ref{tab2}.

\begin{table}
\begin{center}
  \caption{The approximation values of $\hat{m}_{\Delta}$ for various $\gamma$ and $T_{\Delta}$. \label{tab2} }
  \small
{\begin{tabular*}{\textwidth}{@{\extracolsep\fill}ccccccccccccccccccccc} \toprule

$\gamma_{\Delta}$ &{ $8$}&{ $4$}&{ $2$}&{ $1$}&{ $\frac{1}{2}$}&{ $\frac{1}{4}$}&{ $\frac{1}{8}$}\\\hline
$T_{\Delta}$&$1.1$&$1.2$&$2.3$&$6$&$12$&$24$&$50$\\\midrule
$\hat{m}_{\Delta}$&1870&490&64&23&25&24&26\\
\bottomrule
  \end{tabular*}}
\end{center}
\end{table}

\begin{figure}
			\begin{center}
\subfigure[\label{10a} $f(t)=\cos\left(\frac{100}{1+25t^2}\right)$] {
\resizebox*{6cm}{!}{\includegraphics{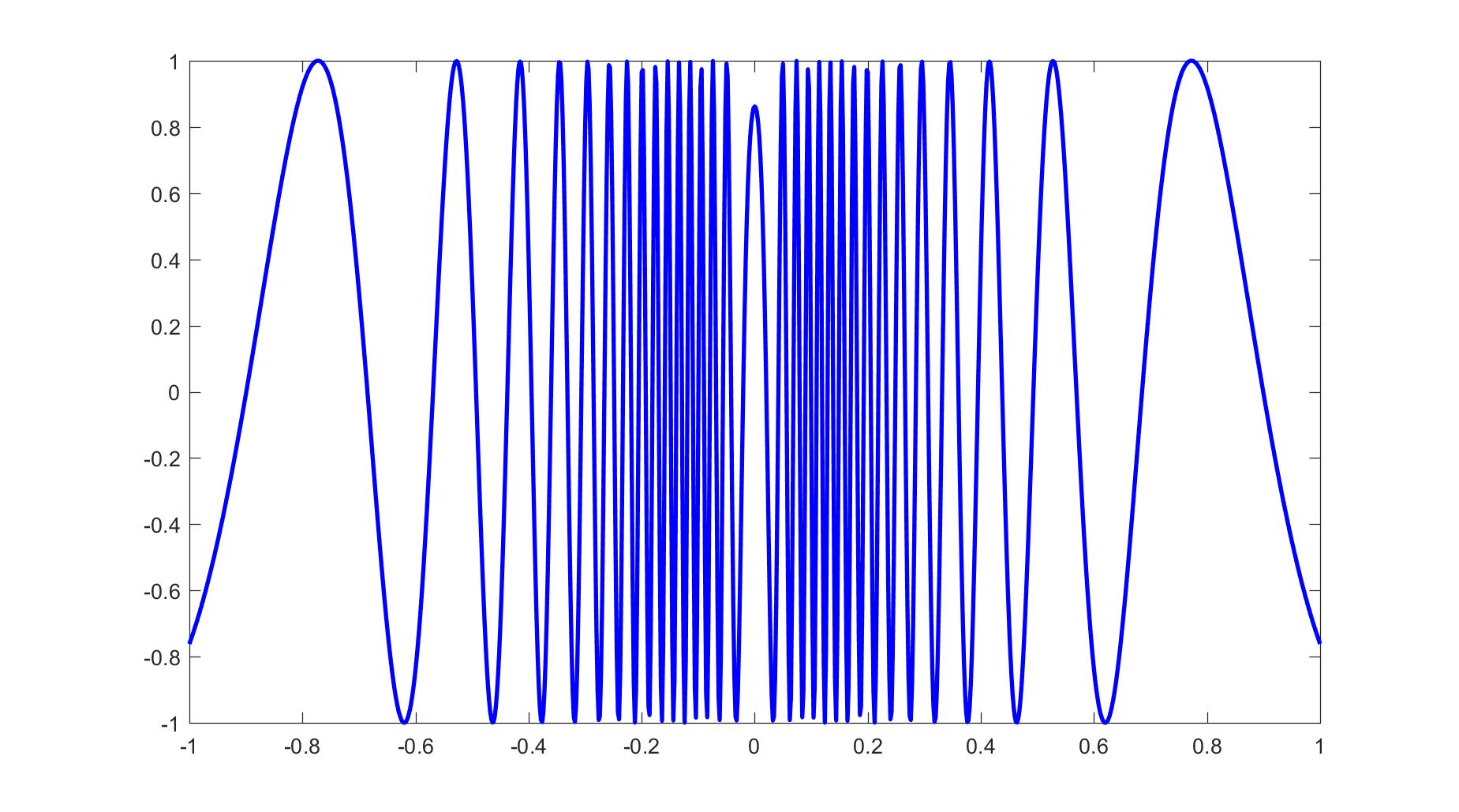}}
}%
\subfigure[\label{10b} approximation error] {
\resizebox*{6cm}{!}{\includegraphics{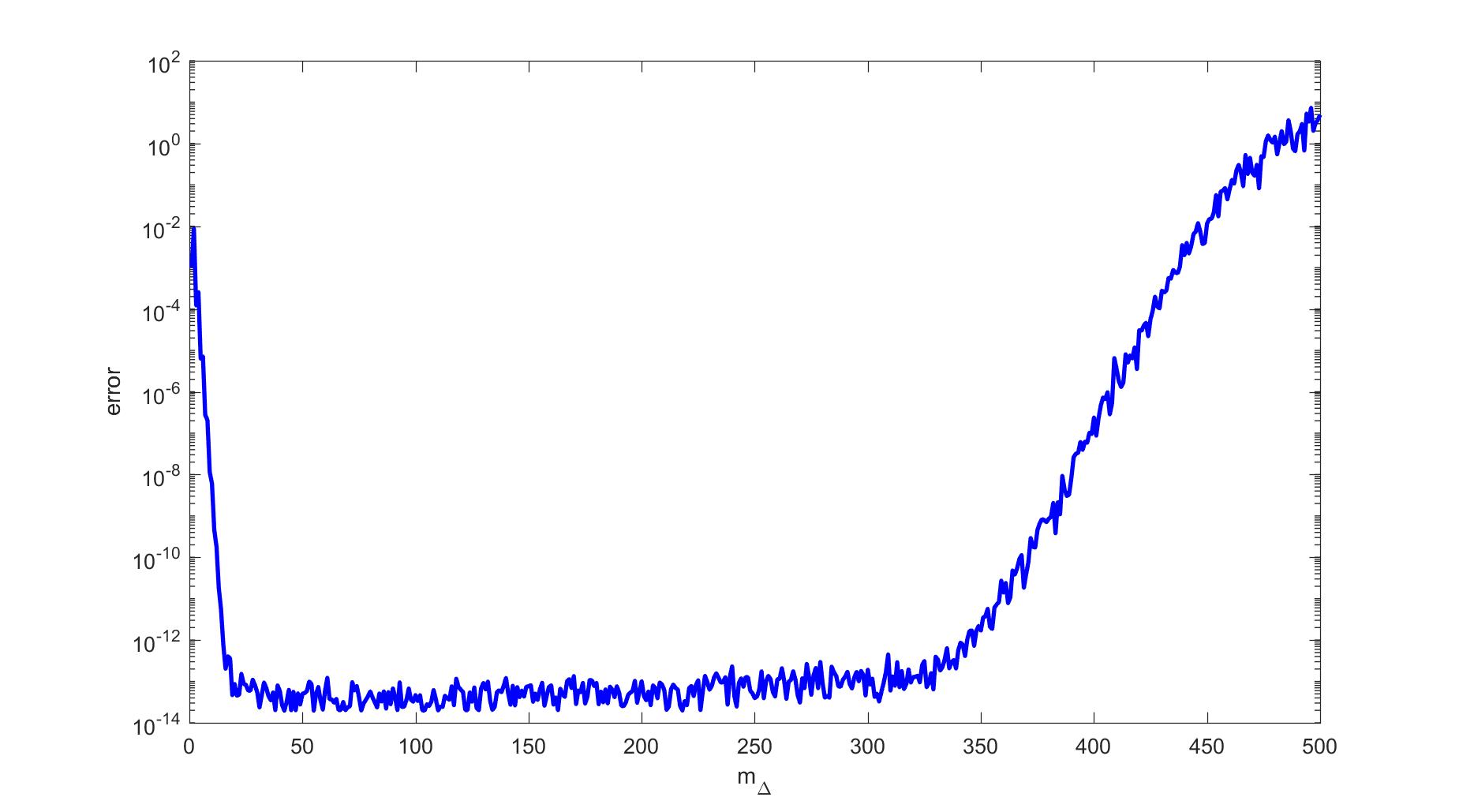}}
}%

	{\caption{ The function and its approximation errors against $m_{\Delta}$ for $M=500$, $\gamma_{\Delta}=1$, $T_{\Delta}=6$.}\label{Fig10}}
	\end{center}
\end{figure}

It is also worth noting that for functions with internal oscillations but are smooth near the boundary, as mentioned in Section \ref{SEC5}, an increase in the number of $m_{\Delta}$ will cause  the value of $\omega$ in \eqref{estrrs3} to increase. If the value of $M$ is not sufficiently large, the approximation result using a larger $m_{\Delta}$ may be worse. This can be seen from the approximation results for the function $f(t)=\cos\left(\frac{100}{1+25t^2}\right)$ in Fig. \ref{Fig10}.

\subsection{Parameter determination}

Combining the test results in the previous two subsections, we know that if we want to approximate the function with frequency $\omega$ with machine precision, the parameters $T_{\Delta}$  and $m_{\Delta}$ must satisfy
\begin{equation}
  T_{\Delta}\geq \hat{T}_1,\quad m_{\Delta}\geq \hat{m}_{\Delta}.
\end{equation}
Accordingly, combined with   \eqref{Tlimt}, the parameters $L_{\Delta}$ and $M$ must satisfy
\begin{equation}
  L_{\Delta}\geq \hat{L}_{\Delta}=:2\times \hat{T}_1\times \hat{m}_{\Delta},\quad M\geq \hat{M}=: \hat{T}_1\gamma_{\Delta}\omega.
\end{equation}
Referring to the results in Table \ref{tab1} and Table \ref{tab2}, the value of $\hat{M}$ is very close for different $\gamma_{\Delta}$ values,  and the value of $\hat{L}_{\Delta}$ is relatively close when $\gamma_{\Delta}=1$ and $2$; however,  the value of $\hat{L}_{\Delta}$ is much larger when $\gamma_{\Delta}$ takes other values. Further combined with the value of $\hat{m}_{\Delta}$, it can be seen that the computational complexity of the method is equivalent when $\gamma_{\Delta}=1$ and $\gamma_{\Delta}=2$, whereas the computational complexity increases significantly in other cases. Furthermore,   a larger $m_{\Delta}$ indicates a greater probability of encountering high frequencies within the selected boundary interval. Therefore, we recommend using the following  parameters
\begin{equation}\label{parameterset}
  \gamma_{\Delta}=1, T_{\Delta}=6, m_{\Delta}=25,
\end{equation}
for actual calculations, and parameters
\begin{equation}
  \gamma_{\Delta}=2, T_{\Delta}= 2.3, m_{\Delta}=65,
\end{equation}
can be used as an alternative. Under such parameter settings, the size of   matrix $A$ in  \eqref{lseq} is small, and  singular value decomposition can be calculated directly. In addition, note that $A$ is fixed in the algorithm, so the SVD of $A$ can be pre-computed and stored.
\section{Testing and further improvement of the algorithm\label{SEC6}}
\subsection{Testing and comparison of algorithm performance}
\begin{figure}
	\begin{center}
		{\resizebox*{7cm}{5cm}{\includegraphics{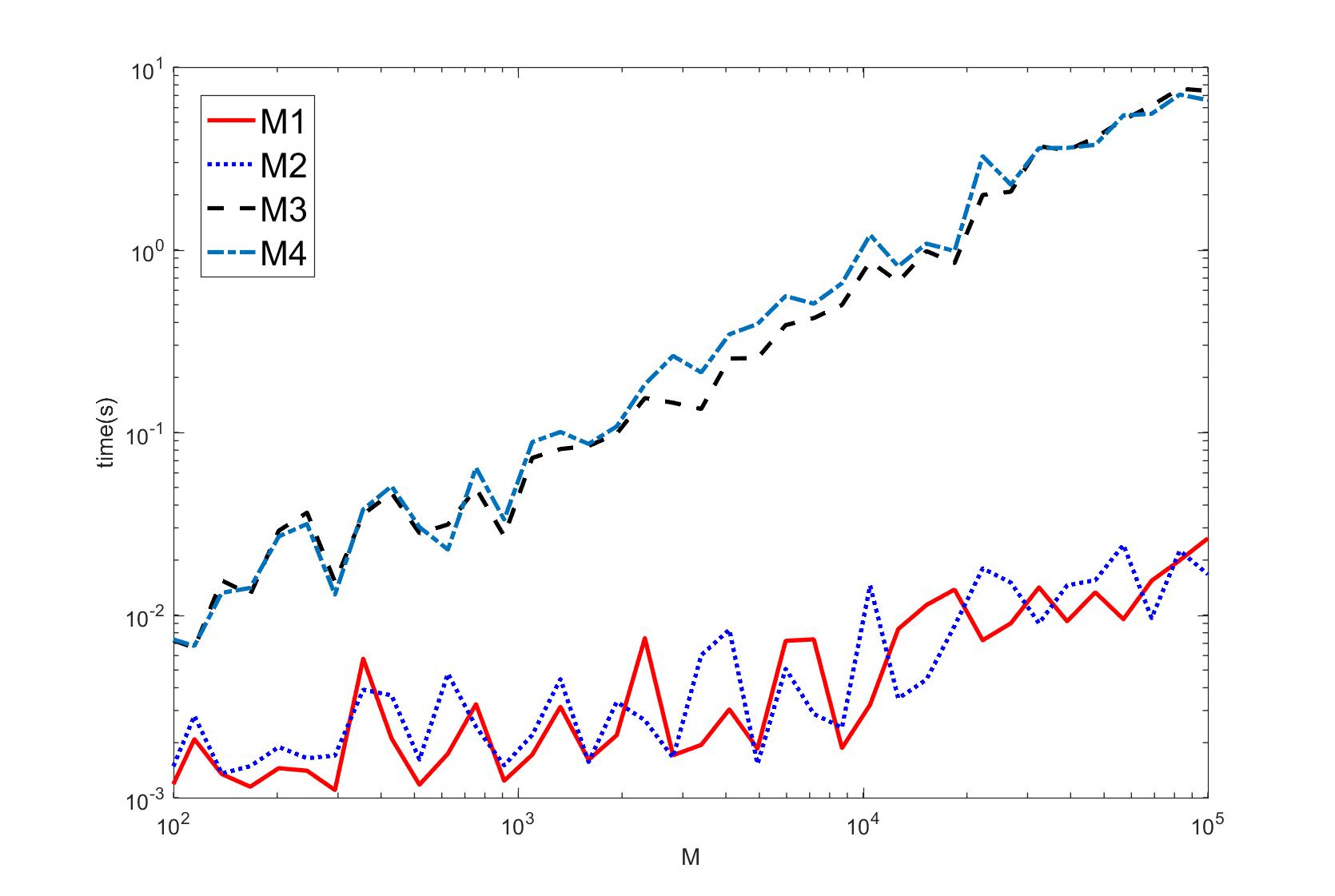}}}
		{\caption{Execution time for increasing the number of nodes $M$, for MATLAB implementations of the Boundary Interval algorithm (M1), the FC-Gram algorithm (M2), the Lyon algorithm (M3) and Implicit projection algorithm (M4). \label{Figtime}}}
	\end{center}
\end{figure}
Based on the previous tests, we have provided specific values for the parameters involved in the algorithm. In this subsection,  the performance of the proposed algorithm is tested. All  tests were computed on a Windows 10 system with  16 GB of memory,
Intel(R) Core(TM)i7-8500U CPU\@1.80GHz using MATLAB 2016b. We   compared the boundary interval algorithm (M1) proposed in this paper with the FC-Gram algorithm (M2) in \cite{bruno2010high}, the Lyons algorithm (M3) in \cite{Lyon2011} and the implicit projection algorithm (M4) in \cite{Matthysen2016FAST}.
{ For M2, we use the codes at https://github.com/oscarbruno/FC and take the parameters $d=13$, $C=33$. For the specific meaning of the parameters, please refer to \cite{2016An}. We also tested other parameter values. As $d$ increases, the algorithm becomes unstable. A large $d$ value (such as $d=30$) will prevent the algorithm from reaching near machine precision. According to our test results, if our approximation goal is to achieve  near machine accuracy when $M$ is relatively small, then the optimal value of $d$ is approximately $12\sim 15$.} For M3, we set parameters  $T=2, \gamma=2$. For M4, we selected the best result from the three parameter combinations of $T=6,\gamma=1$; $T=2,\gamma=2$ and $T=1.2,\gamma=4$ (the performance of the three parameter combinations for different functions is slightly different. We also tested other parameter combinations and found that no combination with a product of $T$ and $\gamma$ less than $4$ can achieve near machine precision).

Algorithms M1 and M2 exhibit computational complexity of  $\mathcal{O}(M\log M)$, in contrast to  the  $\mathcal{O}(M\log^2M)$ complexity inherent to algorithms M3 and M4.  Figure \ref{Figtime} demonstrates the scaling behavior of execution time versus node count $M$  across all four algorithms. These experimental observations align with the above conclusions: M1 and M2 demonstrate nearly identical temporal profiles, as do M3 and M4;  for equivalent
$M$ values, M1/M2 achieve superior computational efficiency compared with M3/M4. It should be noted that  $M$ values required by these four algorithms are different. For M1/M2, the minimum $M$ is governed by the approximation capacity of the boundary interval and   maximum frequency of the approximated function.
 For full-data algorithms such as M3 and M4, $M$ depends on both the function frequency and the critical parameter product $T\gamma$. Notably, our empirical studies and existing literature establish $T\gamma=4$ as the minimum allowable value  for full-data algorithms to attain machine precision. Subthreshold $T\gamma$ values result in  accuracy degradation for such methods.

\begin{figure}
			\begin{center}
\subfigure[\label{11a} $f(t)=\text{erf}(2t)$] {
\resizebox*{6cm}{!}{\includegraphics{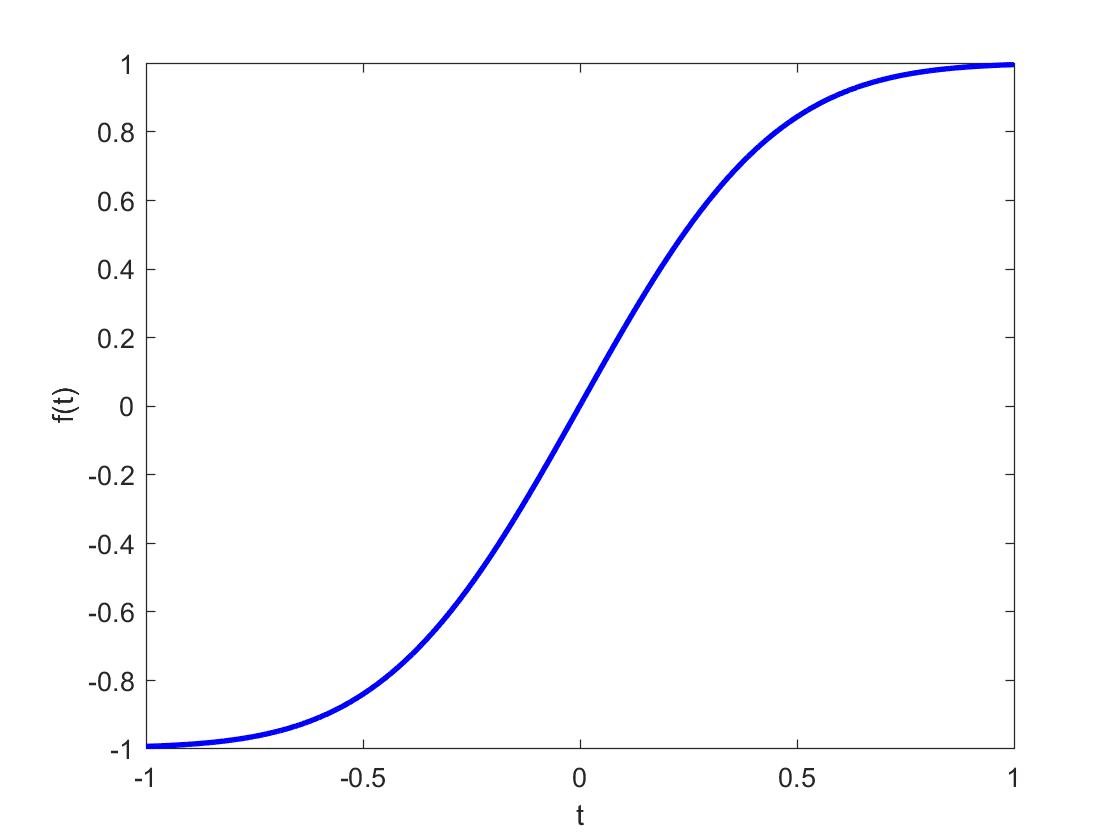}}
}%
\subfigure[\label{11b} approximation error] {
\resizebox*{6cm}{!}{\includegraphics{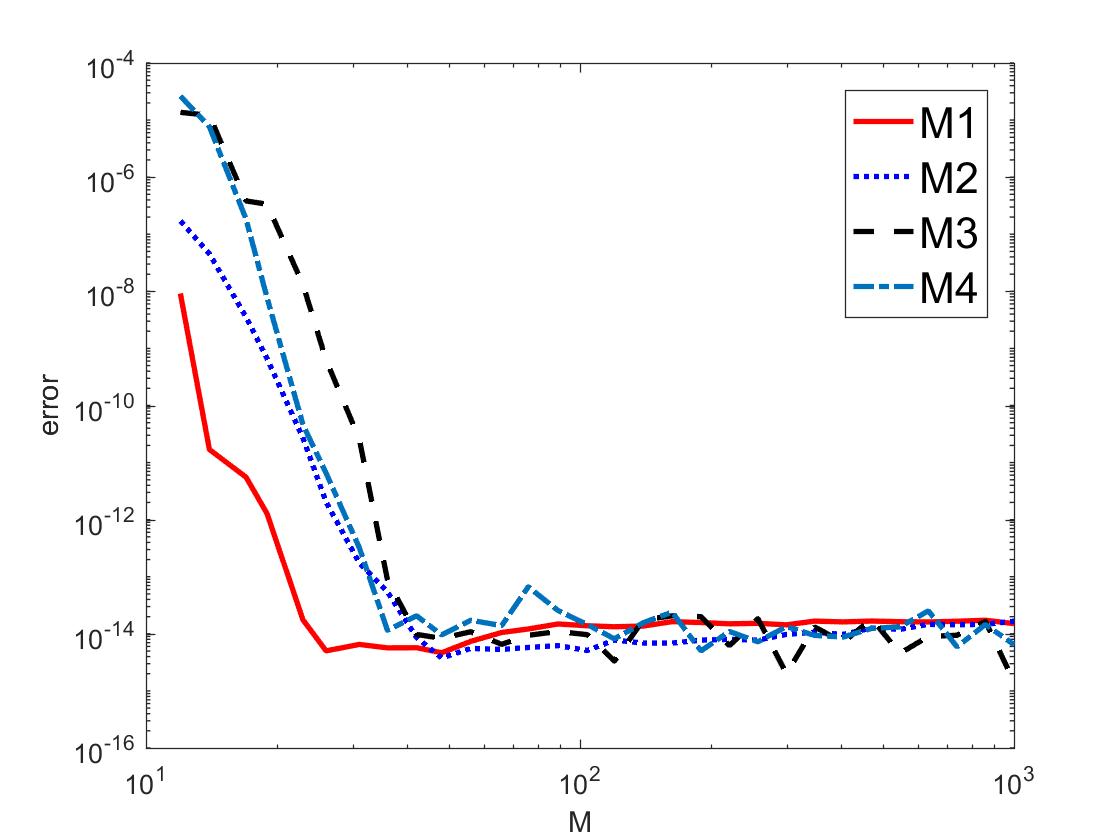}}
}%

\subfigure[\label{11c} $f(t)=\text{Ai}(1+3t)$] {
\resizebox*{6cm}{!}{\includegraphics{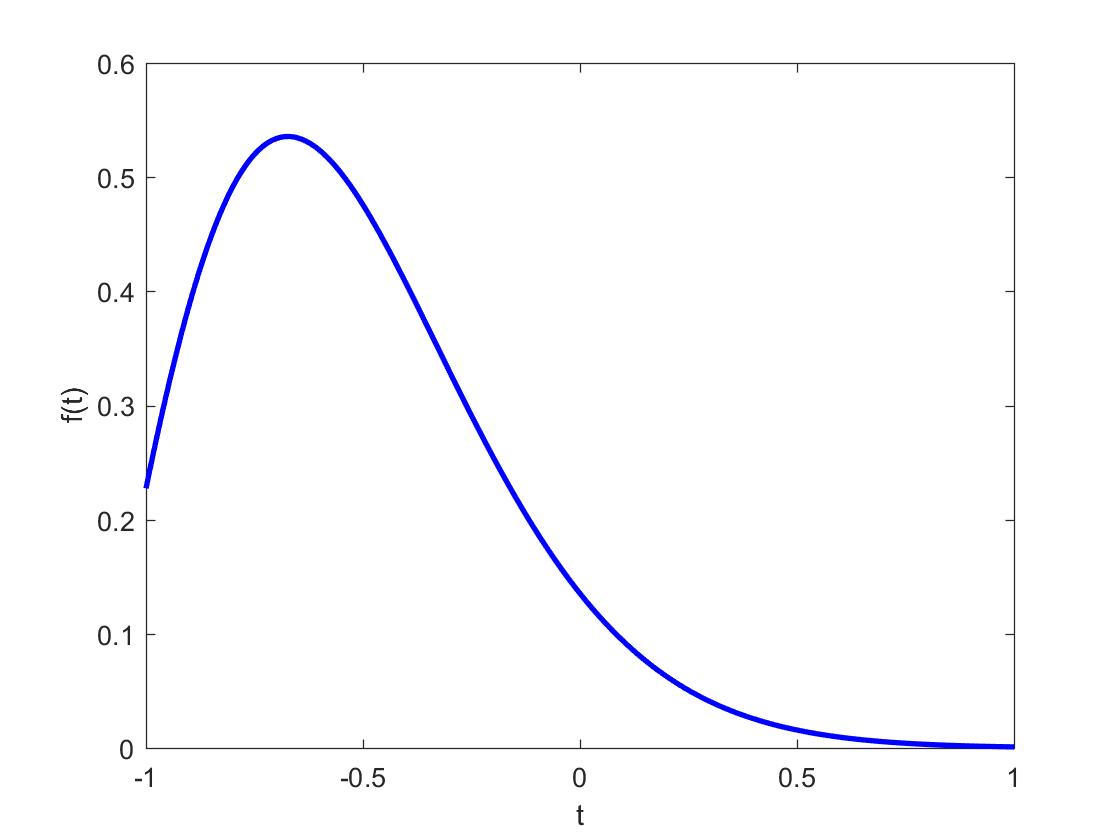}}
}%
\subfigure[\label{11d} approximation error] {
\resizebox*{6cm}{!}{\includegraphics{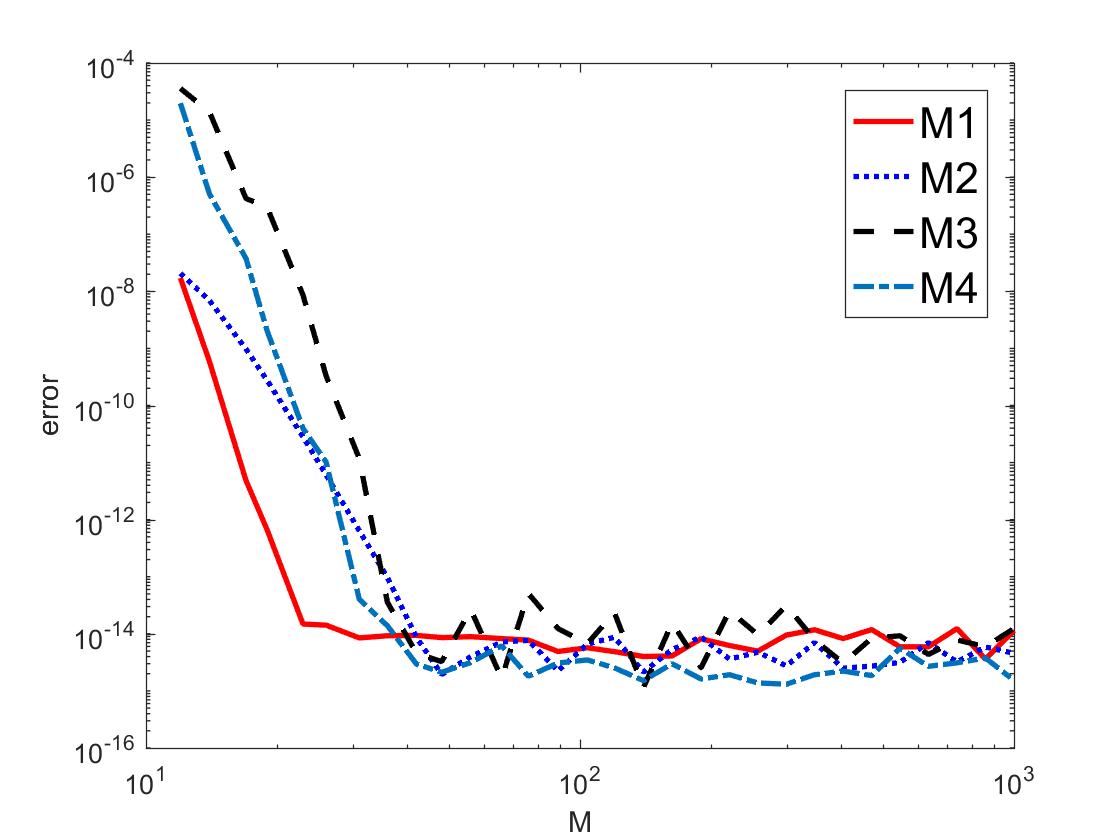}}
}%

\subfigure[\label{11e} $f(t)=\exp(\sin(2.7\pi t)+\cos(\pi t))$] {
\resizebox*{6cm}{!}{\includegraphics{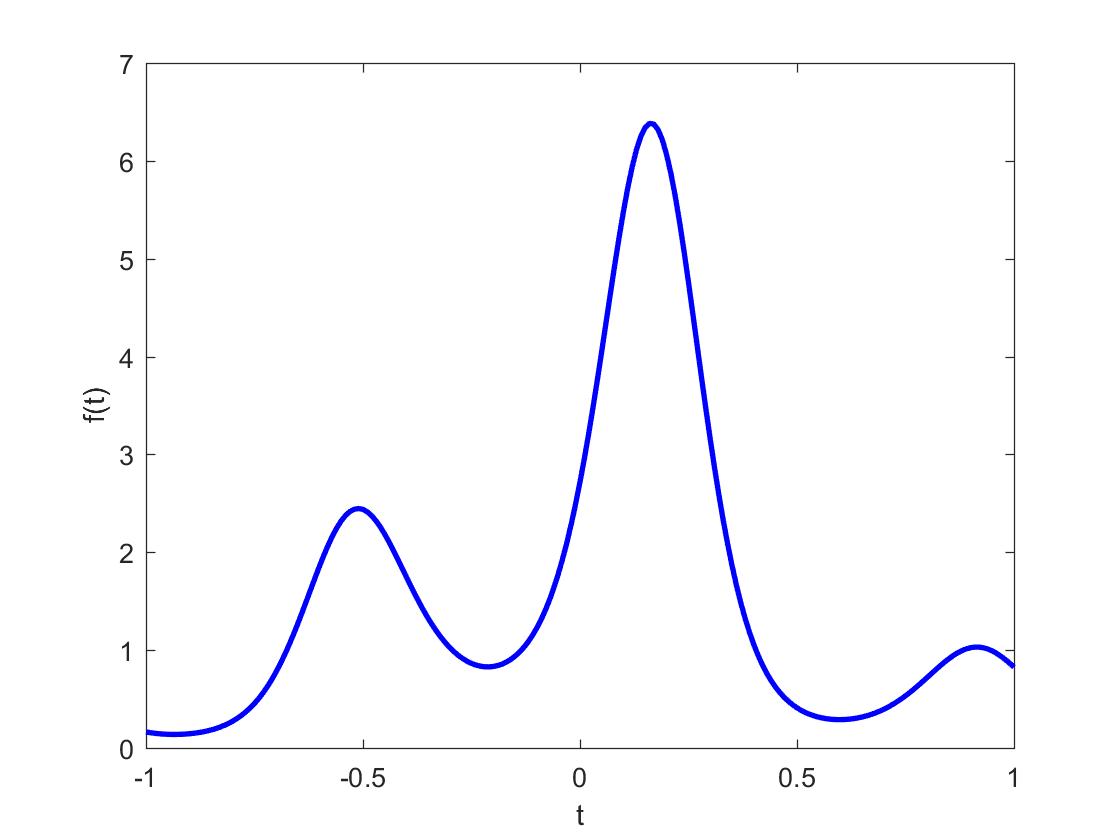}}
}%
\subfigure[\label{11f} approximation error] {
\resizebox*{6cm}{!}{\includegraphics{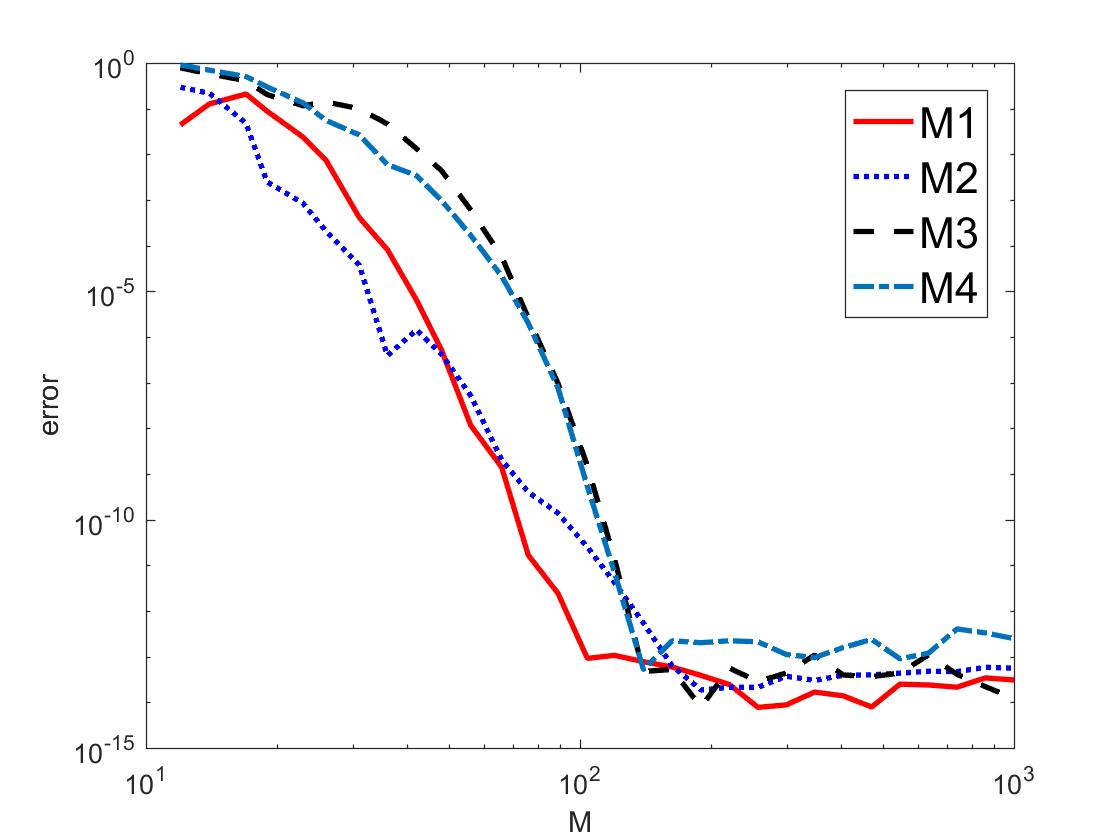}}
}%

	{\caption{ Overall smooth functions  and their approximation errors against $M$.}\label{Fig11}}
	\end{center}
\end{figure}

\begin{figure}
			\begin{center}
\subfigure[\label{12a} $f(t)=\frac{1}{1+100t^2}$] {
\resizebox*{6cm}{!}{\includegraphics{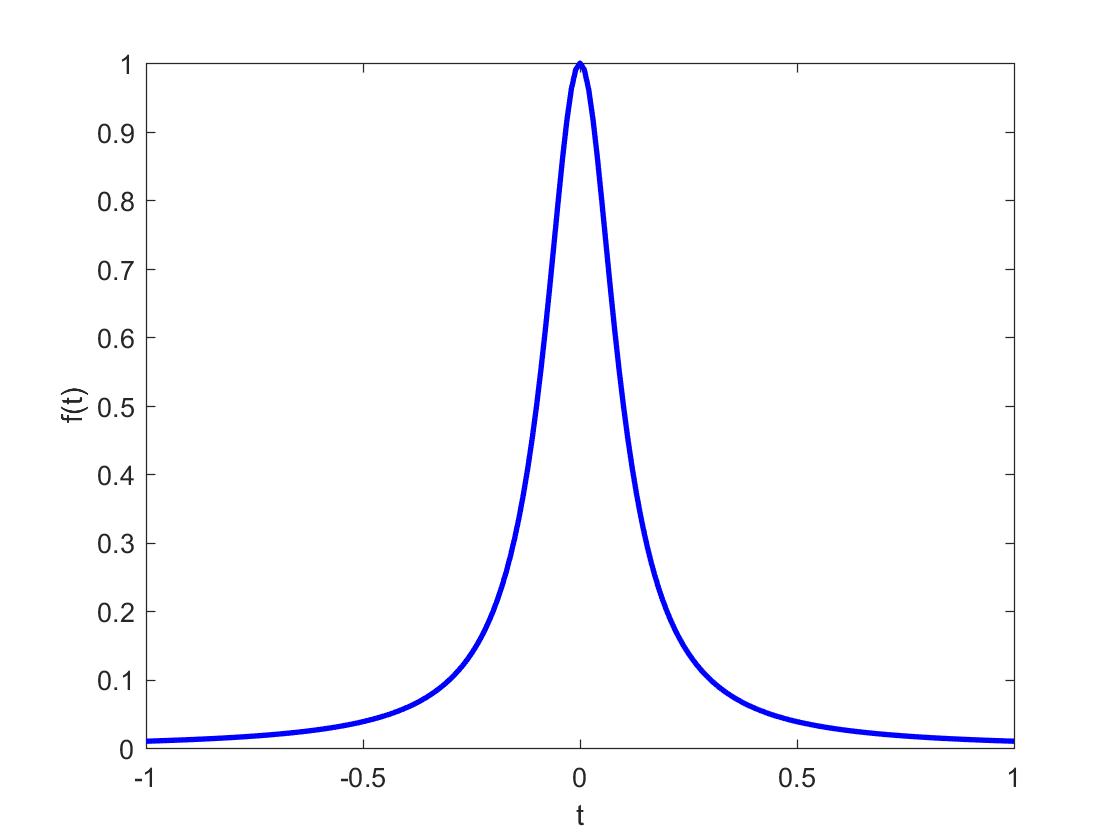}}
}%
\subfigure[\label{12b} approximation error] {
\resizebox*{6cm}{!}{\includegraphics{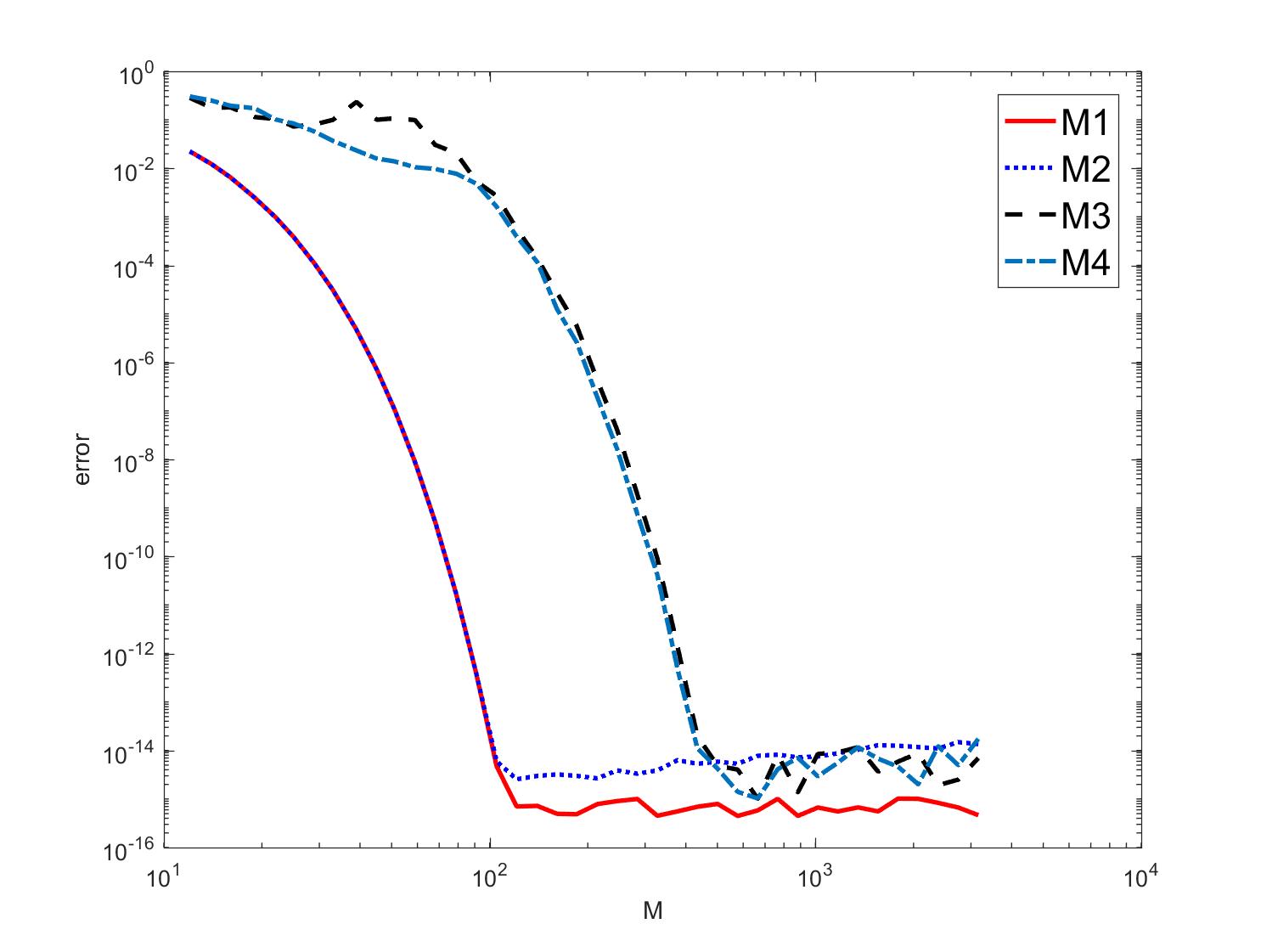}}
}%

\subfigure[\label{12c} $f(t)=\cos\left(\frac{100}{1+25t^2}\right)$] {
\resizebox*{6cm}{!}{\includegraphics{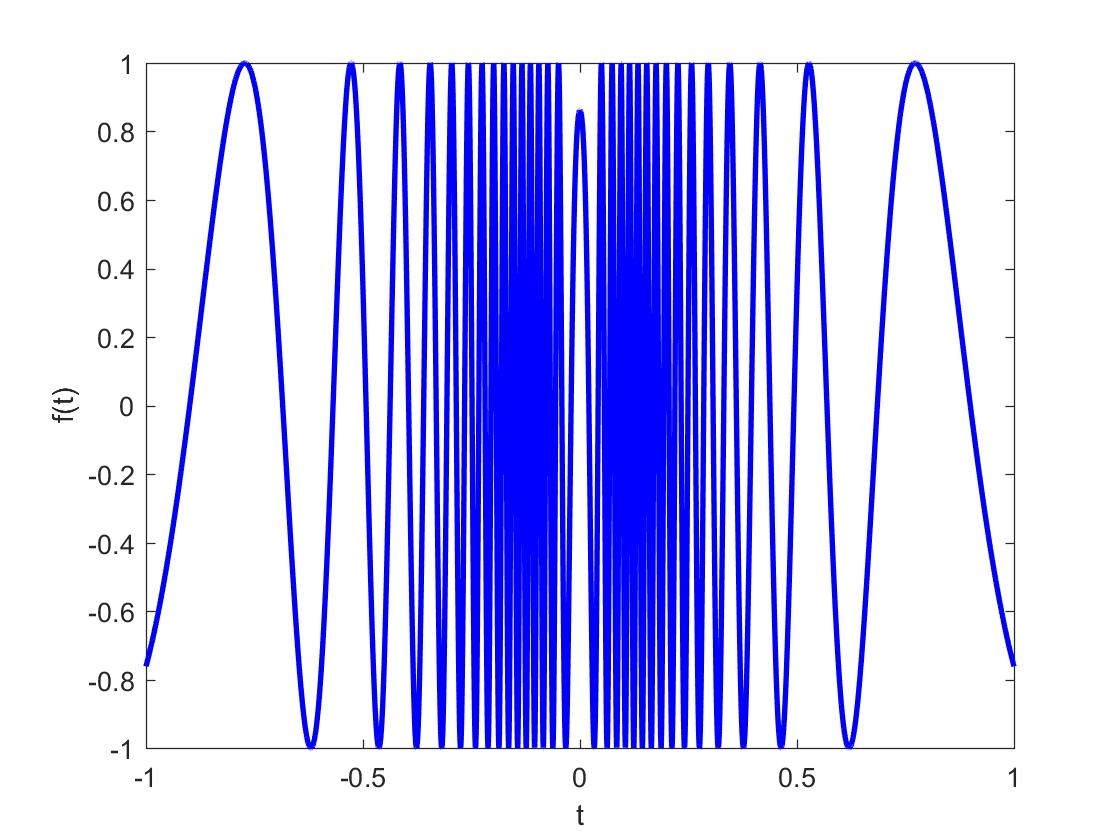}}
}%
\subfigure[\label{12d} approximation error] {
\resizebox*{6cm}{!}{\includegraphics{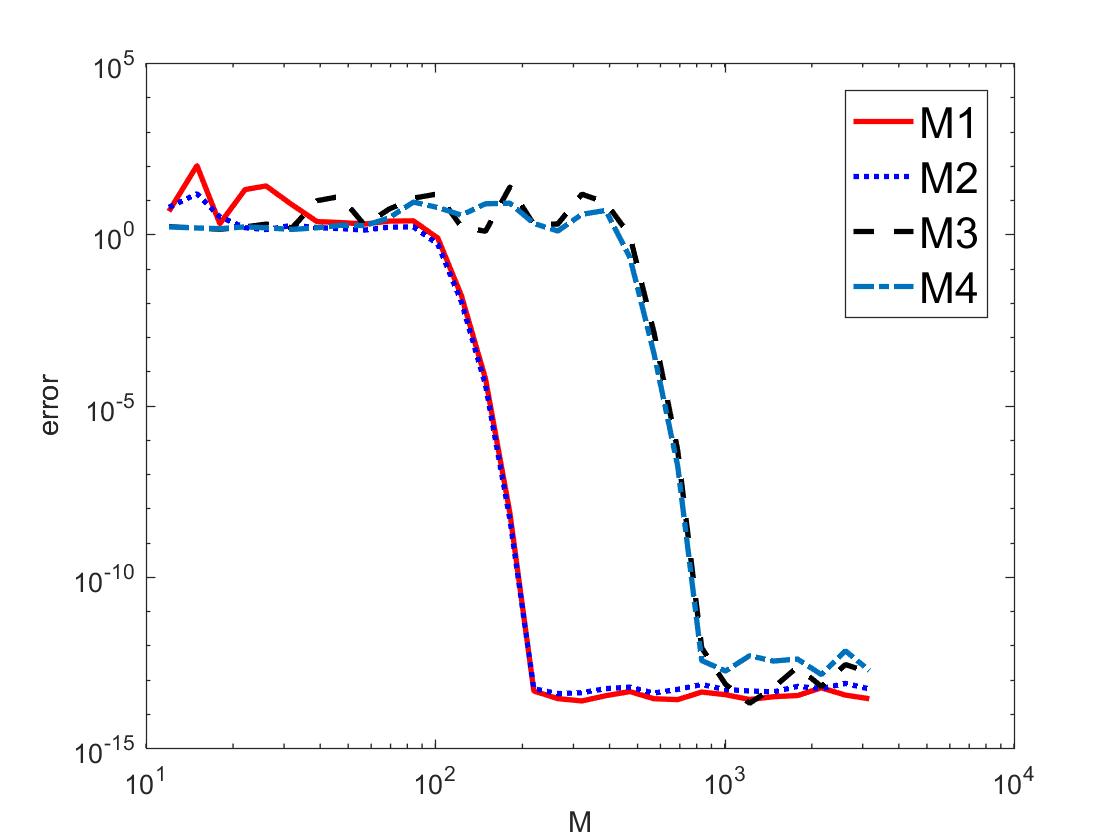}}
}%

\subfigure[\label{12e} $f(t)=\text{erf}(100t)$] {
\resizebox*{6cm}{!}{\includegraphics{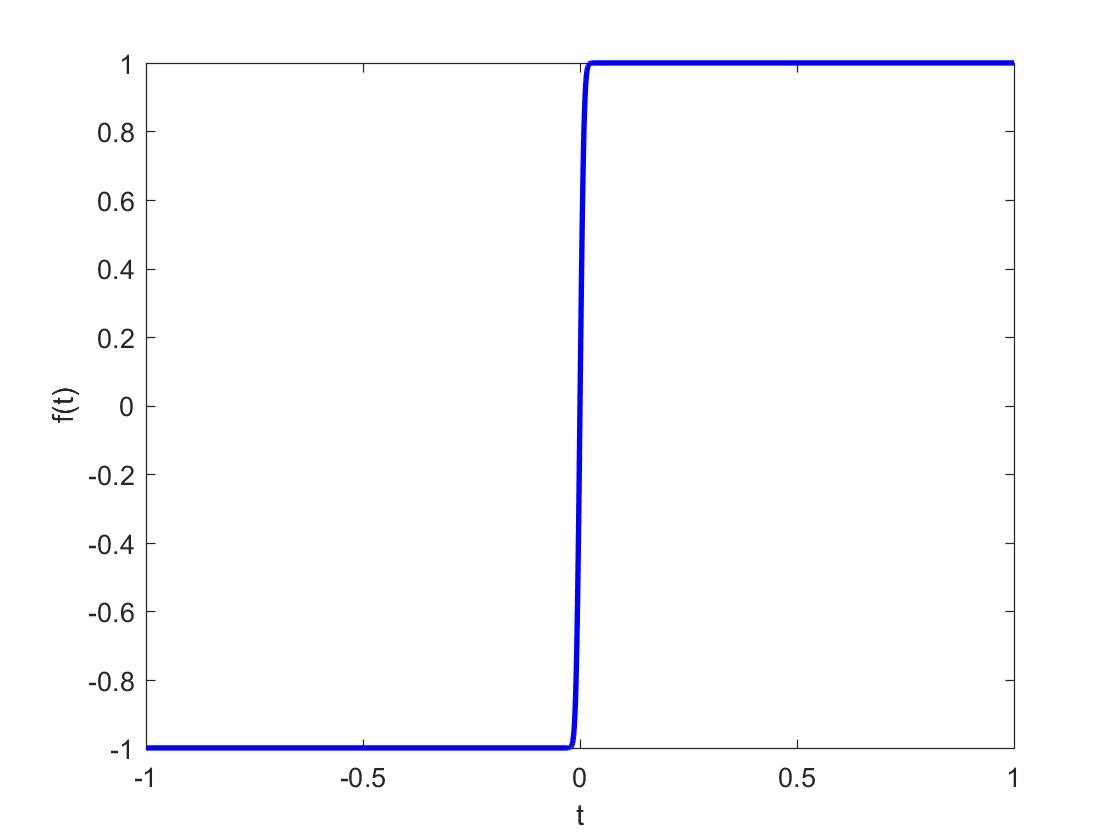}}
}%
\subfigure[\label{12f} approximation error] {
\resizebox*{6cm}{!}{\includegraphics{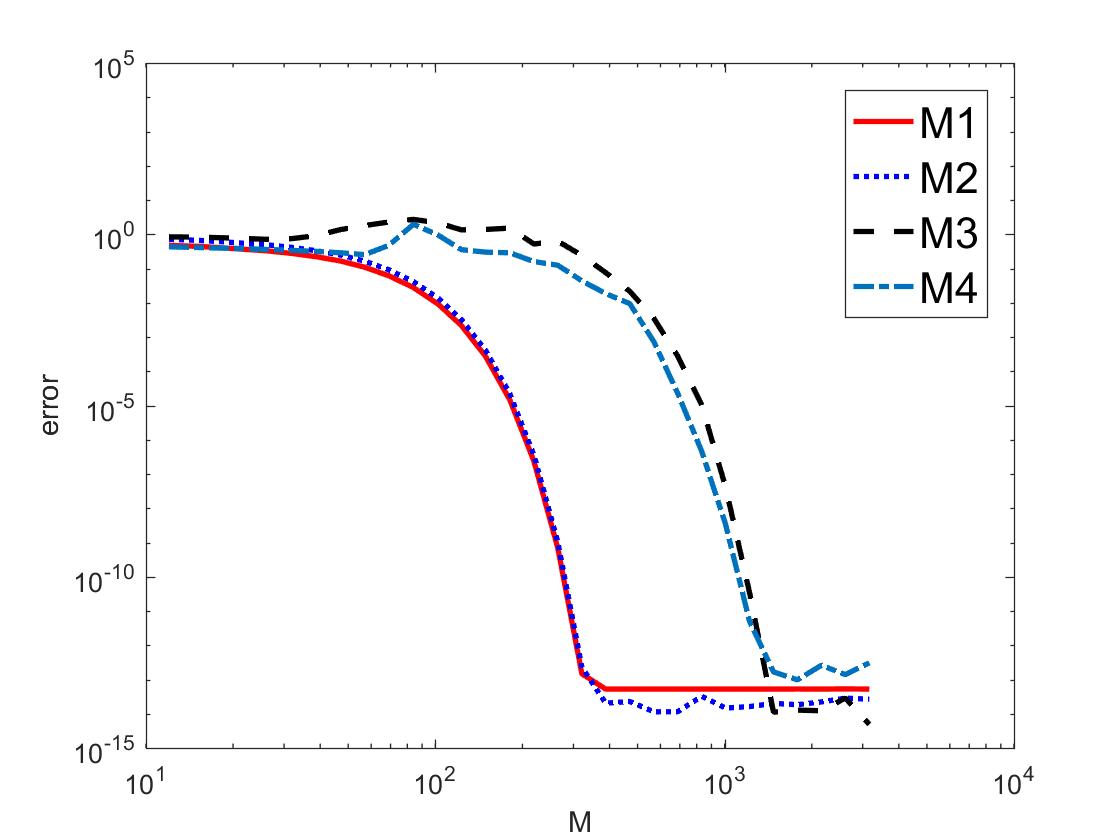}}
}%

	{\caption{ Internal oscillating functions  and their approximation errors against $M$.}\label{Fig12}}
	\end{center}
\end{figure}

\begin{figure}
			\begin{center}
\subfigure[\label{13a} $f(t)=\cos(100t^2)$] {
\resizebox*{6cm}{!}{\includegraphics{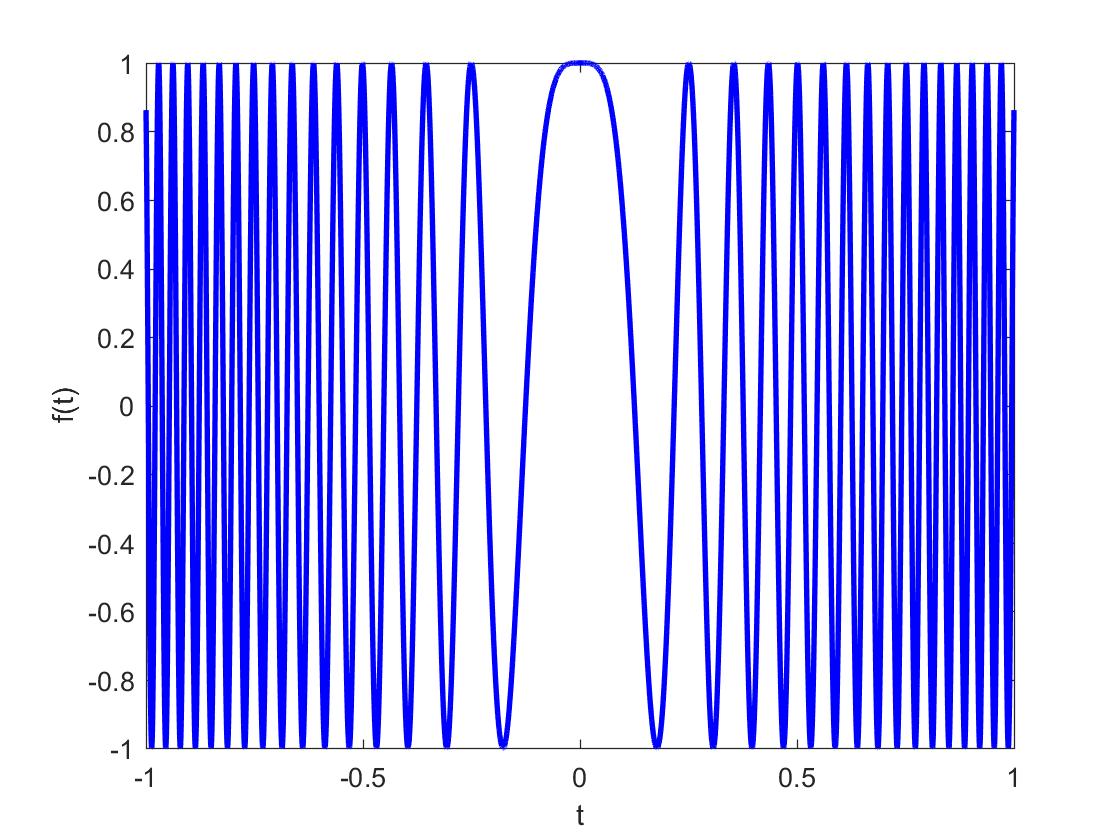}}
}%
\subfigure[\label{13b} approximation error] {
\resizebox*{6cm}{!}{\includegraphics{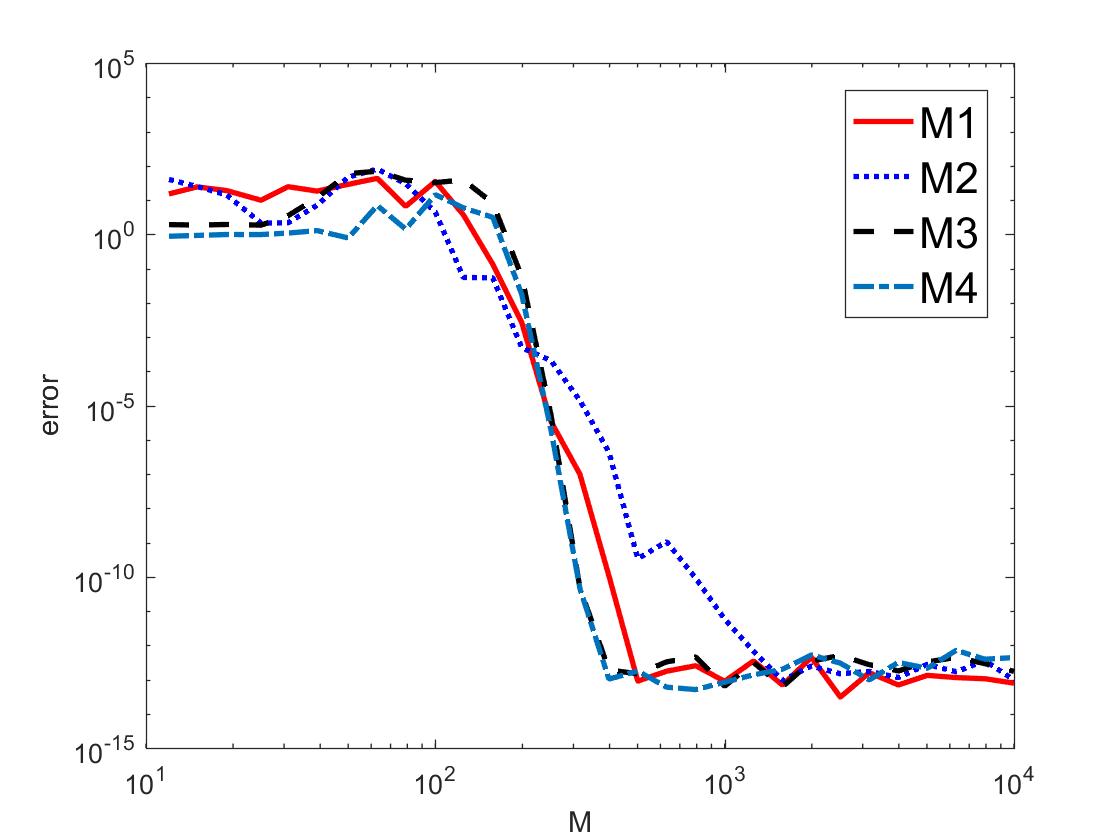}}
}%

\subfigure[\label{13c} $f(t)=\text{Ai}(-66-70t)$] {
\resizebox*{6cm}{!}{\includegraphics{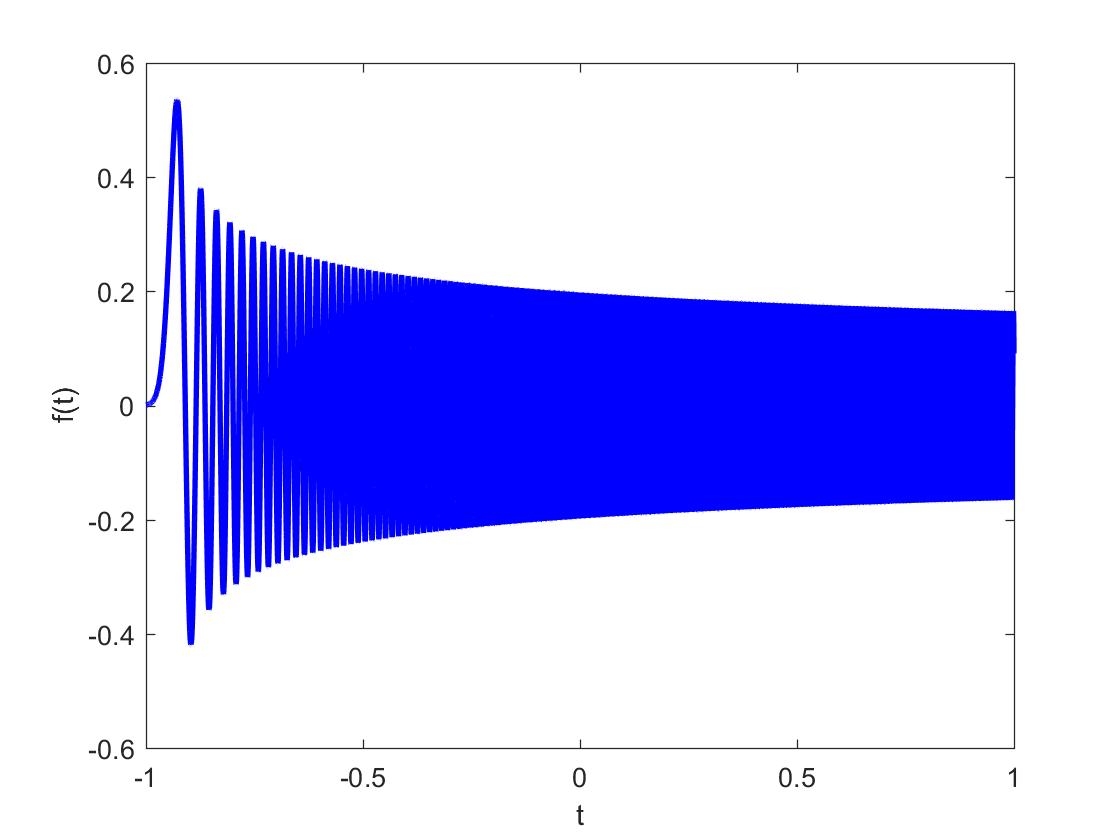}}
}%
\subfigure[\label{13d} approximation error] {
\resizebox*{6cm}{!}{\includegraphics{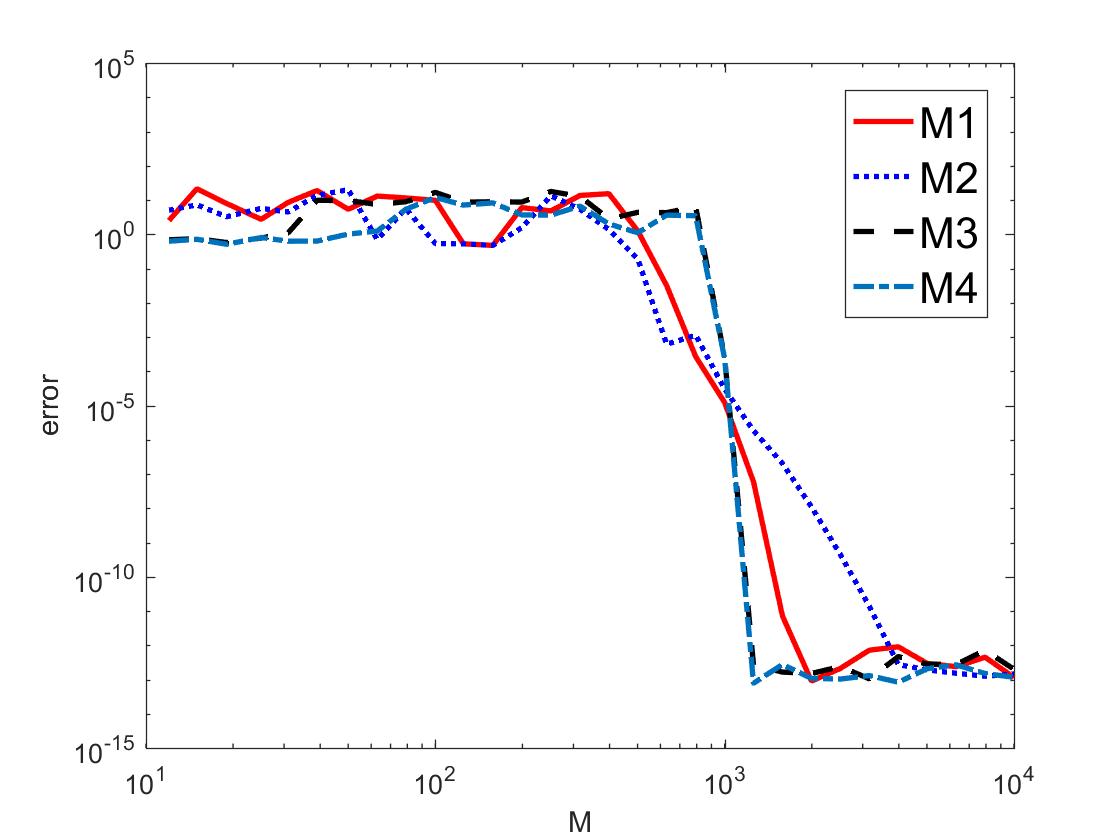}}
}%

\subfigure[\label{13e} $f(t)=\exp(\sin(65.5\pi t-27\pi)-\cos(20.6\pi t))$] {
\resizebox*{6cm}{!}{\includegraphics{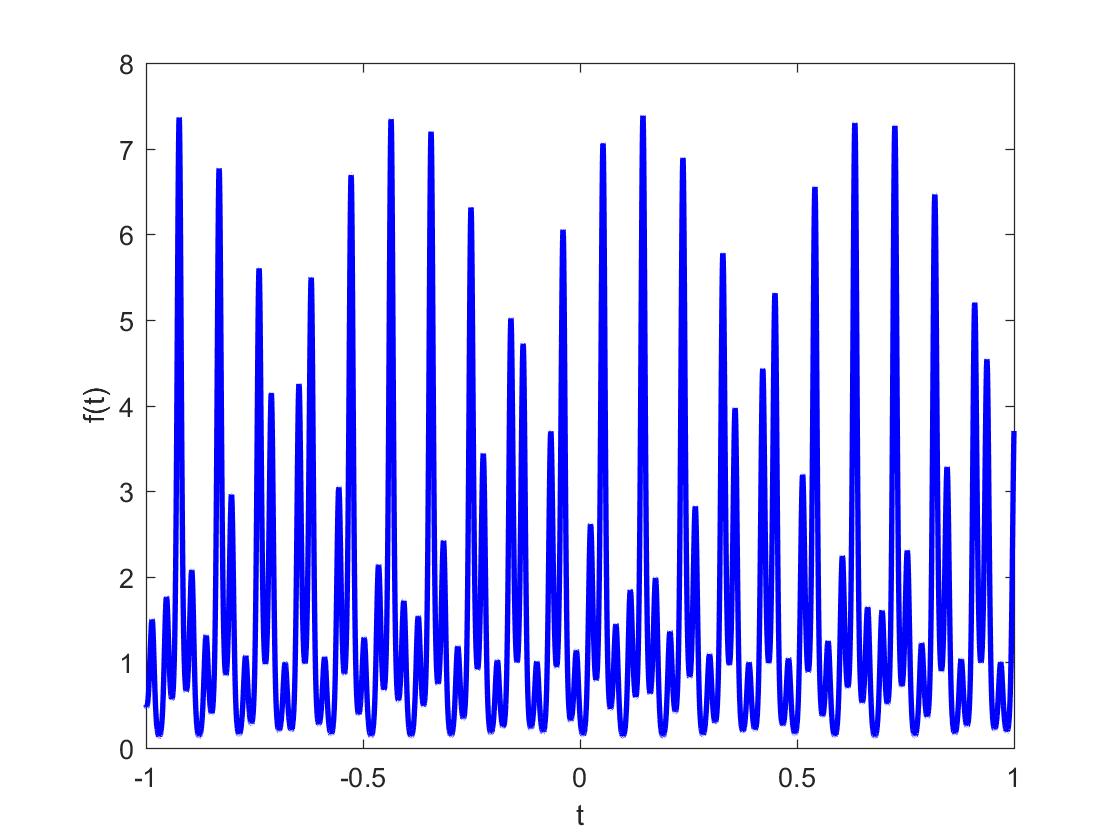}}
}%
\subfigure[\label{13f} approximation error] {
\resizebox*{6cm}{!}{\includegraphics{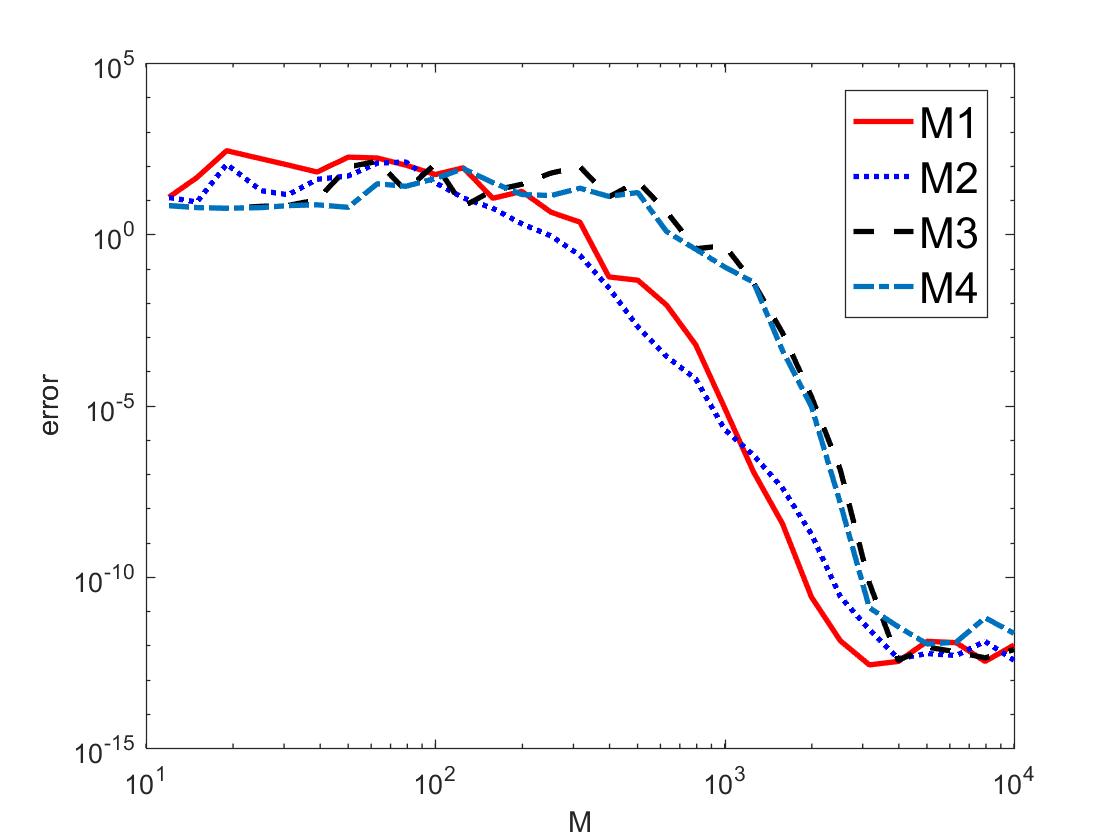}}
}%

	{\caption{ Boundary oscillating functions  and their approximation errors against $M$.}\label{Fig13}}
	\end{center}
\end{figure}

The accuracy of the solution obtained using the four algorithms is shown in Fig. \ref{Fig11} - \ref{Fig13}. To better show the differences, we divided the function into three cases for comparison:
\begin{itemize}
  \item Overall smooth functions:
\begin{equation*}
  f_1(t)=\text{erf}(2t),\quad f_2(t)=\text{Ai}(1+3t),\quad f_3(t)=e^{\sin(2.7\pi t)+\cos(\pi t)}.
\end{equation*}
  \item Internal oscillating functions:
\begin{equation*}
  f_4(t)=\frac{1}{1+100t^2},\quad f_5(t)=\frac{100}{\cos(1+25t^2)},\quad f_6(t)=\text{erf}(100t).
\end{equation*}
  \item  Boundary oscillating functions:
\begin{equation*}
\begin{aligned}
  &f_7(t)=\cos(100t^2),\quad f_8(t)=\text{Ai}(-66-70t),\\
  & f_9(t)=e^{\sin(65.5\pi t-27\pi)-\cos(20.6\pi t)}.
  \end{aligned}
\end{equation*}
\end{itemize}

 { Fig. \ref{Fig11} shows the numerical results for the first case. As mentioned above, the number of nodes $M$ required for M1 and M2 is determined by the frequency of the boundary interval, but the re-solution ability of polynomial approximation is weaker than that of trigonometric approximation, so the value of $M$ required by M2 is  larger than M1; while M3 and M4 are determined by the maximum frequency of the function in the entire interval, and the number of nodes required is  slightly larger than M1. At this time, the number of nodes $M$ required by the four algorithms to achieve machine accuracy was relatively small, therefore, there was no obvious difference in performance.

In the second case, that is, when the function oscillates internally. As the frequency of the function in the boundary interval is much lower than that in the interior, the number of nodes required for M1 and M2 is determined by the resolution constant of the standard Fourier transform. Note that the resolution constants of M3 and M4 are $T\gamma$, so from \eqref{estrrs4} we can see that the number of nodes required for M1 and M2 is approximately $1/4$ of that for M3 and M4, which is consistent with the results given in  Fig. \ref{Fig12}.

In the third case, because the function oscillates near the boundary, the algorithm in this paper cannot avoid this singularity. For functions $f_7$ and $f_8$, the maximum frequency of the function appears in the boundary interval. According to \eqref{estrs2}, the number of nodes required by the proposed algorithm is $T_{\Delta}\gamma_{\Delta}=6$ times the maximum frequency, whereas the number of nodes required by M3 and M4 is $T\gamma=4$ times the maximum frequency. This is consistent with the results presented in Figures 2 and 3. For $f_9$, because the internal frequency is higher, the number of nodes required by the proposed algorithm is less than that of M3 and M4. As for M2, we can see that its approximation ability for the boundary interval data is obviously weaker than that of M1, therefore, more nodes are needed to achieve { near-machine} accuracy.}

{ We observe that in certain cases (e.g., Fig. \ref{11f}), Algorithm M2 exhibits superior accuracy over M1 for small node counts $M$. Since M2 employs fewer boundary nodes, its captured boundary intervals inherently exhibit smoother relative behavior compared to M1 under the same conditions. Notably, reducing the number of boundary nodes in M1 can similarly enhance early-stage accuracy; however, as demonstrated in prior tests (Section 5.2), such parameter configurations compromise M1's ability to achieve near-machine precision with sufficient nodes.}

\subsection{Algorithm improvement: boundary grid refinement}
In certain  situations, different step sizes may be used for interior and boundary sampling.
The test results in the previous subsection indicate that the resolution constant is affected by the oscillation of the boundary interval data. This naturally raises the question: Can the algorithm performance be improved if a finer grid is used near the boundary than inside? In this section, we   combine this concept to  provide a modified algorithm for the boundary oscillation function.

Let $R$ be an positive  integer and
\begin{equation}
  {t}^{R}_{\ell}=\frac{\ell}{RM}, \quad \ell=-RM,\ldots, RM.
\end{equation}
The positive integer $m^{R}_{\Delta}$ denotes the number of nodes selected in the boundary interval. Notably, only $2\times m^{R}_{\Delta}$ nodes near the boundary were used in the actual calculation. The other nodes can be ignored, but they are listed here for  convenience.
For consistency with the boundary interval taken in the previous algorithm, we let
\begin{equation}
  {m}^{R}_{\Delta}=R(m_{\Delta}-1)+1
\end{equation}
and
\begin{equation*}
\begin{aligned}
 & L_{\Delta}^R=2\times \left\lceil T_{\Delta}\times (m^{R}_{\Delta}-1)\right\rceil,\\& h^{R}=\frac{2\pi}{L^{R}_{\Delta}},\quad x_j^R=(j-1)h, \quad j=1,2,\ldots,L^{R}_{\Delta}.
 \end{aligned}
\end{equation*}
Similarly, let
\begin{equation*}
\begin{aligned}
  &J^R_{\Delta,1}=\{1,2,\ldots,m^R_{\Delta}\},~ J^R_{\Delta,2}=\{\frac{L^R_{\Delta}}{2}+1,\frac{L^R_{\Delta}}{2}+2,\ldots,\frac{L^R_{\Delta}}{2}+m^R_{\Delta}\}, ~ \\&J^R_{\Delta}=J^R_{\Delta,1}\cup J^R_{\Delta,2}.
  \end{aligned}
\end{equation*}
and for $j\in J^{R}_{\Delta}$, take
\begin{equation}
  g^R(x^R_j)=\left\{
  \begin{aligned}
    &f(t^R_{RM-m^R_{\Delta}+j}),& j\in J_{\Delta,1},\\
    &f(t^R_{-RM+j-{L^R_{\Delta}}/{2}-1}),&j\in J_{\Delta,2}.
  \end{aligned}\right.
\end{equation}
The extension function  $g^R_c(x^R_j), j=1,2,\ldots, L^R_{\Delta}$ was also obtained according to the method described in section \ref{SEC3}.
And then we  obtain a function $f^R_c$ with a period of $2+\lambda$, where $\lambda$ is same as in \eqref{deflambda}. The data of the extension function $f^{R}_c(t_{\ell}),\ell=-M,\ldots,M, M+1,\ldots, M+\lceil T_{\Delta}-1\rceil\times (m_{\Delta}-1)-1$ in one period $[-1,1+\lambda)$ are
\begin{equation}
  f^R_c(t_{\ell})=\left\{
  \begin{aligned}
    &f(t_{\ell}),&|\ell|\leq M,\\
&g^R_c(x_{m^{R}_{\Delta}+(\ell-M)\cdot R}),&\ell>M.
  \end{aligned}\right.
\end{equation}
The original algorithm is equivalent to a special case of the improved algorithm when the parameter $R=1$.

\begin{figure}[H]
	\begin{center}
\subfigure[\label{R1} the approximation error against $R$] {
\resizebox*{6cm}{!}{\includegraphics{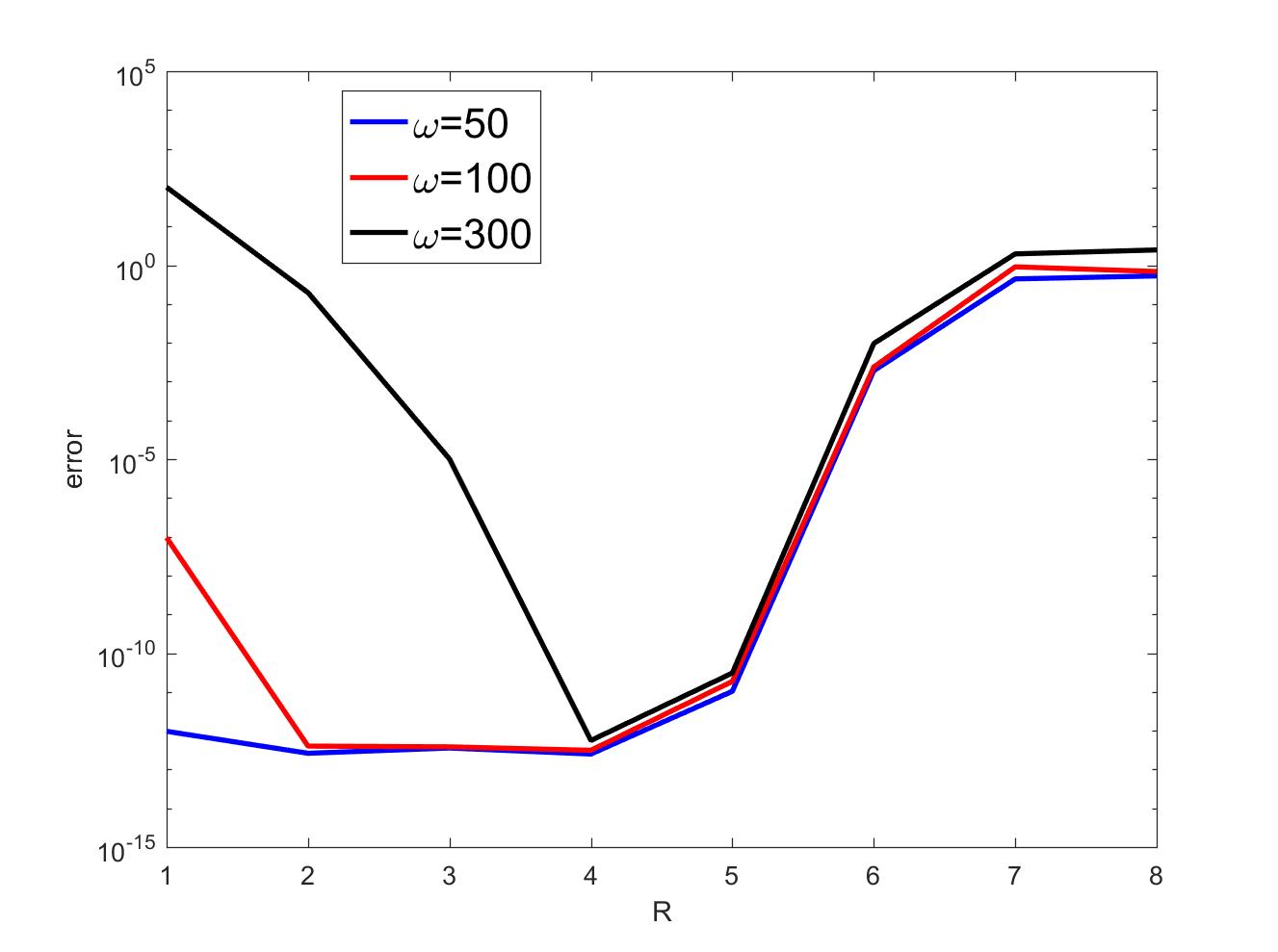}}
}{ \subfigure[\label{R2} the approximation error against $t$  $(\omega=300)$] {
\resizebox*{6cm}{!}{\includegraphics{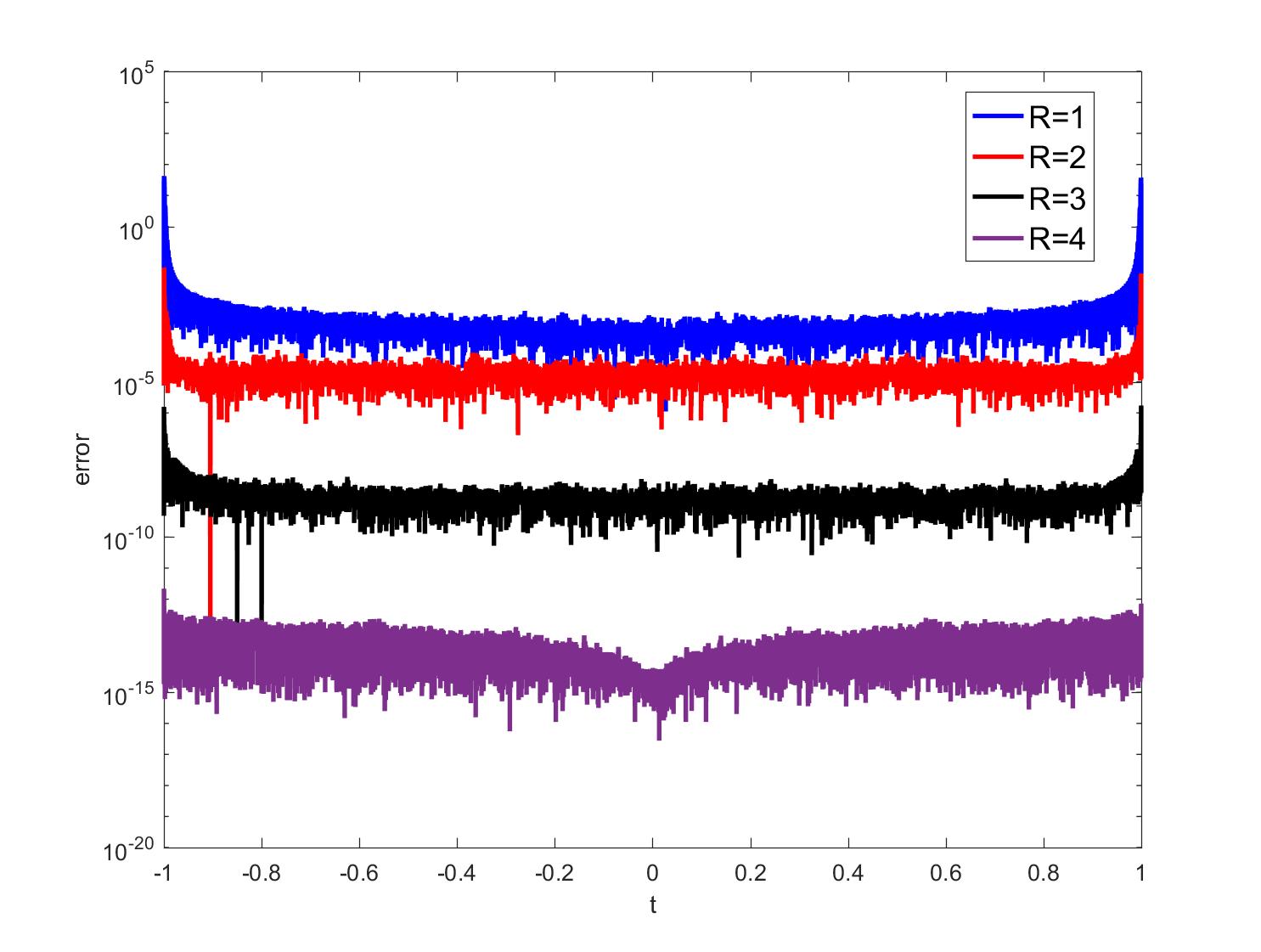}}
}}%
		{\caption{The approximation error with various $R$ for functions $\exp(\text{i}\omega \pi t)$, $M=500$. \label{testR}}}
	\end{center}
\end{figure}
\begin{figure}
			\begin{center}
\subfigure[\label{14a} $f(t)=\frac{1}{1.01-t^2}$] {
\resizebox*{6cm}{!}{\includegraphics{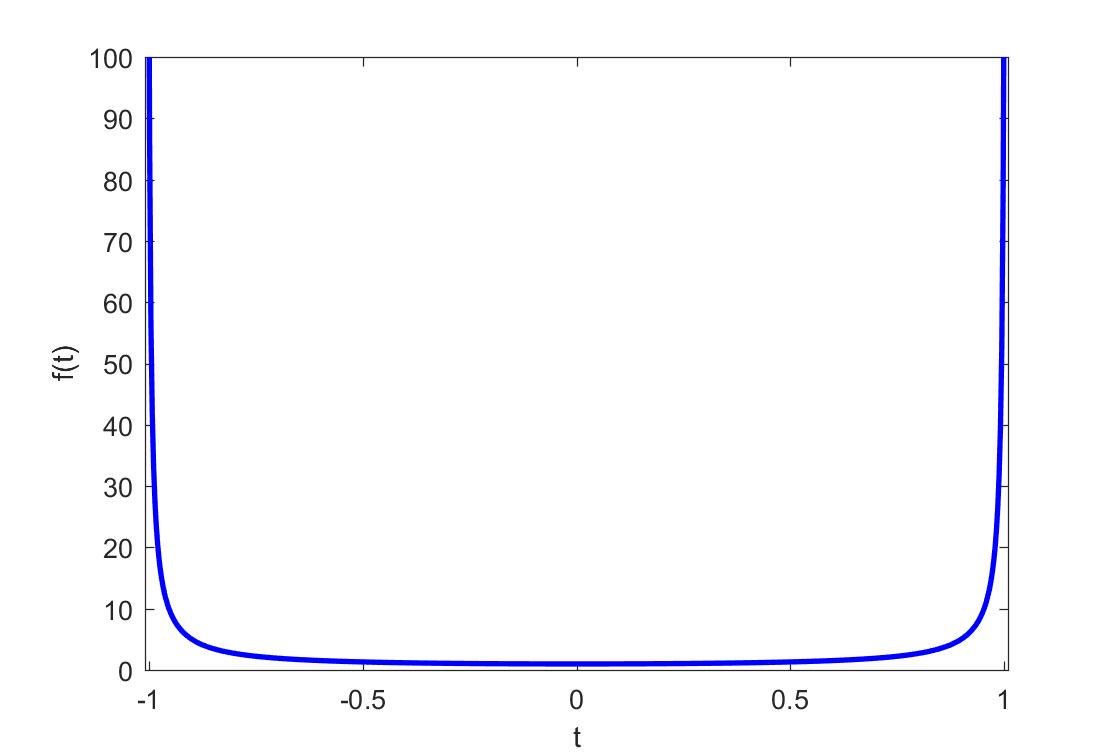}}
}%
\subfigure[\label{14b} approximation error] {
\resizebox*{6cm}{!}{\includegraphics{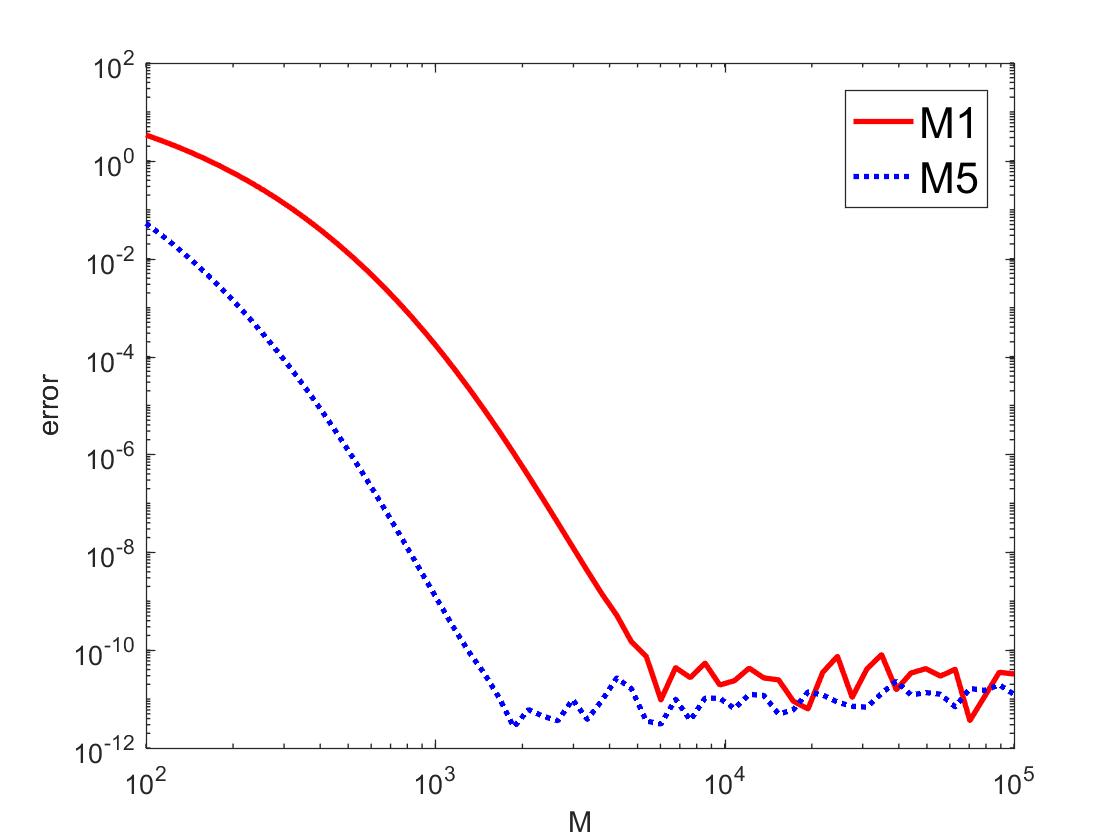}}
}%

\subfigure[\label{14c} $f(t)=\text{Ai}(150t)$] {
\resizebox*{6cm}{!}{\includegraphics{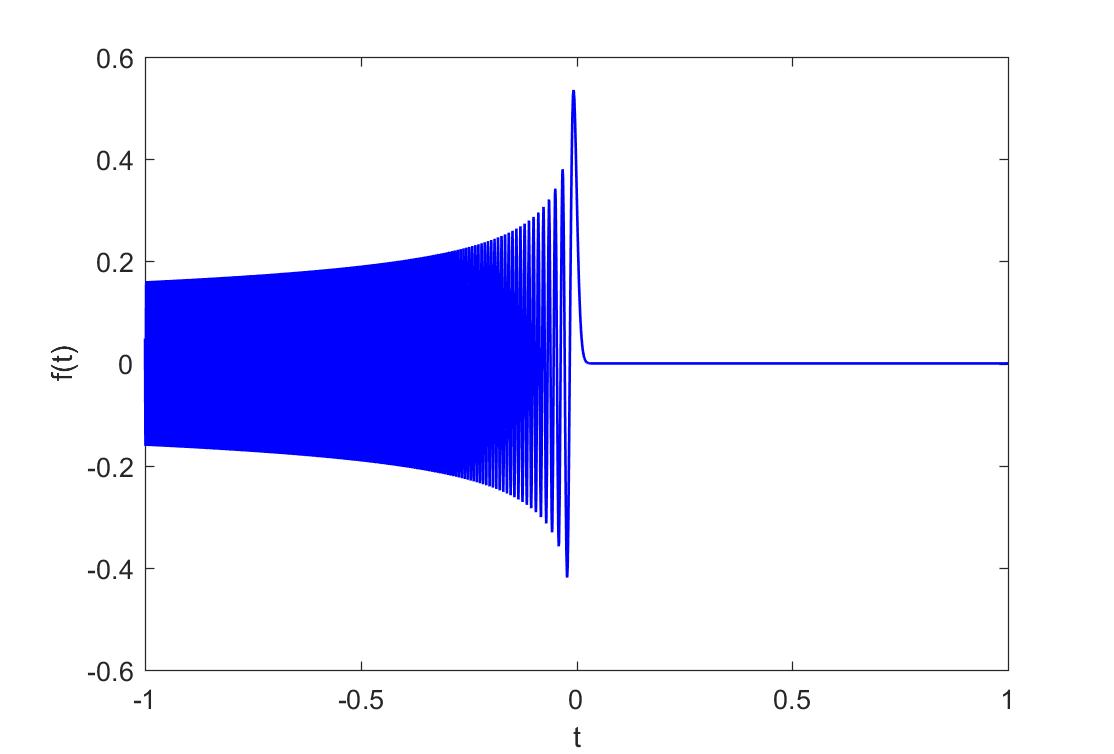}}
}%
\subfigure[\label{14d} approximation error] {
\resizebox*{6cm}{!}{\includegraphics{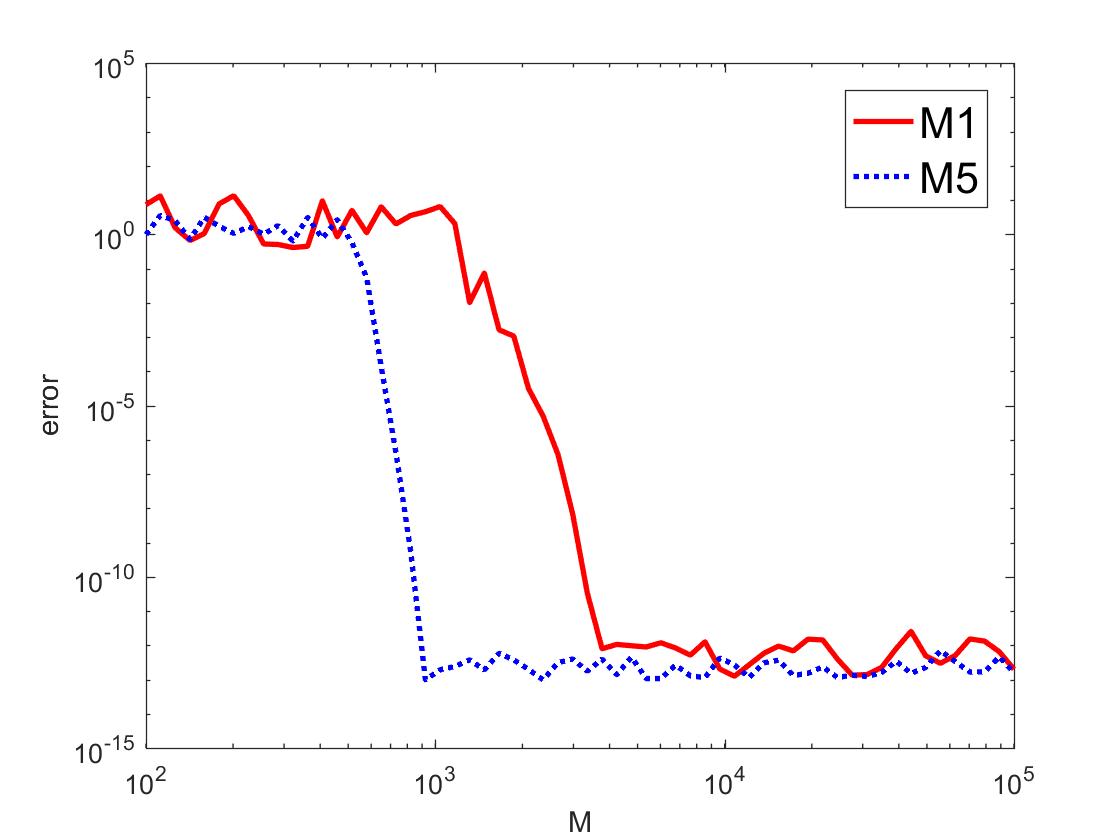}}
}%

\subfigure[\label{14e} $f(t)=\sin(1500t^2)$] {
\resizebox*{6cm}{!}{\includegraphics{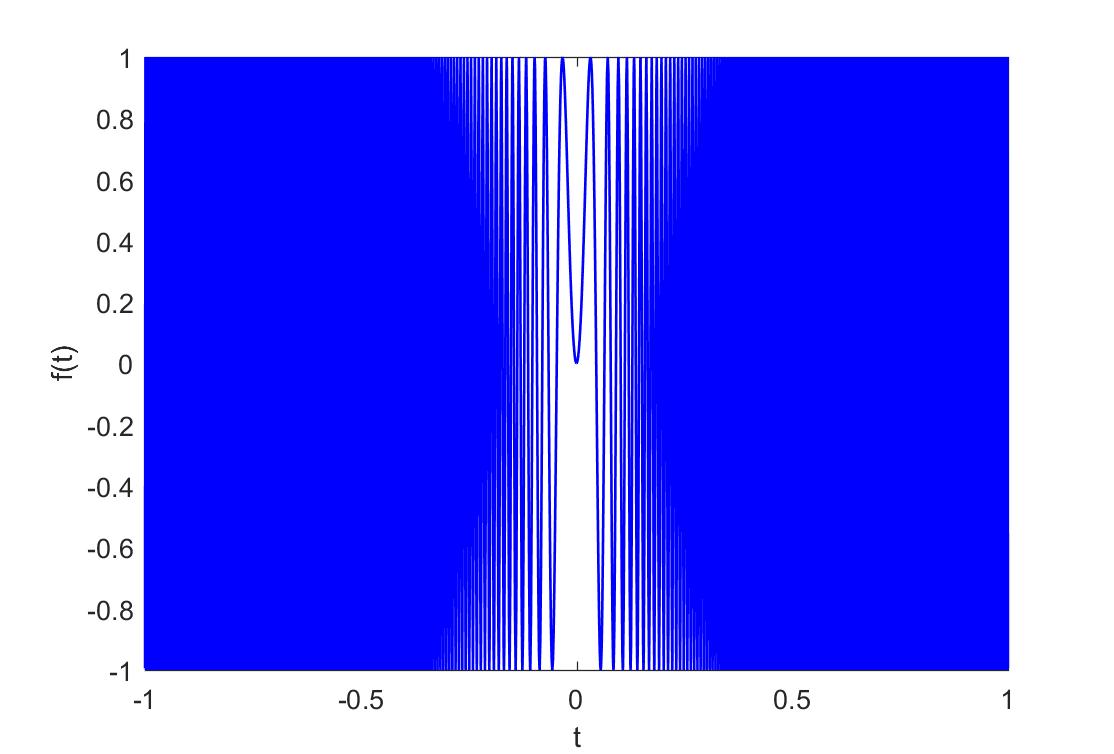}}
}%
\subfigure[\label{14f} approximation error] {
\resizebox*{6cm}{!}{\includegraphics{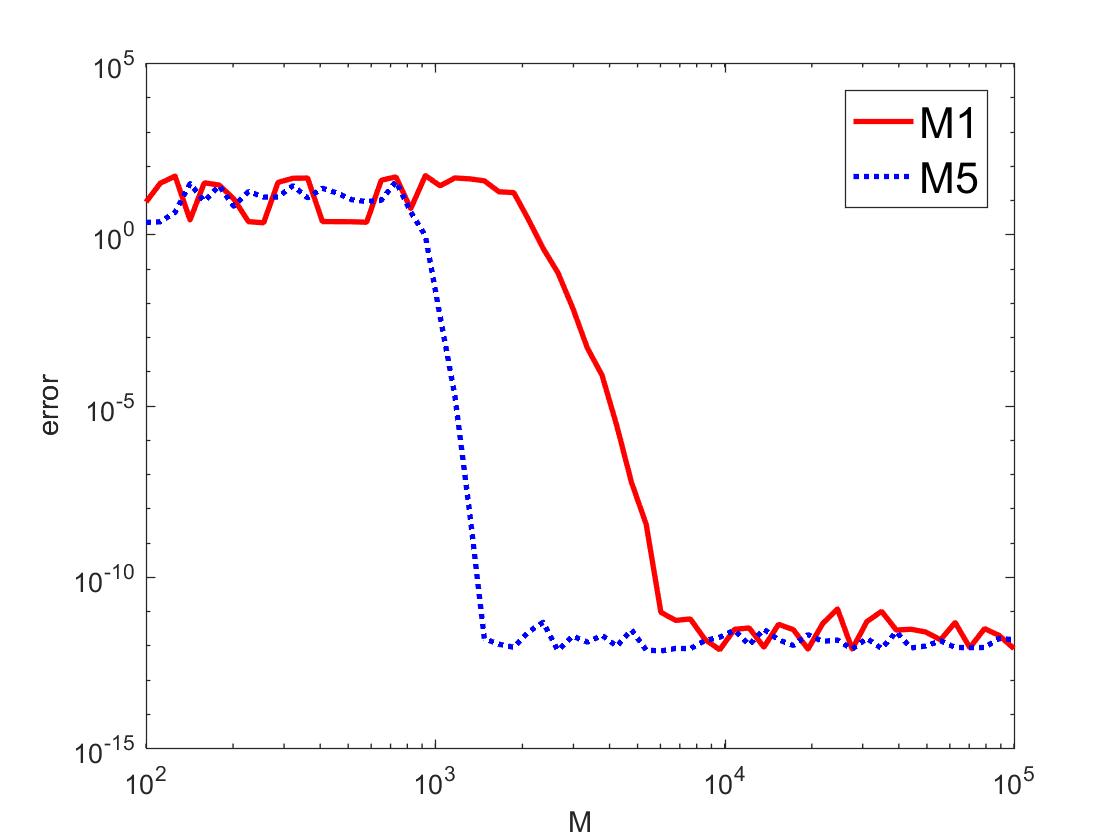}}
}%

	{\caption{ Comparisons of algorithms M1 and M5.}\label{Fig14}}
	\end{center}
\end{figure}
Now we take $\gamma=1, T_{\Delta}=6, m_{\Delta}=25$ to test the influence of parameter $R$ on the algorithm. As can be seen from Fig. \ref{R1}, when $R\geq5$, the approximation error not only does not improve, but becomes worse, which is caused by the instability of the calculation of the extension part. When $R\leq4$, it can be observed that the refinement of grid near the boundary improves the approximation error. When $\omega=50$, the improvement is not obvious, because the number of nodes $M=500>T_{\Delta}\gamma_{\Delta}\omega$ has exceeded the threshold for resolving this frequency function. Therefore, when facing high-frequency functions, we recommend using parameters
$$\gamma_{\Delta}=1, T_{\Delta}=6, m_{\Delta}=25, R=4$$
for calculation. {It can be seen from Fig. \ref{R2} demonstrates that progressive increases in parameter
$R$ from $1$ to $4$ induce uniform error reduction across the entire computational domain. As the boundary grid is refined, the reduction of the boundary error improves the overall approximation error.}

Next, we tested the performance of the improved algorithm (M5) and compared it with the original algorithm (M1). We selected several test functions with significant high-frequency features in the boundary interval:
\begin{equation*}
\begin{aligned}
 f_{10}(t)=\frac{1}{1.01-t^2},& &f_{11}(t)=\text{Ai}(150t),&&f_{12}(t)=\sin(1500t^2).
\end{aligned}
 \end{equation*}

Fig. \ref{Fig14} shows the graphes of the test functions and the corresponding approximation results. We can see that the refinement of grid near the boundary effectively reduces the resolution constant. The number of nodes required to achieve machine accuracy in the improved method was   approximately $1/4$ of that in the original method.

{Next, we compare the performance of the M5 method with Chebyshev interpolation. Polynomial interpolation at Chebyshev nodes is a widely used approach for function approximation, renowned for its high accuracy and favorable resolution constant. Figure \ref{Figcheby} shows the approximation errors for both methods applied to the test functions $f_3$, $f_5$, $f_7$, and $f_9$ as the number of interpolation nodes $M$ increases. The relative performance of the two methods depends primarily on the characteristics of the function, especially its frequency distribution properties. Overall, the performance difference between the methods is modest. Furthermore, both techniques are compatible with FFT-based implementations, ensuring high computational efficiency.

\begin{figure}
			\begin{center}
\subfigure[\label{15a} $f_3(t)$] {
\resizebox*{6cm}{!}{\includegraphics{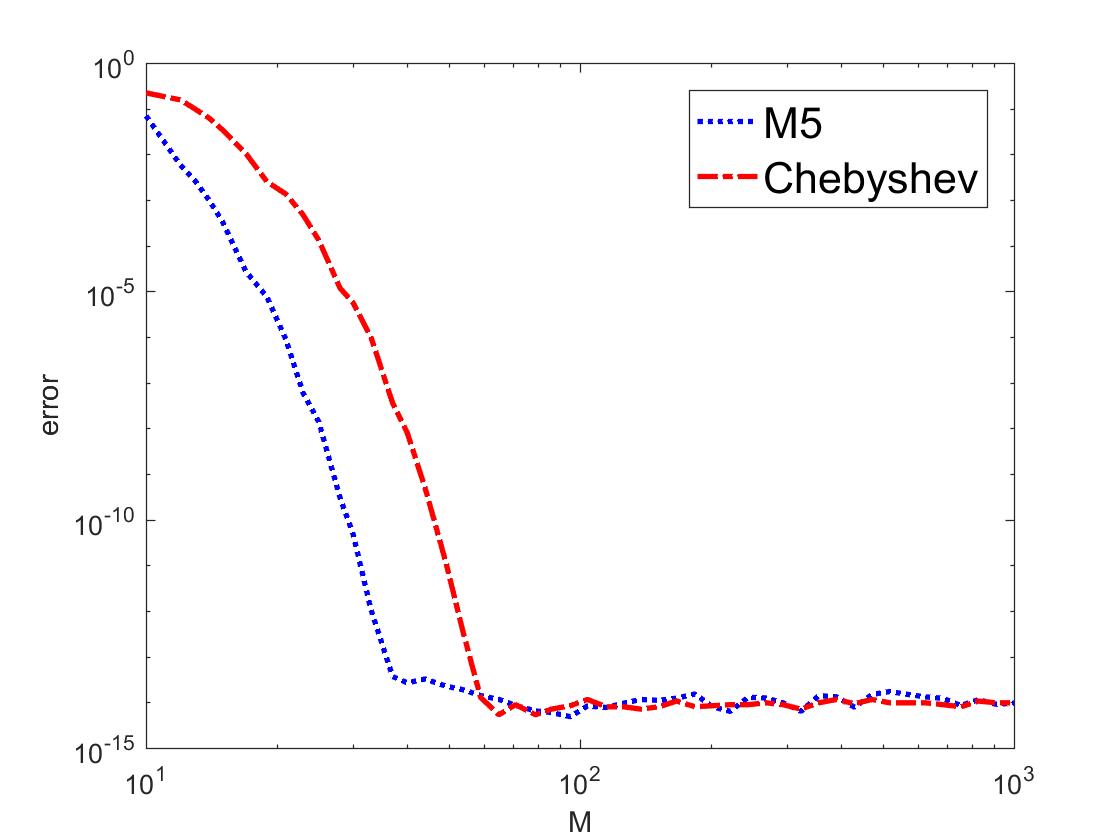}}
}%
\subfigure[\label{15b} $f_5(t)$] {
\resizebox*{6cm}{!}{\includegraphics{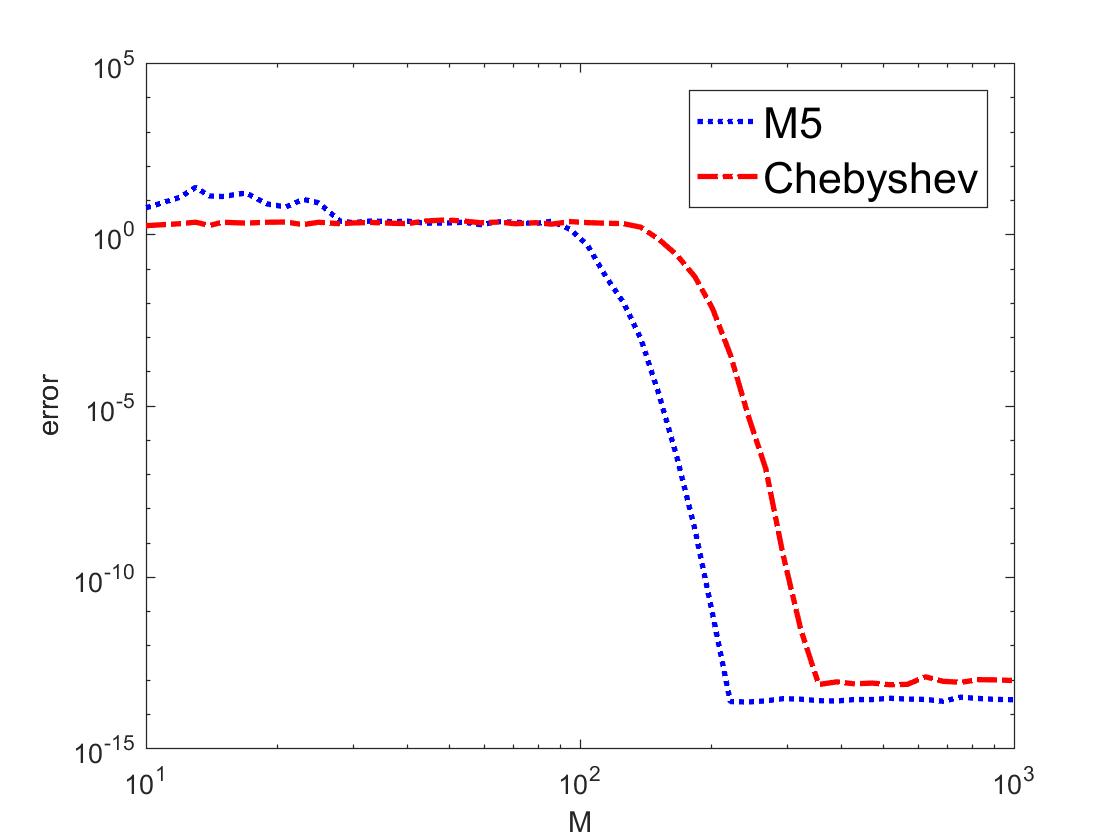}}
}%

\subfigure[\label{15c} $f_7(t)$] {
\resizebox*{6cm}{!}{\includegraphics{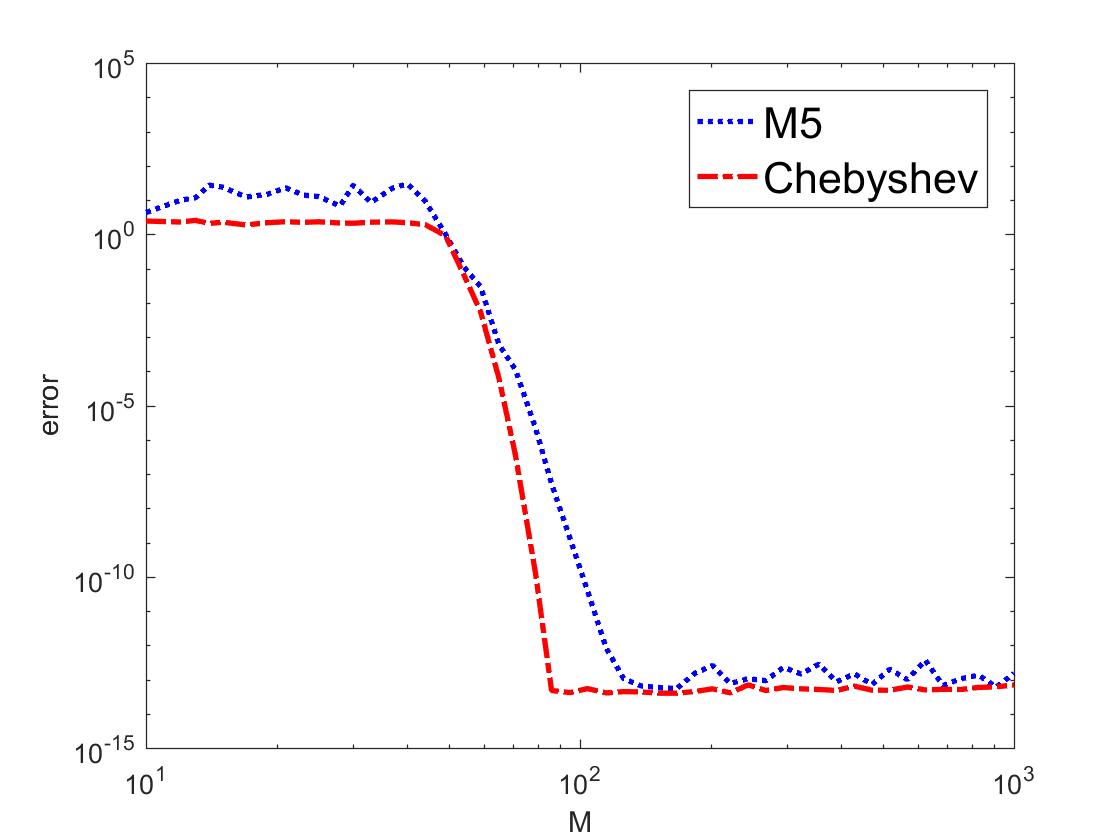}}
}%
\subfigure[\label{15d}$f_9(t)$] {
\resizebox*{6cm}{!}{\includegraphics{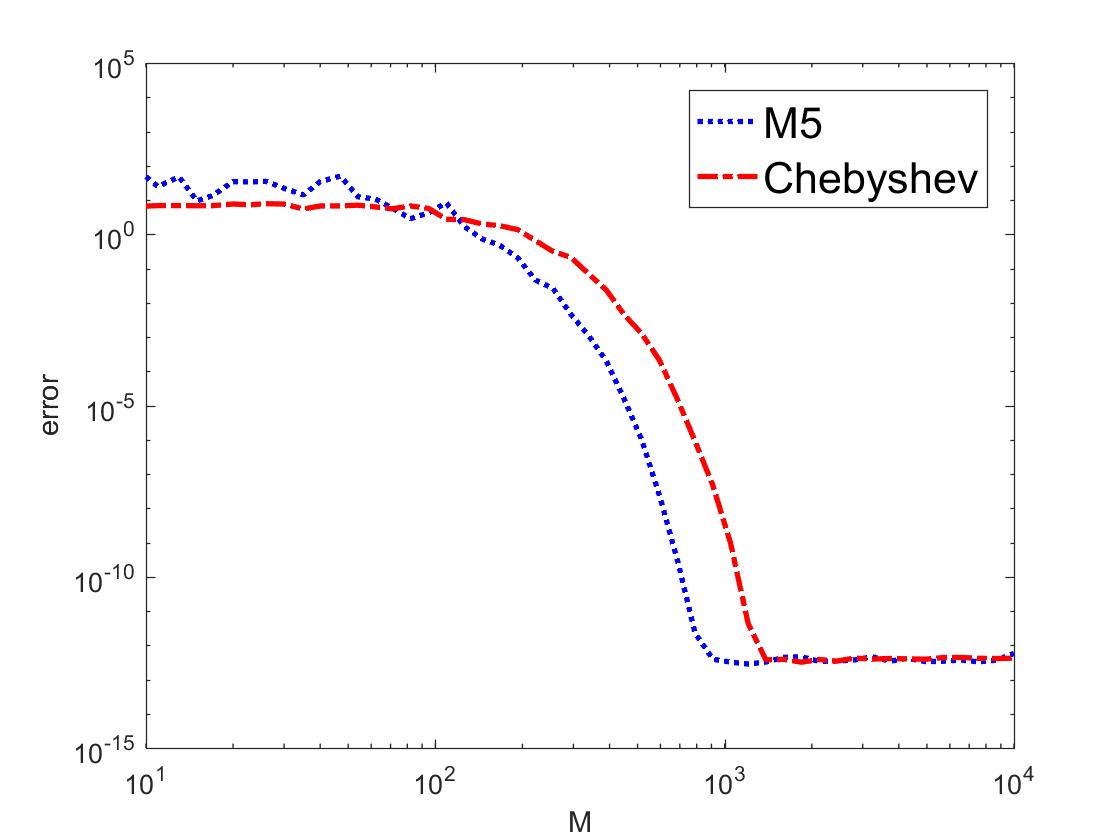}}
}%

	{\caption{ Comparisons of algorithms M5 and the Chebyshev interpolation.}\label{Figcheby}}
	\end{center}
\end{figure}}

\section{Conclusions}
Fast algorithms based on boundary interval data has been proposed for the computation of the Fourier extension approximation to nonperiodic functions. We tested the relevant parameters and provided the setting scheme. The new algorithms have  a computational complexity of $\mathcal{O}(M\log(M))$ and  can be used as basic tools for the calculation of the  Fourier extension.
\section*{Conflict of Interest declaration} The authors declare that they have no
affiliations with or involvement in any organization or entity with any
financial interest in the subject matter or materials discussed in this
manuscript.
{ \section*{Acknowledgments}The authors would like to sincerely thank Mark Lyon for generously sharing the MATLAB codes to facilitate comparative analysis in this work. We are also very grateful to the reviewers and the editor for their constructive comments.}

\end{document}